\documentclass[12pt]{article}

\usepackage{amsfonts, amssymb, amsmath, amsthm, epsf}

\makeatletter

\@addtoreset{figure}{section}

\theoremstyle{plain}
\newtheorem{theorem}{Theorem}
\@addtoreset{theorem}{section}

\newtheorem{lemma}{Lemma}
\@addtoreset{lemma}{section}

\newtheorem{corollary}{Corollary}
\@addtoreset{corollary}{section}

\newtheorem{condition}{Condition}
\@addtoreset{condition}{section}

\theoremstyle{definition}

\newtheorem{definition}{Definition}
\@addtoreset{definition}{section}

\newtheorem{example}{Example}
\@addtoreset{example}{section}

\newtheorem{remark}{Remark}
\@addtoreset{remark}{section}

\@addtoreset{equation}{section}

\makeatother

\renewcommand{\Im}{{\rm Im\,}}
\renewcommand{\Re}{{\rm Re\,}}
\newcommand{\loc}{{\rm loc}}
\newcommand{\ind}{{\rm ind\,}}
\newcommand{\supp}{{\rm supp\,}}

\renewcommand{\ker}{{\rm ker\,}}

\newcommand{\arctg}{{\rm arctg\,}}
\newcommand{\Orb}{{\rm Orb}}

\newcommand{\const}{{\rm const}}
\newcommand{\dist}{{\rm dist}}

\title{Solvability of nonlocal elliptic problems\\ in Sobolev
spaces}
\author{Pavel~Gurevich\thanks{This research was supported by
Russian Foundation for Basic Research (grant No.~03-01-06523),
Russian Ministry for Education (grant No.~E02-1.0-131), and INTAS
(grant~YSF~2002-008).}}

\date{}

\begin{document}

\maketitle
%

\begin{abstract}
We study $2m$ order elliptic equations with nonlocal
boundary-value conditions in plane angles and in bounded domains,
dealing with the case where the support of nonlocal terms
intersects the boundary. We establish necessary and sufficient
conditions under which nonlocal problems are Fredholm in Sobolev
spaces and, respectively, in weighted spaces with small weight
exponents. We also obtain an asymptotics of solutions to nonlocal
problems near the conjugation points on the boundary, where
solutions may have power singularities.
\end{abstract}

\setcounter{tocdepth}{2} \tableofcontents

\section*{Introduction}
\addcontentsline{toc}{section}{\protect\numberline{}Introduction}

Nonlocal problems have been studied since the beginning of the
20th century, but only during the last two decades these problems
have been investigated thoroughly. On the one hand, this can be
explained by significant theoretical achievements in that
direction and, on the other hand, by various applications arising
in the fields such as biophysics, theory of multidimensional
diffusion process~\cite{Feller}, plasma theory~\cite{Sam}, theory
of sandwich shells and plates~\cite{OnSk}, and so on.

In one-dimensional case, nonlocal problems were studied by
Sommerfeld~\cite{Sommerfeld}, Tamarkin~\cite{Tamarkin},
Picone~\cite{Picone}, etc. In two-dimensional case, one of the
first works is due to Carleman~\cite{Carleman}.
In~\cite{Carleman}, Carleman considered the problem of finding a
harmonic function in a plane bounded domain, satisfying the
following nonlocal condition on the boundary $\Upsilon$:
$u(y)+bu\bigl(\Omega(y)\bigr)=g(y)$, $y\in\Upsilon$, with
$\Omega:\Upsilon\to \Upsilon$ being a transformation on the
boundary such that $\Omega\bigl(\Omega(y)\bigr)\equiv y,\
y\in\Upsilon$. Such a statement of the problem originated further
investigation of nonlocal problems with transformations mapping a
boundary onto itself.

In~1969, Bitsadze and Samarskii~\cite{BitzSam} considered
essentially different kind of nonlocal problem arising in the
plasma theory: to find a function $u(y_1, y_2)$ which is harmonic
in the rectangular $G=\{y\in\mathbb R^2: -1<y_1<1,\ 0<y_2<1\}$,
continuous in $\bar G$, and satisfies the relations
 \begin{align*}
  u(y_1, 0)&=f_1(y_1),\quad u(y_1, 1)=f_2(y_1),\quad -1<y_1<1,\\
  u(-1, y_2)&=f_3(y_2),\quad u(1, y_2)=u(0, y_2),\quad 0<y_2<1,
 \end{align*}
where $f_1, f_2, f_3$ are given continuous functions. This problem
was solved in~\cite{BitzSam} by reduction to a Fredholm integral
equation and using the maximum principle. In case of arbitrary
domains and general nonlocal conditions, such a problem was
formulated as an unsolved one. Different generalizations of
nonlocal problems with transformations mapping a boundary inside
the closure of a domain were studied by Eidelman and
Zhitarashu~\cite{ZhEid}, Roitberg and Sheftel'~\cite{RSh},
Kishkis~\cite{Kishk}, Gushchin and Mikhailov~\cite{GM}, etc.

The most complete theory for $2m$ order elliptic equations with
general nonlocal conditions in multidimensional domains was
developed by Skubachevskii and his
pupils~\cite{SkMs83,SkMs86,SkDu90,SkDu91,SkBook,KovSk,SkRJMP,GurGiessen}:
classification with respect to types of nonlocal conditions was
suggested, Fredholm solvability in corresponding spaces and index
properties were studied, asymptotics of solutions near special
conjugation points was obtained. It turns out that the most
difficult situation is that where the support of nonlocal terms
intersects with the boundary. In that case, generalized solutions
to nonlocal problems may have power singularities near some points
even if the boundary and right-hand sides are infinitely
smooth~\cite{SkMs86, SkRJMP}. That is why, to investigate such
problems, weighted spaces (introduced by Kondrat'ev for
boundary-value problems in nonsmooth domains~\cite{KondrTMMO67})
are naturally applied.

\smallskip

\begin{figure}[ht]
{ \hfill\epsfbox{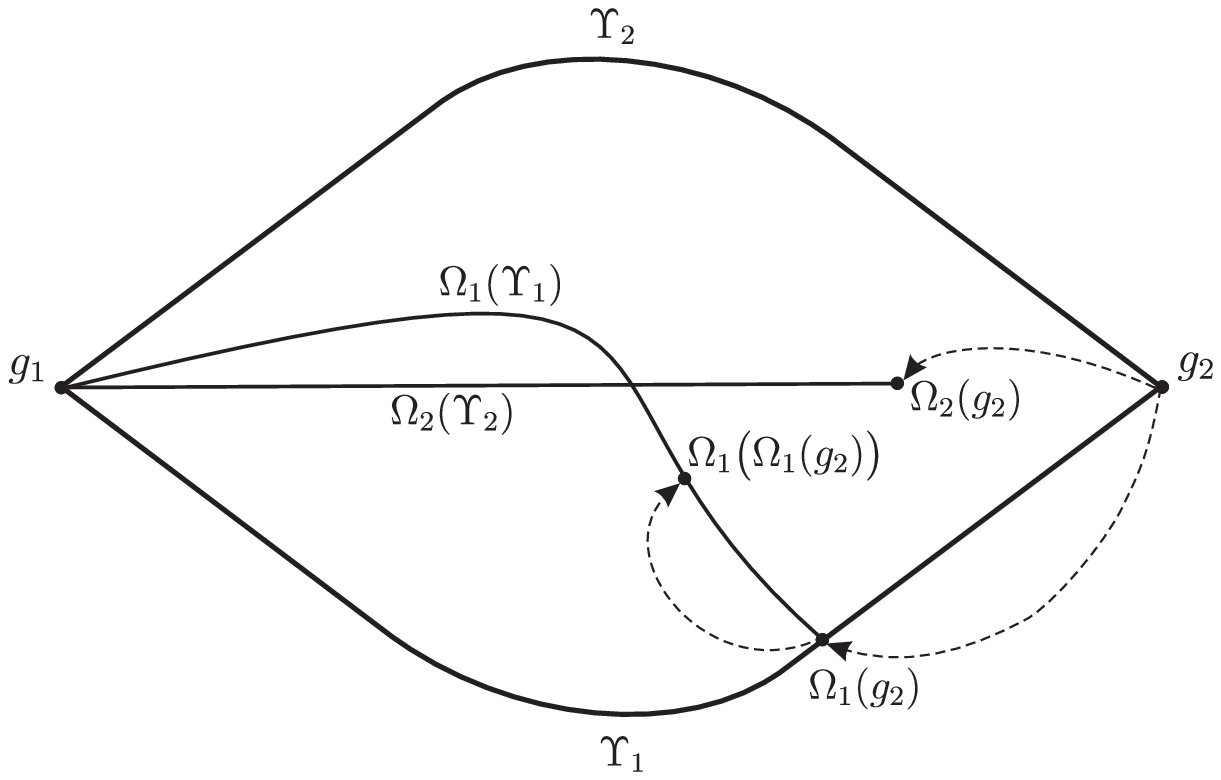}\hfill\ } \caption{Domain $G$ with the
boundary $\partial G=\bar\Upsilon_1\cup\bar\Upsilon_2$.}
   \label{figDomG}
\end{figure}

In the present paper, we study nonlocal elliptic problems in plane
domains in Sobolev spaces $W^l(G)=W_2^l(G)$ (with no weight),
dealing with the situation where the support of nonlocal terms may
intersect a boundary. Let us consider the following example. We
denote by $G\subset\mathbb R^2$ a bounded domain with boundary
$\partial G=\Upsilon_1\cup\Upsilon_2\cup\{g_1, g_2\}$, where
$\Upsilon_i$ are open (in the topology of $\partial G$)
$C^\infty$-curves, $g_1$ and $g_2$ are the end points of the
curves $\bar\Upsilon_1$, $\bar\Upsilon_2$. Let, in some
neighborhoods of $g_1$ and $g_2$, the domain $G$ coincides with
plane angles. We consider the following nonlocal problem in $G$:
\begin{align}
 \Delta u&=f_0(y)\quad (y\in G),\label{eqIntro1}\\
 u|_{\Upsilon_i}-b_i u\big(\Omega_i(y)\big)\big|_{\Upsilon_i}&=f_i(y)\quad
 (y\in\Upsilon_i;\ i=1, 2).\label{eqIntro2}
\end{align}
Here $b_1, b_2\in\mathbb R$; $\Omega_i$ is an infinitely
differentiable nondegenerate transformation mapping some
neighbourhood ${\mathcal O}_i$ of the curve $\Upsilon_i$ onto
$\Omega({\mathcal O}_i)$ so that $\Omega_i(\Upsilon_i)\subset G$
and $\overline{\omega_i(\Upsilon_i)}\cap\partial G\ne\varnothing$
(see Fig.~\ref{figDomG}). We seek for a solution $u\in W^{l+2}(G)$
under the assumption that $f_0\in W^{l}(G)$, $f_i\in
W^{l+3/2}(\Upsilon_i)$.

In this work, we will obtain necessary and sufficient condition
under which problem of type~\eqref{eqIntro1}, \eqref{eqIntro2} is
Fredholm. It will be shown that the solvability of such a problem
is influenced by {\rm(I)}~spectral properties of model nonlocal
problems with a parameter and {\rm(II)} fulfilment of some
algebraic relations between the differential operator and nonlocal
boundary-value operators at the points of conjugation of nonlocal
conditions (points $g_1$ and $g_2$ at Fig.~\ref{figDomG}). We will
consider nonlocal problems both with nonhomogeneous and with
homogeneous boundary-value conditions, which turn out to be not
equivalent ones in terms of Fredholm solvability. Near the
conjugation points, asymptotics of solutions will be obtained.

We note that nonlocal problems in Sobolev spaces in the case where
the support of nonlocal terms does not intersect the boundary was
thoroughly investigated by Skubachevskii~\cite{SkMs83, SkBook}.
However, $2m$ order elliptic equations with general nonlocal
conditions in the case where the support of nonlocal terms
intersects the boundary is being studied in Sobolev spaces for the
first time.

Our paper is organized as follows. The statement of the problem is
given in Sec.~\ref{sectStatement}. In the same section, we define
model problems in plane angles and problems with a parameter,
corresponding to the points of conjugation of nonlocal conditions.
Properties of the original problem crucially depend on whether or
not some line
\begin{equation}\label{eqLineinC}
\{\lambda\in\mathbb C: \Im\lambda=\Lambda\}
\end{equation}
(where $\Lambda\in\mathbb R$ is defined by the order of
differential equation and the order of the corresponding Sobolev
spaces) contains eigenvalues of model problems with a parameter.
In Sec.~\ref{sectLinKNoEigen}, we study nonlocal problems in plane
angles in the case where the line~\eqref{eqLineinC} contains no
eigenvalues, and in Sec.~\ref{sectLinKProperEigen} we deal with
the case where this line contains only the proper eigenvalue (see
Definition~\ref{defRegEigVal}). We use the results of
Sec.~\ref{sectLinKNoEigen} in Sec.~\ref{sectLFred} to investigate
the Fredholm solvability of the original problem in a bounded
domain, and in Sec.~\ref{sectAsymp} to obtain an asymptotics of
solutions to nonlocal problems near the conjugation points.

In~\cite{SkMs86,SkDu91,KovSk}, the authors consider nonlocal
problems in weighted spaces $H_a^l(G)$ with the norm
$$
 \|u\|_{H_a^k(G)}=
 \left(\sum\limits_{|\alpha|\le k}\int\limits_G \rho^{2(a-k+|\alpha|)}
 |D^\alpha u|^2\right)^{1/2}.
$$
Here $k\ge 0$ is an integer, $a\in\mathbb R$, and $\rho=\rho(y)$
is the distance between the point $y$ and the set of conjugation
points. For problem~\eqref{eqIntro1}, \eqref{eqIntro2}, we have
$\rho(y)=\dist(y,\{g_1,g_2\})$. In~\cite{SkDu91,KovSk}, it is
proved that if $f_0\in H_a^{l}(G)$, $f_i\in
H_a^{l+3/2}(\Upsilon_i)$, $a>l+1$, and the function $\{f_0,f_i\}$
satisfies finitely many orthogonality conditions, then
problem~\eqref{eqIntro1}, \eqref{eqIntro2} admits a solution $u\in
H_a^{l+2}(G)$. If $a\le l+1$, the following difficulty arises: the
inclusion $u\in H_a^{l+2}(G)$ does not, in general, imply that
$u\big(\Omega_i(y)\big)\big|_{\Upsilon_i}\in
H_a^{l+3/2}(\Upsilon_i)$. To eliminate this difficulty, one can
introduce the spaces (for problem~\eqref{eqIntro1},
\eqref{eqIntro2}) with the weight function
$$
\hat\rho(y)=
\dist\left(y,\big\{g_1,g_2,\Omega_1(g_2),\Omega_1\big(\Omega_1(g_2)\big),\Omega_2(g_2)\big\}\right)
$$
and arbitrary $a\in\mathbb R$ and prove, in these spaces, the
Fredholm solvability of nonlocal problems (see~\cite{SkMs86}).
However, the presence of the weight function $\hat\rho(y)$ means
that we impose a restriction both on the right-hand side and on
the solution not only near the conjugation points $g_1, g_2$ but
also near the point $\Omega_1(g_2)$ lying on a smooth part of the
boundary and near the points $\Omega_1\big(\Omega_1(g_2)\big)$ and
$\Omega_2(g_2)$ lying inside the domain (see Fig.~\ref{figDomG}).

In Sec.~\ref{sectLFredH_a}, we show: in spite of the fact that,
for $a\le l+1$, the inclusion $u\in H_a^{l+2}(G)$ does not imply
the inclusion $u\big(\Omega_i(y)\big)\big|_{\Upsilon_i}\in
H_a^{l+3/2}(\Upsilon_i)$, if $a>0$, $f_0\in H_a^{l}(G)$, $f_i\in
H_a^{l+3/2}(\Upsilon_i)$, and $\{f_0,f_i\}$ satisfies finitely
many orthogonality conditions, then problem~\eqref{eqIntro1},
\eqref{eqIntro2} yet admits a solution $u\in H_a^{l+2}(G)$. In
this case, as before, the line~\eqref{eqLineinC} (with $\Lambda$
depending now on the exponent $a$ as well) must not contain
eigenvalues of model problems with a parameter.

In Sec.~\ref{sectHatLFred}, with the help of the results from
Sec.~\ref{sectLinKProperEigen}, we study nonlocal problems in
bounded domains in the special case where the
line~\eqref{eqLineinC} contains only a proper eigenvalue of model
problems with a parameter. In this case, to provide the existence
of solutions, we impose additional consistency conditions on the
right-hand side at the conjugation points.

We note that the most complicated considerations in
sections~\ref{sectLFred}, \ref{sectLFredH_a},
and~\ref{sectHatLFred} are related to constructing right
regularizers for nonlocal problems in bounded domains. In all
these sections, we use the same scheme to construct the
regularizer, which is described in detail in Sec.~\ref{sectLFred}.
This allows us to dwell only on the most important moments in
sections~\ref{sectLFredH_a} and~\ref{sectHatLFred}.

Finally, in Sec.~\ref{secL_B}, by using the results of
sections~\ref{sectLFred} and~\ref{sectHatLFred}, we obtain a
criteria of Fredholm solvability of elliptic problems with
homogeneous nonlocal conditions. Here algebraic relations between
the differential operator and nonlocal boundary-value operators
play essential role. Two examples illustrating the results of this
paper are given in Sec.~\ref{sectEx}.

\section{Statement of Nonlocal Problems in Bounded Domains}\label{sectStatement}

\subsection{Statement of nonlocal problem}\label{subsectStatement}
Let $G\subset{\mathbb R}^2$ be a bounded domain with the boundary
$\partial G$. We introduce a set ${\mathcal K}\subset\partial G$
consisting of finitely many points and assume that $\partial
G\setminus{\mathcal K}=\bigcup\limits_{i=1}^{N_0}\Upsilon_i$,
where $\Upsilon_i$ are open (in the topology of $\partial G$)
$C^\infty$-curves. In a neighborhood of each point $g\in{\mathcal
K}$, the domain $G$ is supposed to coincide with some plane angle.

We denote by ${\bf P}(y, D_y)$ and $B_{i\mu s}(y, D_y)$
differential operators of orders $2m$ and $m_{i\mu}$ respectively
with complex-valued $C^\infty$-coefficients ($i=1, \dots, N_0;$
$\mu=1, \dots, m;$ $s=0, \dots, S_i$). Throughout the paper, we
assume that the operator ${\bf P}(y, D_y)$ is properly elliptic
for all $y\in \bar G$ and the system of operators $\{B_{i\mu 0}(y,
D_y)\}_{\mu=1}^m$ covers ${\bf P}(y, D_y)$ for all $i=1, \dots,
N_0$ and $y\in\bar\Upsilon_i$ (see, e.g.,~\cite[Ch.~2, \S~1]{LM}).

For integer $k\ge0$, we denote by $W^k(G)=W_2^k(G)$ the Sobolev
space with the norm
$$
\|u\|_{W^k(G)}=\left(\sum\limits_{|\alpha|\le k}\int\limits_G
|D^\alpha u|^2\,dy\right)^{1/2}
$$
(we put $W^0(G)=L_2(G)$ for $k=0$). For integer $k\ge1$, we
introduce the space $W^{k-1/2}(\Upsilon)$ of traces on a smooth
curve $\Upsilon\subset\bar G$, with the norm
\begin{equation}\label{eqTraceNormW}
\|\psi\|_{W^{k-1/2}(\Upsilon)}=\inf\|u\|_{W^k(G)}\quad (u\in
W^k(G): u|_\Upsilon=\psi).
\end{equation}

\smallskip

We consider the operators $\mathbf P: W^{l+2m}(G)\to W^l(G)$ and
$\mathbf B_{i\mu}^0: W^{l+2m}(G)\to
W^{l+2m-m_{i\mu}-1/2}(\Upsilon_i)$ given by ${\bf P}u={\bf P}(y,
D_y)u$ and $\mathbf B_{i\mu}^0 u=B_{i\mu 0}(y,
D_y)u(y)|_{\Upsilon_i}$. Hereinafter we assume that
$l+2m-m_{i\mu}\ge1$. The operators $\mathbf P$ and $\mathbf
B_{i\mu}^0$ will \textit{correspond to a ``local'' boundary-value
problem}.

\smallskip

Now we proceed to define the operators corresponding to nonlocal
conditions near the set $\mathcal K$. Let $\Omega_{is}$ ($i=1,
\dots, N_0;$ $s=1, \dots, S_i$) be an infinitely differentiable
nondegenerate transformation mapping some neighborhood ${\mathcal
O}_i$ of the curve $\overline{\Upsilon_i\cap\mathcal
O_{{2\varepsilon_0}}(\mathcal K)}$ onto the set
$\Omega_{is}({\mathcal O}_i)$ so that
$\Omega_{is}(\Upsilon_i)\subset G$ and
\begin{equation}\label{eqOmega}
\Omega_{is}(g)\in\mathcal K\quad\text{for }
g\in\bar\Upsilon_i\cap\mathcal K.
\end{equation}
Here ${\varepsilon_0}>0$ and $\mathcal
O_{{2\varepsilon_0}}(\mathcal K)=\{y\in \mathbb R^2: \dist(y,
\mathcal K)<{2\varepsilon_0}\}$ is the
$2{\varepsilon_0}$-neighborhood of the set $\mathcal K$. Thus,
under the transformations $\Omega_{is}$, the curves $\Upsilon_i$
map strictly inside the domain $G$ while the set of end points of
$\Upsilon_i$ maps to itself.

Let $\varepsilon_0$ be so small (see Remark~\ref{remSmallEps}
below) that, in the $2{\varepsilon_0}$-neighborhood $\mathcal
O_{2{\varepsilon_0}}(g)$ of each point $g\in\mathcal K$, the
domain $G$ coincides with a plane angle. Let us specify the
structure of the transformation $\Omega_{is}$ near the set
$\mathcal K$.

We denote by $\Omega_{is}^{+1}$ the transformation
$\Omega_{is}:{\mathcal O}_i\to\Omega_{is}({\mathcal O}_i)$ and by
$\Omega_{is}^{-1}$ the transformation
$\Omega_{is}^{-1}:\Omega_{is}({\mathcal O}_i)\to{\mathcal O}_i$
being inverse to $\Omega_{is}.$ The set of all points
$\Omega_{i_qs_q}^{\pm1}(\dots\Omega_{i_1s_1}^{\pm1}(g))\in{\mathcal
K}$ ($1\le s_j\le S_{i_j},\ j=1, \dots, q$) (i.e., points which
can be obtained by consecutive applying to the point $g$ the
transformations $\Omega_{i_js_j}^{+1}$ or $\Omega_{i_js_j}^{-1}$
taking the points of ${\cal K}$ to those of ${\cal K}$) is called
an {\it orbit} of $g\in{\cal K}$ and denoted by $\Orb(g)$.

Clearly, for any $g, g'\in{\mathcal K}$, either $\Orb(g)=\Orb(g')$
or $\Orb(g)\cap\Orb(g')=\varnothing$. Thus, we have ${\mathcal
K}=\bigcup\limits_{p=1}^{N_1}\Orb_p$, where
$\Orb_{p_1}\cap\Orb_{p_2}=\varnothing$ ($p_1\ne p_2$), and, for
each $p=1, \dots, N_1$, the set $\Orb_p$ coincides with an orbit
of some point $g\in {\mathcal K}$. Let each orbit $\Orb_p$ consist
of points $g^{p}_j$, $j=1, \dots, N_{1p}$.

For every point $g\in\mathcal K$, we consider neighborhoods
\begin{equation}\label{eqVepsilon}
\hat{\mathcal V}(g)\supset{\mathcal V}(g)\supset\mathcal
O_{2\varepsilon_0}(g)
\end{equation}
such that
\begin{enumerate}
\item[(1)] in the neighborhood $\hat{\mathcal V}(g)$, the boundary
$\partial G$ coincides with a plane angle;
\item[(2)]
$\overline{\hat{\mathcal V}(g)}\cap\overline{\hat{\mathcal
V}(g')}=\varnothing$ for any $g,g'\in\mathcal K$, $g\ne g'$;
\item[(3)] if $g^{p}_j\in\bar\Upsilon_i\cap\Orb_p$ and
$\Omega_{is}(g^{p}_j)=g^{p}_k,$ then ${\mathcal
V}(g^{p}_j)\subset\mathcal
 O_i$ and
 $\Omega_{is}\big({\mathcal V}(g^{p}_j)\big)\subset\hat{\mathcal V}(g^{p}_k).$
\end{enumerate}

For each $g^{p}_j\in\bar\Upsilon_i\cap\Orb_p$, we fix an argument
transformation $y\mapsto y'(g^{p}_j)$ which is a composition of
the shift by the vector $-\overrightarrow{Og^{p}_j}$ and rotation
by some angle so that the set ${\mathcal V}(g^{p}_j)$
($\hat{\mathcal V}(g^{p}_j)$) maps onto a neighborhood ${\mathcal
V}^{p}_j(0)$ ($\hat{\mathcal V}^{p}_j(0)$) of the origin while the
sets $G\cap{\mathcal V}(g^{p}_j)$ ($G\cap\hat{\mathcal
V}(g^{p}_j)$) and $\Upsilon_i\cap{\mathcal V}(g^{p}_j)$
($\Upsilon_i\cap\hat{\mathcal V}(g^{p}_j)$) map to the
intersection of the plane angle $K^{p}_j=\{y\in{\mathbb R}^2:\
r>0,\ |\omega|<b^{p}_j<\pi\}$ with ${\mathcal V}^{p}_j(0)$
 ($\hat{\mathcal V}^{p}_j(0)$) and the intersection of
a side of the angle $K^{p}_j$ with
 ${\mathcal V}^{p}_j(0)$ ($\hat{\mathcal V}^{p}_j(0)$) respectively.

\begin{condition}\label{condK1}
The argument transformation $y\mapsto y'(g)$ for $y\in{\mathcal
V}(g)$, $g\in\mathcal K\cap\bar\Upsilon_i$, described above
reduces the transformation $\Omega_{is}(y)$ ($i=1,\dots N_0$,
$s=1,\dots, S_i$) to a composition of rotation and expansion in
new variables $y'$.
\end{condition}

\begin{remark}\label{remK1}
Condition~\ref{condK1} combined with the assumption that
$\Omega_{is}(\Upsilon_i)\subset G$, in particular, means that {\em
if
$g\in\Omega_{is}(\bar\Upsilon_i\setminus\Upsilon_i)\cap\bar\Upsilon_j\cap{\mathcal
K}\ne\varnothing$, then the curves $\Omega_{is}(\bar\Upsilon_i)$
and $\bar\Upsilon_j$ are nontangent to each other at the point
$g$.}
\end{remark}

We introduce the bounded operators $\mathbf B_{i\mu}^1:
W^{l+2m}(G)\to W^{l+2m-m_{i\mu}-1/2}(\Upsilon_i)$ by the formula
$$
 \mathbf B_{i\mu}^1u=\sum\limits_{s=1}^{S_i}
   \big(B_{i\mu s}(y,
   D_y)(\zeta u)\big)\big(\Omega_{is}(y)\big)\big|_{\Upsilon_i}.
$$
Here $\big(B_{i\mu s}(y,
D_y)v\big)\big(\Omega_{is}(y)\big)=B_{i\mu s}(y',
D_{y'})v(y')|_{y'=\Omega_{is}(y)}$ and the function $\zeta\in
C^\infty(\mathbb R^2)$ is such that
\begin{equation}\label{eqZeta}
 \zeta(y)=1\ (y\in\mathcal O_{{\varepsilon_0}/2}(\mathcal K)),\quad
 \zeta(y)=0\ (y\notin\mathcal O_{{\varepsilon_0}}(\mathcal K)).
\end{equation}
Since $\mathbf B_{i\mu}^1u=0$ whenever $\supp u\subset\bar
G\setminus\overline{\mathcal O_{{\varepsilon_0}}(\mathcal K)}$, we
say that the operator $\mathbf B_{i\mu}^1$ \textit{corresponds to
nonlocal terms with the support near the set} $\mathcal K$.

\smallskip

We also introduce the bounded operator $\mathbf B_{i\mu}^2:
W^{l+2m}(G)\to W^{l+2m-m_{i\mu}-1/2}(\Upsilon_i)$ satisfying the
following condition.
\begin{condition}\label{condSeparK23}
There exist numbers $\varkappa_1>\varkappa_2>0$ and $\rho>0$ such
that, for all $u\in W^{l+2m}(G\setminus\overline{\mathcal
O_{\varkappa_1}(\mathcal K)})\cup W^{l+2m}(G_\rho)$, the following
inequalities hold:
\begin{equation}\label{eqSeparK23'}
  \|\mathbf B^2_{i\mu}u\|_{W^{l+2m-m_{i\mu}-1/2}(\Upsilon_i)}\le c_1
  \|u\|_{W^{l+2m}(G\setminus\overline{\mathcal O_{\varkappa_1}(\mathcal
  K)})},
\end{equation}
\begin{equation}\label{eqSeparK23''}
  \|\mathbf B^2_{i\mu}u\|_{W^{l+2m-m_{i\mu}-1/2}
   (\Upsilon_i\setminus\overline{\mathcal O_{\varkappa_2}(\mathcal K)})}\le
  c_2 \|u\|_{W^{l+2m}(G_\rho)},
\end{equation}
where $i=1, \dots, N_0$; $\mu=1, \dots, m$; $c_1,c_2>0$;
$G_\rho=\{y\in G: \dist(y, \partial G)>\rho\}$.
\end{condition}

From~\eqref{eqSeparK23'}, it follows that $\mathbf B_{i\mu}^2u=0$
whenever $\supp u\subset \mathcal O_{\varkappa_1}(\mathcal K)$.
Therefore, we say that the operator $\mathbf B_{i\mu}^2$
\textit{corresponds to nonlocal terms with the support outside the
set} $\mathcal K$.

\smallskip

We suppose that Conditions~\ref{condK1} and~\ref{condSeparK23} are
fulfilled throughout.

\smallskip

Notice that we a priori assume no connection between the numbers
$\varkappa_1,\varkappa_2,\rho$ in Condition~\ref{condSeparK23} and
the number $\varepsilon_0$ in Condition~\ref{condK1}.

We study the following nonlocal elliptic problem:
\begin{align}
 {\bf P} u&=f_0(y) \quad (y\in G),\label{eqPinG}\\
     \mathbf B_{i\mu}^0 u+\mathbf B_{i\mu}^1 u+\mathbf B_{i\mu}^2 u&=
   f_{i\mu}(y)\quad
    (y\in \Upsilon_i;\ i=1, \dots, N_0;\ \mu=1, \dots, m).\label{eqBinG}
\end{align}

Let us introduce the following operator corresponding to
problem~(\ref{eqPinG}), (\ref{eqBinG}):
$$
{\bf L}=\{{\bf P},\ {\bf B}_{i\mu}^0+{\bf B}_{i\mu}^1+{\bf
B}_{i\mu}^2 \}: W^{l+2m}(G)\to  {\mathcal W}^l(G,\Upsilon),
$$
where ${\mathcal
W}^l(G,\Upsilon)=W^l(G)\times\prod\limits_{i=1}^{N_0}\prod\limits_{\mu=1}^m
W^{l+2m-m_{i\mu}-1/2}(\Upsilon_i)$.

\begin{remark}\label{remSmallEps}
Further, we will need that $\varepsilon_0$ be sufficiently small
(while $\varkappa_1,\varkappa_2,\rho$ may be arbitrary). Let us
show that this does not lead to the loss of generality.

Let us have a number $\hat\varepsilon_0$ such that
$0<\hat\varepsilon_0<\varepsilon_0$. We consider a function
$\hat\zeta\in C^\infty(\mathbb R^2)$ satisfying
$$
 \hat\zeta(y)=1\ (y\in\mathcal O_{{\hat\varepsilon_0}/2}(\mathcal K)),\quad
 \hat\zeta(y)=0\ (y\notin\mathcal O_{{\hat\varepsilon_0}}(\mathcal K))
$$
and introduce the operator $\hat{\mathbf B}_{i\mu}^1:
W^{l+2m}(G)\to W^{l+2m-m_{i\mu}-1/2}(\Upsilon_i)$ by the formula
$$
 \hat{\mathbf B}_{i\mu}^1u=\sum\limits_{s=1}^{S_i}
   \big(B_{i\mu s}(y,
   D_y)(\hat\zeta
   u)\big)\big(\Omega_{is}(y)\big)\big|_{\Upsilon_i}.
$$
Clearly, we have
$$
\mathbf B_{i\mu}^0+\mathbf B_{i\mu}^1+\mathbf B_{i\mu}^2= \mathbf
B_{i\mu}^0+\hat{\mathbf B}_{i\mu}^1+\hat{\mathbf B}_{i\mu}^2,
$$
where $\hat{\mathbf B}_{i\mu}^2=\mathbf B_{i\mu}^1-\hat{\mathbf
B}_{i\mu}^1+\mathbf B_{i\mu}^2$. From
example~\ref{exGeneralProblem} (see Sec.~\ref{subsecExinG}), it
follows that the operator $\mathbf B_{i\mu}^1-\hat{\mathbf
B}_{i\mu}^1$ satisfies Condition~\ref{condSeparK23} for some
$\varkappa_1,\varkappa_2,\rho$. Therefore, we can always choose
$\varepsilon_0$ being as small as necessary (maybe at the expense
of the change of the operator $\mathbf B_{i\mu}^2$ and values of
$\varkappa_1,\varkappa_2,\rho$).
\end{remark}

\subsection{Example of nonlocal problem}\label{subsecExinG}
In the following example, we present a concrete realization for
the abstract nonlocal operators $\mathbf B_{i\mu}^2$.
\begin{example}\label{exGeneralProblem}
Let the operators ${\bf P}(y, D_y)$ and $B_{i\mu s}(y, D_y)$ be
the same as before. Let $\Omega_{is}$ ($i=1,\ \dots,\ N_0;$ $s=1,\
\dots,\ S_i$) be an infinitely differentiable nondegenerate
transformation mapping some neighborhood ${\mathcal O}_i$ of the
curve $\Upsilon_i$ onto $\Omega_{is}({\mathcal O}_i)$ so that
$\Omega_{is}(\Upsilon_i)\subset G$. Notice that in this example
\textit{assumption~\eqref{eqOmega} is not necessarily supposed to
hold for each $\Omega_{is}$.}

We consider the following nonlocal problem:
\begin{equation}\label{eqPinGEx}
 {\bf P}(y, D_y) u=f_0(y) \quad (y\in G),
\end{equation}
\begin{equation}\label{eqBinGEx}
 \begin{aligned}
 B_{i\mu 0}(y, D_y)u(y)|_{\Upsilon_i}+
 \sum\limits_{s=1}^{S_i}
   \big(B_{i\mu s}(y, D_y)u\big)\big(\Omega_{is}(y)\big)\big|_{\Upsilon_i}=
   f_{i\mu}(y)\\
    (y\in \Upsilon_i;\ i=1, \dots, N_0;\ \mu=1, \dots, m).
 \end{aligned}
\end{equation}
Let us choose ${\varepsilon_0}$ so small that, for any point
$g\in\mathcal K$, the set $\overline{\mathcal
O_{{\varepsilon_0}}(g)}$ intersects with the curve
$\overline{\Omega_{is}(\Upsilon_i)}$ only if $g\in\mathcal
K\cap\overline{\Omega_{is}(\Upsilon_i)}$.

Let a point $g\in\mathcal K\cap\bar\Upsilon_i$ be such that
$\Omega_{is}(g)\in\mathcal K$. Then we define the orbit $\Orb(g)$
of the point $g$ analogously to the above and assume that, for
each point of this orbit $\Orb(g)$, Condition~\ref{condK1} holds.

\begin{remark}
According to Remark~\ref{remK1}, Condition~\ref{condK1} is a
restriction upon a geometrical structure of the support of
nonlocal terms near the set $\mathcal K$. However, if
$\Omega_{is}(\bar\Upsilon_i\setminus\Upsilon_i)\subset\partial
G\setminus{\mathcal K}$, then we impose no restrictions upon a
geometrical structure of the curve $\Omega_{is}(\bar\Upsilon_i)$
near $\partial G$ (cf.~\cite{SkMs86, SkDu91}).
\end{remark}

We put
\begin{align*}
{\bf P}u&={\bf P}(y, D_y)u\\
 \mathbf B_{i\mu}^0 u&=B_{i\mu 0}(y,
D_y)u(y)|_{\Upsilon_i},\\
\mathbf B_{i\mu}^1u&=\sum\limits_{s=1}^{S_i} \big(B_{i\mu s}(y,
D_y)(\zeta u)\big)\big(\Omega_{is}(y)\big)\big|_{\Upsilon_i},\\
 \mathbf B_{i\mu}^2u&=\sum\limits_{s=1}^{S_i}
   \big(B_{i\mu s}(y,
   D_y)((1-\zeta)
   u)\big)\big(\Omega_{is}(y)\big)\big|_{\Upsilon_i},
\end{align*}
where $\zeta$ is defined by~\eqref{eqZeta} (see
figures~\ref{figB1} and~\ref{figB2}). Then
problem~(\ref{eqPinGEx}), (\ref{eqBinGEx}) assumes the
form~(\ref{eqPinG}), (\ref{eqBinG}).
\begin{figure}[p]
{ \hfill\epsfbox{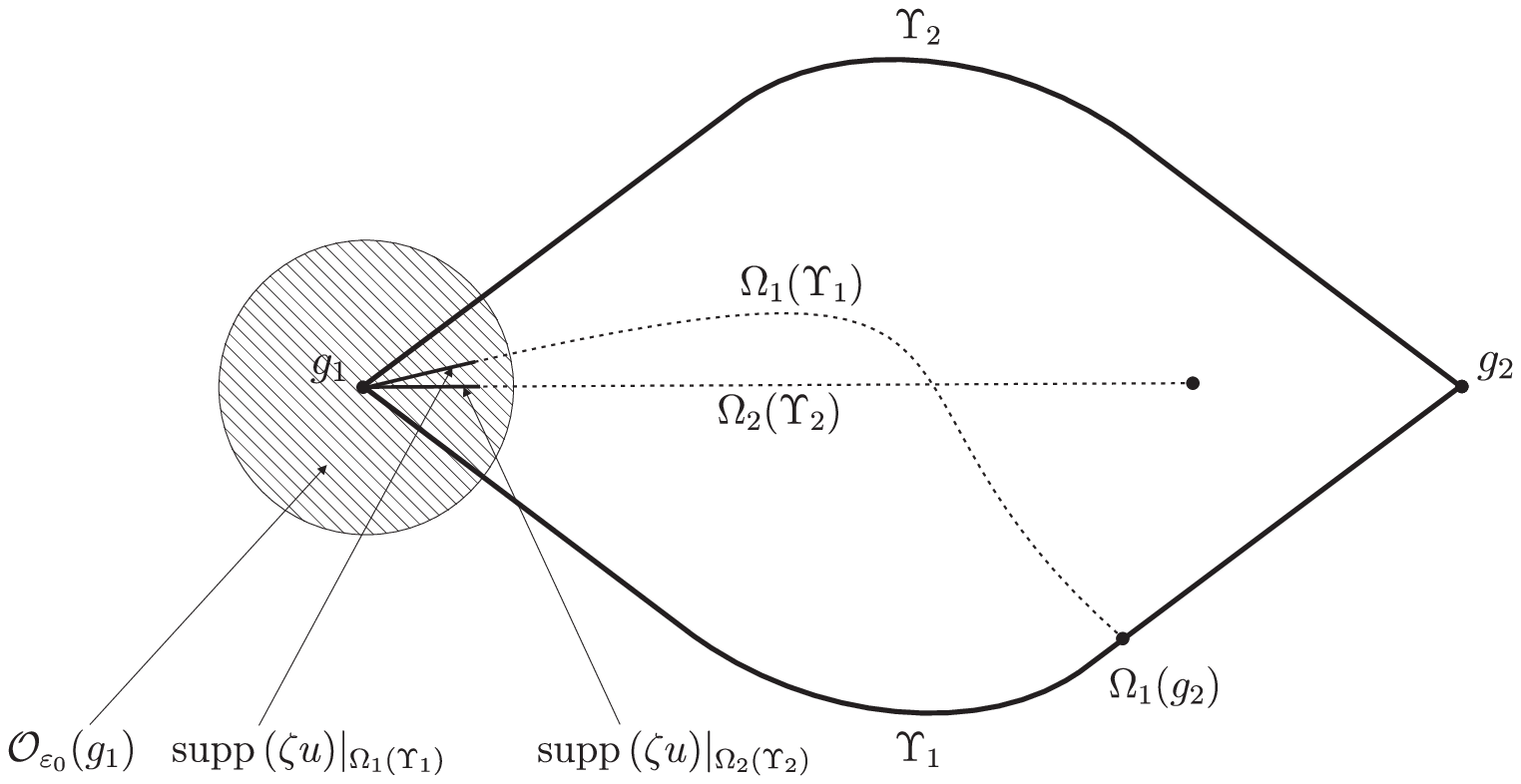}\hfill\ } \caption{Dotted lines denote the
support of nonlocal terms corresponding to the operator~$\mathbf
B_{i\mu}^2$.}
   \label{figB1}
\end{figure}
\begin{figure}[p]
{ \hfill\epsfbox{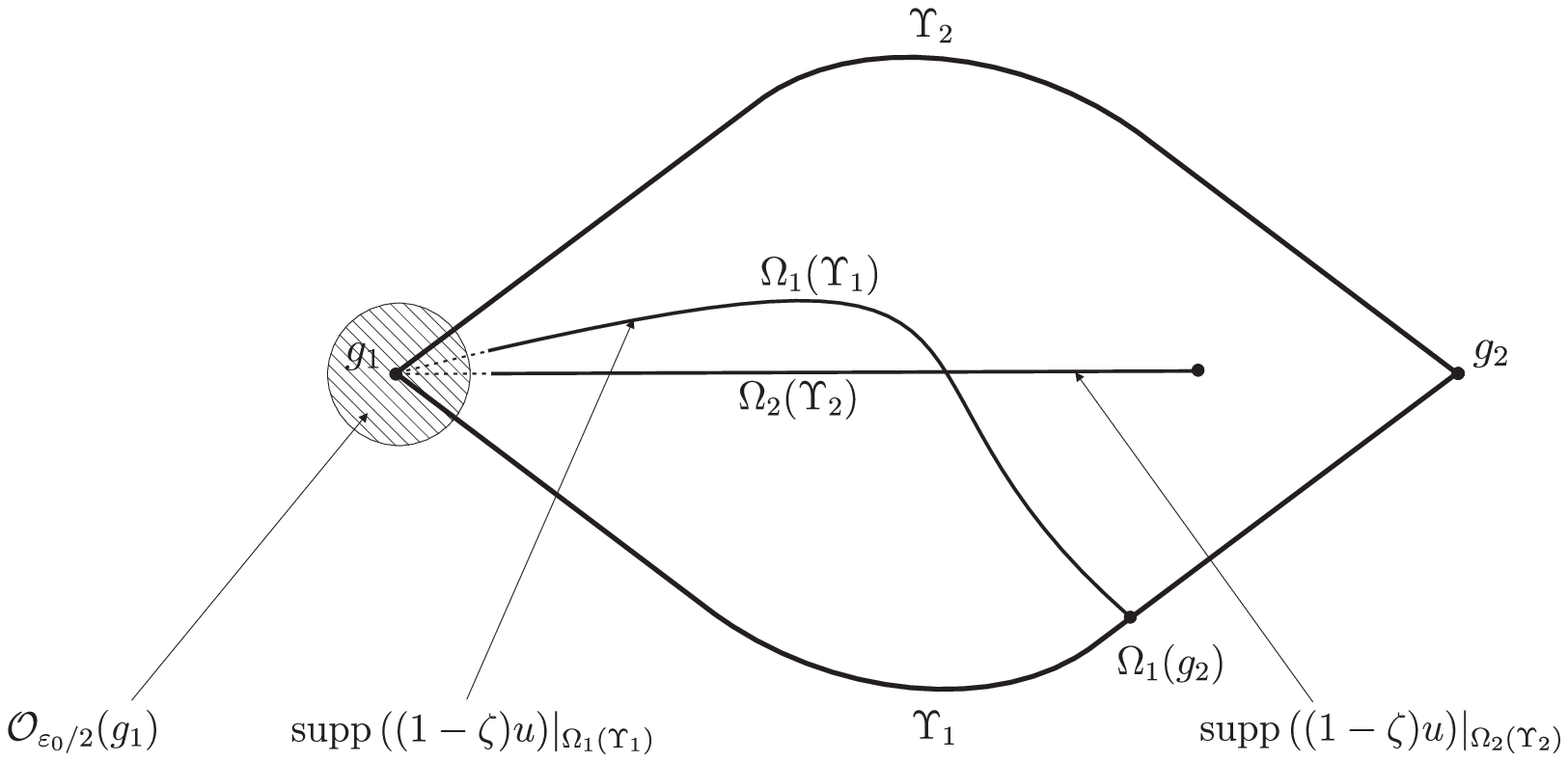}\hfill\ } \caption{Dotted lines denote the
support of nonlocal terms corresponding to the operator~$\mathbf
B_{i\mu}^1$.}
   \label{figB2}
\end{figure}

Analogously to the proof of Lemma~2.5~\cite{SkDu91} (where
weighted spaces should be replaced by corresponding Sobolev
spaces), one can show that the operator $\mathbf B_{i\mu}^2$
satisfies Condition~\ref{condSeparK23}. Let us prove, for example,
inequality~\eqref{eqSeparK23'}. Clearly, it suffices to consider
an arbitrary term $\psi=\big(B_{i\mu s}(y,
   D_y)((1-\zeta)
   u)\big)\big(\Omega_{is}(y)\big)\big|_{\Upsilon_i}$.
We introduce a function $v\in C_0^\infty\big(\Omega_{is}({\mathcal
O}_i)\big)$ such that
\begin{gather}
v|_{\Omega_{is}(\Upsilon_i)}=\big(B_{i\mu s}(y,
D_y)((1-\zeta)u)\big)\big|_{\Omega_{is}(\Upsilon_i)},\label{eqContinuation1}\\
\|v\|_{W^{l+2m-m_{i\mu}}(\Omega_{is}({\mathcal O}_i))}\le
2\big\|\big(B_{i\mu s}(y, D_y)
((1-\zeta)u)\big)\big|_{\Omega_{is}(\Upsilon_i)}\big\|_{W^{l+2m-m_{i\mu}-1/2}(\Omega_{is}(\Upsilon_i))},\label{eqContinuation2}
\end{gather}
From~\eqref{eqContinuation1}, it follows that
$$
 v\big(\Omega_{is}(y)\big)\big|_{\Upsilon_i}=\psi.
$$
Combining this with the boundedness of the trace operator in
Sobolev spaces and inequality~\eqref{eqContinuation2}, we get
\begin{multline}\label{eqproveSeparK23'}
\|\psi\|_{W^{l+2m-m_{i\mu}-1/2}(\Upsilon_i)}=
\big\|v\big(\Omega_{is}(y)\big)\big|_{\Upsilon_i}\big\|_{W^{l+2m-m_{i\mu}-1/2}(\Upsilon_i)}
\le\left\|v\big(\Omega_{is}(y)\big)\right\|_{W^{l+2m-m_{i\mu}}(\mathcal
O_i)}\\ \le k_1 \|v\|_{W^{l+2m-m_{i\mu}}(\Omega_{is}(\mathcal
O_i))} \le 2k_1\|B_{i\mu s}(y,
D_y)((1-\zeta)u)|_{\Upsilon_i}\|_{W^{l+2m-m_{i\mu}-1/2}(\Omega_{is}(\Upsilon_i))}\\
\le k_2\|(1-\zeta)u\|_{W^{l+2m}(G)}.
\end{multline}
Thus, putting $\varkappa_1={\varepsilon_0}/2$, we see
that~\eqref{eqproveSeparK23'} implies
estimate~\eqref{eqSeparK23'}. Notice that, in this case, the
numbers $\varkappa_1$ and ${\varepsilon_0}$ turn out to be
connected with each other.

Analogous considerations allow one to obtain
estimate~\eqref{eqSeparK23''}. The proof is based on the
boundedness of the trace operator, smoothness of the
transformations $\Omega_{is}$, and relation
$$
 \Omega_{is}(\Upsilon_i\setminus\overline{\mathcal O_{\varkappa_2}(\mathcal
 K)})\subset G_\rho
$$
(which is valid for any $\varkappa_2<\varkappa_1$ and sufficiently
small $\rho=\rho(\varkappa_2)$). The latter relation follows from
the embedding $\Omega_{is}(\Upsilon_i)\subset G$ and continuity of
$\Omega_{is}$.
\end{example}

\subsection{Nonlocal problems near the set $\mathcal
K$}\label{subsectStatementNearK}

While studying problem~(\ref{eqPinG}), (\ref{eqBinG}), one must
pay especial attention to a behavior of solutions in a
neighborhood of the set ${\mathcal K}$, which consists of the
conjugation points. Let us consider corresponding model problems
in plane angles. To this end, we formally assume that
\begin{equation}\label{eqB2=0}
 {\mathbf B}_{i\mu}^2=0,\quad i=1,\dots,N_0,\ \mu=1,\dots,m.
\end{equation}

Let us fix some orbit $\Orb_p\subset{\mathcal K}$ ($p=1, \dots,
N_1$) and suppose that $\supp
u\subset\Big(\bigcup\limits_{j=1}^{N_{1p}}{\mathcal
V}(g^{p}_j)\Big)\cap\bar G.$ We denote by $u_j(y)$ the function
$u(y)$ for $y\in\hat{\mathcal V}(g^{p}_j)\cap G$. If
$g^{p}_j\in\bar\Upsilon_i,$ $y\in{\mathcal V}(g^{p}_j)$, and
$\Omega_{is}(y)\in\hat{\mathcal V}(g^{p}_k),$ we denote
$u(\Omega_{is}(y))$ by $u_k(\Omega_{is}(y))$. Then, by virtue of
assumption~\eqref{eqB2=0}, nonlocal problem~(\ref{eqPinG}),
(\ref{eqBinG}) assumes the following form:
\begin{gather*}
 {\bf P}(y, D_y) u_j=f_0(y) \quad (y\in{\mathcal V}(g^{p}_j)\cap
 G),\\
\begin{aligned}
B_{i\mu 0}(y, D_y)u_j(y)|_{{\mathcal V}(g^{p}_j)\cap\Upsilon_i}+
 \sum\limits_{s=1}^{S_i}
   \big(B_{i\mu s}(y,
   D_y)(\zeta u_k)\big)\big(\Omega_{is}(y)\big)\big|_{{\mathcal V}(g^{p}_j)\cap\Upsilon_i}=f_{i\mu}(y) \\
   \big(y\in {\mathcal V}(g^{p}_j)\cap\Upsilon_i;\ i\in\{1\le i\le N_0:
   g^{p}_j\in\bar\Upsilon_i\};\
  j=1, \dots, N_{1p};\ \mu=1, \dots, m\big).
\end{aligned}
\end{gather*}

Let $y\mapsto y'(g^{p}_j)$ be the argument transformation
described above. We introduce the function $U_j(y')=u_j(y(y'))$
and denote $y'$ again by $y$. For $p$ being fixed, we put
$N=N_{1p}$, $b_j=b^{p}_j$, $K_j=K^p_j$ (see
Sec.~\ref{subsectStatement}),  and
$\gamma_{j\sigma}=\{y\in{\mathbb R}^2: r>0,\ \omega=(-1)^\sigma
b_{j}\}$ ($\sigma=1, 2$), where $(\omega, r)$ are polar
coordinates with pole at the origin. Now, using
Condition~\ref{condK1}, we can rewrite problem~(\ref{eqPinG}),
(\ref{eqBinG}) as follows:
\begin{gather}
  {\bf P}_{j}(y, D_y)U_j=f_{j}(y) \quad (y\in
  K_j),\label{eqPinK}\\
  {\bf B}_{j\sigma\mu}(y, D_y)U|_{\gamma_{j\sigma}}\equiv \sum\limits_{k,s}
       (B_{j\sigma\mu ks}(y, D_y)U_k)({\mathcal G}_{j\sigma ks}y)|_{\gamma_{j\sigma}}
    =f_{j\sigma\mu}(y) \quad (y\in\gamma_{j\sigma}).\label{eqBinK}
\end{gather}
Here (and further, unless the contrary is specified) $j, k=1,
\dots, N=N_{1p};$ $\sigma=1, 2;$ $\mu=1, \dots, m;$ $s=0, \dots,
S_{j\sigma k}$; ${\mathbf P}_j(y, D_y)$ and $B_{j\sigma\mu ks}(y,
D_y)$ are operators of orders $2m$ and $m_{j\sigma\mu}$
respectively with variable $C^\infty$-coefficients; ${\mathcal
G}_{j\sigma ks}$ is the operator of rotation by an
angle~$\omega_{j\sigma ks}$ and expansion $\chi_{j\sigma ks}$
($\chi_{j\sigma ks}>0$) times in $y$-plane. Furthermore,
$|(-1)^\sigma b_{j}+\omega_{j\sigma ks}|<b_{k}$ for $(j, 0)\ne(k,
s)$, $\omega_{j\sigma j0}=0$, and $\chi_{j\sigma j0}=1$.

Since $\mathcal V(0)\supset\mathcal O_{\varepsilon_0}(0)$
(see.~\eqref{eqVepsilon}), it follows that, for any function $v$
(which need not be compactly supported), we have
\begin{equation}\label{eqCoefficientsB}
B_{j\sigma\mu ks}(y, D_y)v(y)=0\quad \text{for
 } |y|\ge\varepsilon_0,\ (k,s)\ne(j,0).
\end{equation}
Moreover, since we consider problem~(\ref{eqPinK}), (\ref{eqBinK})
for functions $U$ with compact support, we may assume that the
coefficients of the operators ${\bf P}_{j}(y, D_y)$ and
$B_{j\sigma\mu j0}(y, D_y)$ are equal to zero outside a disk of
sufficiently large radius.

Let us introduce the following spaces of vector-functions:
\begin{gather*}
 W^{l+2m,N}(K)=\prod_{j=1}^{N} W^{l+2m}(K_j),\
 \mathcal W^{l,N}(K,\gamma)=\prod_{j=1}^{N} {\mathcal
 W}^l(K_j,\gamma_j),\\
 {\mathcal W}^l(K_j,\gamma_j)=
 W^l(K_j)\times\prod\limits_{\sigma=1,2}\prod\limits_{\mu=1}^m W^{l+2m-m_{j\sigma\mu}-1/2}
   (\gamma_{j\sigma}).
\end{gather*}

We consider the operator ${\bf L}_p: W^{l+2m,N}(K)\to \mathcal
W^{l,N}(K,\gamma)$ given by
$$
 {\bf L}_p U=\{{\bf P}_j(y, D_y)U_j,\ {\bf B}_{j\sigma\mu}(y, D_y)U|_{\gamma_{j\sigma}}\}
$$
and corresponding to problem~(\ref{eqPinK}), (\ref{eqBinK}).
Subindex $p$ means that the operator ${\bf L}_p$ is related to the
orbit $\Orb_p$.

We denote by ${\mathcal P}_j(D_y)$ and $B_{j\sigma\mu ks}(D_y)$
the principal homogeneous parts of the operators ${\bf P}_j(0,
D_y)$ and $B_{j\sigma\mu ks}(0, D_y)$ respectively. Along with
problem~(\ref{eqPinK}), (\ref{eqBinK}), we study the model
nonlocal problem
\begin{equation}\label{eqPinK0}
  {\mathcal P}_{j}(D_y)U_j=f_{j}(y) \quad (y\in K_j),
\end{equation}
\begin{equation}\label{eqBinK0}
  {\mathcal B}_{j\sigma\mu}(D_y)U|_{\gamma_{j\sigma}}\equiv \sum\limits_{k,s}
       (B_{j\sigma\mu ks}(D_y)U_k)({\mathcal G}_{j\sigma ks}y)|_{\gamma_{j\sigma}}
    =f_{j\sigma\mu}(y) \quad (y\in\gamma_{j\sigma}).
\end{equation}
We introduce the operator ${\mathcal L}_p: W^{l+2m,N}(K)\to
\mathcal W^{l,N}(K,\gamma)$ given by
$$
 {\mathcal L}_pU=\{{\mathcal P}_j(D_y)U_j,\
 {\mathcal B}_{j\sigma\mu}(D_y)U|_{\gamma_{j\sigma}}\}
$$
and corresponding to problem~(\ref{eqPinK0}), (\ref{eqBinK0}).

Let us write the operators ${\mathcal P}_j(D_y)$ and
$B_{j\sigma\mu ks}(D_y)$ in polar coordinates: ${\mathcal
P}_j(D_y)=r^{-2m}\tilde{\mathcal P}_j (\omega, D_\omega, rD_r),$
$B_{j\sigma\mu ks}(D_y)=r^{-m_{j\sigma\mu}}\tilde B_{j\sigma\mu
ks}(\omega, D_\omega, rD_r)$.

We introduce the spaces of vector-functions
\begin{gather*}
 W^{l+2m,N}(-b, b)=\prod_{j=1}^{N} W^{l+2m}(-b_j, b_j),\
 \mathcal W^{l,N}[-b, b]=\prod_{j=1}^{N} {\mathcal W}^{l}[-b_j,
 b_j],\\
 {\mathcal W}^{l}[-b_j, b_j]=W^l(-b_j, b_j) \times{\mathbb C}^{2m}
\end{gather*}
and consider the analytic operator-valued function
$\tilde{\mathcal L}_p(\lambda):W^{l+2m,N}(-b, b)\to {\mathcal
W}^{l,N}[-b, b]$ given by
$$
 \tilde{\mathcal L}_p(\lambda)\varphi=\{\tilde{\mathcal P}_j(\omega, D_\omega, \lambda)\varphi_j,\
  \sum\limits_{k,s} (\chi_{j\sigma ks})^{i\lambda-m_{j\sigma\mu}}
 {\tilde B}_{j\sigma\mu ks}(\omega, D_\omega, \lambda)
              \varphi_k(\omega+\omega_{j\sigma ks})|_{\omega=(-1)^\sigma
              b_j}\}.
$$
Main definitions and facts concerning eigenvalues, eigenvectors,
and associate vectors of analytic operator-valued functions can be
found in~\cite{GS}. In the sequel, it will be on principle that
the spectrum of the operator $\tilde{\mathcal L}_p(\lambda)$ is
discrete (see Lemma~2.1~\cite{SkDu90}).

\smallskip

Further, we will show that the Fredholm solvability of
problem~(\ref{eqPinG}), (\ref{eqBinG}) in Sobolev spaces depends
on the location of eigenvalues of model operators~$\tilde{\mathcal
L}_p(\lambda)$ corresponding to the points of~${\mathcal K}$.
Notice that the solvability of the same problem in weighted spaces
depends on the location of eigenvalues of model operators
corresponding not only to the points of ${\mathcal K}$ but also
$\Omega_{is}({\mathcal K})\subset \bar G$ and
$\Omega_{i's'}(\Omega_{is}({\mathcal K})\cap\Upsilon_{i'})\subset
G$ (see~\cite{SkMs86, SkDu91}). This can be explained as follows:
the points of the sets indicated are connected by means of the
transformations $\Omega_{is}$. That is why singularities of
solutions appearing near the set ${\mathcal K}$ may be ``carried''
to other points both on the boundary and strictly inside the
domain. But in our case we will prove that if the right-hand side
of problem~(\ref{eqPinG}), (\ref{eqBinG}) is subject to finitely
many orthogonality conditions in the Sobolev space ${\mathcal
W}^l(G,\Upsilon)$, then the solutions belong to the Sobolev space
$W^{l+2m}(G)$. Therefore, such solutions have no singularities.

\section{Nonlocal Problems in Plane Angles in the Case where the Line
$\Im\lambda=1-l-2m$ Contains no Eigenvalues of $\tilde{\mathcal
L}_p(\lambda)$}\label{sectLinKNoEigen}

In this section, we construct an operator acting in Sobolev
spaces, defined for compactly supported functions, and being the
right inverse for the operator ${\bf L}_p$ up to the sum of small
and compact perturbations. (We remind that ${\bf L}_p$ corresponds
to model problem~\eqref{eqPinK}, \eqref{eqBinK}.)

\subsection{Weighted spaces $H_a^k(Q)$}

Throughout this section, we suppose that the orbit $\Orb_p$ is
fixed; therefore, for short, we denote the operators ${\bf L}_p$,
$\mathcal L_p$, and $\tilde{\mathcal L}_p(\lambda)$ by ${\mathfrak
L}$, $\mathcal L$, and $\tilde{\mathcal L}(\lambda)$ respectively.

The investigation of the solvability for problem~\eqref{eqPinK},
\eqref{eqBinK} in Sobolev spaces will be based upon the results on
the solvability of problem~\eqref{eqPinK0}, \eqref{eqBinK0} in
weighted spaces. Let us introduce these spaces and present some of
their properties.

\smallskip

For any set $X\in \mathbb R^n$ ($n\ge1$), we denote by
$C_0^\infty(X)$ the set of functions infinitely differentiable in
$\bar X$ and compactly supported in $X$. Let either $Q=K_j$ or
$Q=K_j\cap\{y\in\mathbb R^2:\ |y|<d\}\ (d>0)$, or $Q={\mathbb
R}^2$. Denote by $H_a^k(Q)$ the completion of the set
$C_0^\infty(\bar Q\setminus \{0\})$ with respect to the norm
$$
 \|w\|_{H_a^k(Q)}=\left(
    \sum_{|\alpha|\le k}\int\limits_Q r^{2(a-k+|\alpha|)} |D_y^\alpha w|^2 dy
                                       \right)^{1/2},
$$
where $a\in \mathbb R$, $k\ge 0$ is an integer. For $k\ge1$, we
denote by $H_a^{k-1/2}(\gamma)$ the space of traces on a smooth
curve $\gamma\subset\bar Q$ with the norm
$$
\|\psi\|_{H_a^{k-1/2}(\gamma)}=\inf\|w\|_{H_a^k(Q)} \quad (w\in
H_a^l(Q): w|_\gamma = \psi).
$$
We introduce the following spaces of vector-functions:
$$
 \begin{matrix}
 H_a^{l+2m,N}(K)=\prod\limits_{j=1}^{N} H_a^{l+2m}(K_j),\
 \mathcal H_a^{l,N}(K,\gamma)=\prod\limits_{j=1}^{N} \mathcal
 H_a^l(K_j,\gamma_j),\\
 \mathcal H_a^l(K_j,\gamma_j)=
 H_a^l(K_j)\times\prod\limits_{\sigma=1,2}\prod\limits_{\mu=1}^m
 H_a^{l+2m-m_{j\sigma\mu}-1/2}(\gamma_{j\sigma}).
 \end{matrix}
$$
The bounded operator ${\mathcal L}_{a}: H_a^{l+2m,N}(K)\to
\mathcal H_a^{l,N}(K,\gamma)$ given by
\begin{equation}\label{eqL1ap}
{\mathcal L}_{a}U=\{{\mathcal P}_j(D_y)U_j,\
 {\mathcal B}_{j\sigma\mu}(D_y)U|_{\gamma_{j\sigma}}\}
\end{equation}
corresponds to problem~(\ref{eqPinK0}), (\ref{eqBinK0}) in the
weighted spaces. From Theorem~2.1~\cite{SkDu90}, it follows that
the operator ${\mathcal L}_{a}$ has a bounded inverse if and only
if the line $\Im\lambda=1-l-2m$ contains no eigenvalues of the
operator~$\tilde{\mathcal L}(\lambda)$. Using the invertibility of
${\mathcal L}_{a}$, in this section and next one, we will study
the solvability of problems~(\ref{eqPinK0}), (\ref{eqBinK0})
and~(\ref{eqPinK}), (\ref{eqBinK}) in Sobolev spaces. To this end,
we need some auxiliary results (Lemmas~\ref{lWinHa}
and~\ref{lu-uG}) concerning the relation between the spaces
$H_a^k(\cdot)$ and $W^k(\cdot)$.

\begin{lemma}\label{lWinHa}
Let $u\in W^k(Q)$ ($k\ge2$), $u(y)=0$ for $|y|\ge 1$, and
$D^\alpha u|_{y=0}=0$ ($|\alpha|\le k-2$). Then we have
\begin{equation}\label{eqWinHa'}
 \|u\|_{H_a^k(Q)}\le c_a\|u\|_{W^k(Q)},\quad a>0.
\end{equation}
If we additionally suppose that\footnote{If some assertion is
formulated for a function $D^l u$, then it is meant to hold for
all the functions $D^\alpha u$, $|\alpha|=l$.} $D^{k-1}u\in
H_0^1(Q)$, then we have
\begin{equation}\label{eqWinHa''}
 \|u\|_{H_0^k(Q)}\le c\sum\limits_{|\alpha|=k-1}\|D^\alpha u\|_{H_0^1(Q)}.
\end{equation}
Here $Q$ is the same domain as before, $c_a>0$ is independent of
$u$.
\end{lemma}
\begin{proof}
From Lemma~4.9~\cite{KondrTMMO67}, it follows that, for each
$a>0$,
$$
 \|D^{k-1}u\|_{H_a^1(Q)}\le c \|D^{k-1}u\|_{W^1(Q)}\le c
 \|u\|_{W^k(Q)}.
$$
Combining this estimate (or the inclusion $D^{k-1}u\in H_0^1(Q)$)
with Lemma\footnote{Lemma~4.12~\cite{KondrTMMO67} is proved by
Kondrat'ev for $a=0$; however, his proof remains true, with slight
modifications, for all $a<1$.}~4.12~\cite{KondrTMMO67} yields
inequality~\eqref{eqWinHa'} for $0<a<1$ (or
inequality~\eqref{eqWinHa''} respectively). Since the support of
$u$ is compact, it follows that inequality~\eqref{eqWinHa'} holds
for all $a>0$.
\end{proof}
\begin{lemma}\label{lu-uG}
Let $u\in W^1({\mathbb R}^2)$ and $u(y)=0$ for $|y|\ge 1$. Then we
have
 $$
 \|u(y)-u({\mathcal G}_0y)\|_{H_0^1({\mathbb R}^2)}\le c\|u\|_{W^1({\mathbb
 R}^2)},
 $$
where ${\mathcal G}_0$ is a composition of rotation by an angle
$\omega_0$ ($-\pi<\omega_0\le\pi$) and expansion $\chi_0$
($\chi_0>0$) times.
\end{lemma}
\begin{proof}
Writing a function $u$ in polar coordinates $(\omega, r)$ yields
$$
 u(y)-u({\mathcal G}_0y)=u(\omega, r)-u(\omega+\omega_0, \chi_0 r)=v_1+v_2,
$$
where $v_1(\omega, r)=u(\omega, r)-u(\omega+\omega_0, r)$,
$v_2(\omega, r)=u(\omega+\omega_0, r)-u(\omega+\omega_0, \chi_0
r)$.

Let us consider the function $v_1$. By
Lemma~4.15~\cite{KondrTMMO67}, we obtain
$$
 \int\limits_0^\infty r^{-1}|v_1(0, r)|^2dr\le k_1\|u\|_{W^1({\mathbb R}^2)}.
$$
From this and Lemma~4.8~\cite{KondrTMMO67}, it follows that
$v_1\in H_0^1({\mathbb R}^2)$ and
\begin{equation}\label{equ-uG1}
 \|v_1\|_{H_0^1({\mathbb R}^2)}\le k_2\|u\|_{W^1({\mathbb R}^2)}.
\end{equation}
To prove the lemma, it remains to show that
\begin{equation}\label{equ-uG2}
 \int\limits_{{\mathbb R}^2} r^{-2}|v_2|^2dy\le k_3\|u\|_{W^1({\mathbb R}^2)}.
\end{equation}
For $\chi_0>1$ (the case where $0<\chi_0<1$ can be considered
analogously), we have
$$
 \int\limits_{{\mathbb R}^2} r^{-2}|v_2|^2dy=\int\limits_{-\pi}^{\pi}d\omega\int\limits_0^\infty r^{-1}|v_2(\omega, r)|^2dr
=\int\limits_{-\pi+\omega_0}^{\pi+\omega_0}d\omega\int\limits_0^\infty
r^{-1}dr\left|\int\limits_r^{\chi_0r}\frac{\partial u(\omega,
t)}{\partial t}dt\right|^2.
$$
Using the Schwarz inequality, followed by the change of
integration limits, we get estimate~\eqref{equ-uG2}:
\begin{multline}\notag
\int\limits_{\mathbb R^2} r^{-2}|v_2|^2dy\le
(\chi_0-1)\int\limits_{-\pi+\omega_0}^{\pi+\omega_0}d\omega\int\limits_0^\infty
dr
\int\limits_r^{\chi_0r}\left|\frac{\partial u(\omega, t)}{\partial t}\right|^2 dt\\
=\frac{(\chi_0-1)^2}{\chi_0}\int\limits_{-\pi+\omega_0}^{\pi+\omega_0}d\omega\int\limits_0^\infty
\left|\frac{\partial u(\omega, t)}{\partial t}\right|^2
t\,dt\le\frac{(\chi_0-1)^2}{\chi_0}\|u\|^2_{W^1({\mathbb R}^2)}.
\end{multline}
\end{proof}

Let us prove one more auxiliary result.
\begin{lemma}\label{lSmallComp}
Let $H,\ H_1,$ and $H_2$ be Hilbert spaces, ${\mathfrak A}: H\to
H_1$ a linear bounded operator, and ${\cal T}: H\to H_2$ a compact
operator. Suppose that, for some $\varepsilon>0,\ c>0$, and $f\in
H$, the following inequality holds:
 \begin{equation}\label{eqSmallComp}
  \|{\mathfrak A}f\|_{H_1}\le \varepsilon\|f\|_H+c\|{\cal T} f\|_{H_2}.
 \end{equation}
Then there exist operators ${\cal M}, {\cal F}:H\to H_1$ such that
 $$
  {\mathfrak A}={\cal M}+{\cal F},
 $$
$\|{\cal M}\|\le 2\varepsilon$, and the operator ${\cal F}$ is
finite-dimensional.
\end{lemma}
\begin{proof}
As is well known (see, e.g., \cite[Ch.~5, \S~85]{RS}), any compact
operator is the limit of a uniformly convergent sequence of
finite-dimensional operators. Therefore, there exist bounded
operators ${\cal M}_0, {\cal F}_0:H\to H_2$ such that ${\cal
T}={\cal M}_0+{\cal F}_0$, $\|{\cal M}_0\|\le c^{-1}\varepsilon$,
and the operator ${\cal F}_0$ is finite-dimensional. This
and~(\ref{eqSmallComp}) imply
 \begin{equation}\label{eqSmallComp1}
\|{\mathfrak A}f\|_{H_1}\le 2\varepsilon\|f\|_H+c\|{\cal F}_0
f\|_{H_2}\quad \mbox{for all } f\in H.
 \end{equation}
Denote by $\ker({\cal F}_0)^\bot$ the orthogonal supplement in $H$
to the kernel of the operator ${\cal F}_0$. Since the
finite-dimensional operator ${\cal F}_0$ maps $\ker({\cal
F}_0)^\bot$ onto its image in a one-to-one manner, the subspace
$\ker({\cal F}_0)^\bot$ is of finite dimension. Let ${\cal I}$
denote the identity operator in $H$ and ${\cal P}_0$ the
orthogonal projector onto $\ker({\cal F}_0)^\bot$. Clearly,
${\mathfrak A}{\cal P}_0:H\to H_1$ is a finite-dimensional
operator. Furthermore, since ${\cal I}-{\cal P}_0$ is the
orthogonal projector onto $\ker({\cal F}_0)$, it follows that
${\cal F}_0({\cal I}-{\cal P}_0)=0$. Therefore, substituting the
function $({\cal I}-{\cal P}_0)f$ instead of $f$
in~(\ref{eqSmallComp1}), we get
$$
\|{\mathfrak A}({\cal I}-{\cal P}_0) f\|_{H_1}\le
2\varepsilon\|({\cal I}-{\cal
P}_0)f\|_H\le2\varepsilon\|f\|_H\quad \text{for all }
  f\in H.
$$
Denoting ${\cal M}={\mathfrak A}({\cal I}-{\cal P}_0)$ and ${\cal
F}={\mathfrak A}{\cal P}_0$ completes the proof.
\end{proof}

\subsection{Construction of the operator $\mathfrak R$}
In this subsection, we construct the operator $\mathfrak R$ acting
in a subspace $\mathcal S^{l,N}(K,\gamma)$ of the space
$W^{l,N}(K,\gamma)$, defined for compactly supported functions,
and being the right inverse for the operator $\mathfrak L$ up to
the sum of small and compact perturbations (see
Theorem~\ref{thRegL1p}). To construct such an operator, we assume
that the following condition holds.
\begin{condition}\label{condNoEigen}
The line $\Im\lambda=1-l-2m$ contains no eigenvalues of the
operator $\tilde{\mathcal L}(\lambda)$.
\end{condition}

We denote by $\mathcal S^{l,N}(K,\gamma)$ the subspace of
$\mathcal W^{l,N}(K,\gamma)$, consisting of the functions $\{f_j,
f_{j\sigma\mu}\}$ such that
\begin{equation}\label{eqSlN1}
 D^\alpha f_j|_{y=0}=0,\quad |\alpha|\le l-2,
\end{equation}
\begin{equation}\label{eqSlN2}
 \left.\frac{\partial^\beta f_{j\sigma\mu}}
 {\partial\tau_{j\sigma}^\beta}\right|_{y=0}=0,\quad \beta\le
 l+2m-m_{j\sigma\mu}-2,
\end{equation}
where $\tau_{j\sigma}$ is the unit vector directed along the
ray~$\gamma_{j\sigma}$. If $l-2<0$ or $l+2m-m_{j\sigma\mu}-2<0$,
the corresponding conditions are absent. From Sobolev's embedding
theorem and Riesz' theorem on a general form of linear continuous
functionals in Hilbert spaces, it follows that the set $\mathcal
S^{l,N}(K,\gamma)$ is closed and of finite codimension in
$\mathcal W^{l,N}(K,\gamma)$.

\smallskip

Let us consider the operators
$$
 \frac{\partial^{l+2m-m_{j\sigma\mu}-1}} {\partial
  \tau_{j\sigma}^{l+2m-m_{j\sigma\mu}-1}}{\mathcal B}_{j\sigma\mu}(D_y)U\equiv
  \frac{\partial^{l+2m-m_{j\sigma\mu}-1}} {\partial
  \tau_{j\sigma}^{l+2m-m_{j\sigma\mu}-1}} \Big(\sum\limits_{k,s}(B_{j\sigma\mu
  ks}(D_y)U_k)({\mathcal G}_{j\sigma ks}y)\Big).
$$
Using the chain rule, we can write
\begin{equation}\label{eqDiffB}
 \frac{\partial^{l+2m-m_{j\sigma\mu}-1}}
 {\partial \tau_{j\sigma}^{l+2m-m_{j\sigma\mu}-1}}
 {\mathcal B}_{j\sigma\mu}(D_y)U\equiv
 \sum\limits_{k,s}(\hat B_{j\sigma\mu ks}(D_y)U_k)({\mathcal G}_{j\sigma ks}y),
\end{equation}
where $\hat B_{j\sigma\mu ks}(D_y)$ are some homogenous
differential operators of order $l+2m-1$ with constant
coefficients. In particular, we have $\hat B_{j\sigma\mu j0}(D_y)=
 \dfrac{\partial^{l+2m-m_{j\sigma\mu}-1}}
 {\partial \tau_{j\sigma}^{l+2m-m_{j\sigma\mu}-1}}
B_{j\sigma\mu j0}(D_y)$ since ${\mathcal G}_{j\sigma j0}y\equiv
y$. Formally replacing the nonlocal operators in~(\ref{eqDiffB})
by the corresponding local ones, we introduce the operators
\begin{equation}\label{eqSystemB}
 \hat{\mathcal B}_{j\sigma\mu}(D_y)U\equiv
 \sum\limits_{k,s}\hat B_{j\sigma\mu ks}(D_y)U_k(y).
\end{equation}
Along with system~(\ref{eqSystemB}), we consider (for $l\ge 1$)
the operators
\begin{equation}\label{eqSystemP}
 D^\xi{\mathcal P}_j(D_y)U_j(y),\quad |\xi|=l-1.
\end{equation}
The system of operators~(\ref{eqSystemB}) and~(\ref{eqSystemP})
plays an essential in the proof of the following lemma, which is
used for the construction of the operator $\mathfrak R$.

\begin{lemma}\label{lA}
Let Condition~\ref{condNoEigen} hold. Then, for any $\varepsilon$,
$0<\varepsilon<1$, there exists a bounded operator
$$
\mathcal A:
\{f\in \mathcal S^{l,N}(K,\gamma): \supp f\subset\mathcal
O_{\varepsilon}(0)\}\to W^{l+2m,N}(K)
$$
such that, for any $f=\{f_j, f_{j\sigma\mu}\}\in {\rm
Dom}(\mathcal A)$, the function $V=\mathcal A f$ satisfies the
following conditions: $V=0$ for $|y|\ge 1$,
\begin{equation}\label{eqA'}
\|{\mathcal L}V-f\|_{\mathcal H_0^{l,N}(K)}\le c\|f\|_{\mathcal
W^{l,N}(K,\gamma)},
\end{equation}
\begin{equation}\label{eqA''}
 \|V\|_{H_a^{l+2m,N}(K)}\le c_a\|f\|_{\mathcal W^{l,N}(K,\gamma)}
 \quad \mbox{for any}\ a>0.
\end{equation}
\end{lemma}
\begin{proof}
1. We introduce the operator
\begin{equation}\label{eqA1}
 f_{j\sigma\mu}\mapsto \Phi_{j\sigma\mu}
\end{equation}
taking a function $f_{j\sigma\mu}\in
W^{l+2m-m_{j\sigma\mu}-1/2}(\gamma_{j\sigma})$ to its extension
$\Phi_{j\sigma\mu}\in W^{l+2m-m_{j\sigma\mu}}({\mathbb R}^2)$ to
${\mathbb R}^2$, satisfying $\Phi_{j\sigma\mu}=0$ for $|y|\ge2$.
We also consider an extension of the function $f_j$ from $K_j$ to
${\mathbb R}^2$ so that the extended function (which we also
denote by $f_j$) is equal to zero for $|y|\ge2$. The corresponding
extension operators can be chosen linear and bounded
(see~\cite[Ch.~6, \S~3]{Stein}).

Let us consider the following linear algebraic system for all
partial derivatives $D^\alpha W_j$, $|\alpha|=l+2m-1$, $j=1,\dots,
N$:
\begin{equation}\label{eqA2}
 \hat{\mathcal B}_{j\sigma\mu}(D_y)W=
   \frac{\partial^{l+2m-m_{j\sigma\mu}-1}}
        {\partial\tau_{j\sigma}^{l+2m-m_{j\sigma\mu}-1} }
    \Phi_{j\sigma\mu},
\end{equation}
\begin{equation}\label{eqA3}
  D^\xi\mathcal P_j(D_y)W_j=D^\xi f_j
\end{equation}
($j=1,\dots, N$; $\sigma=1, 2$; $\mu=1,\dots, m$; $|\xi|=l-1$). We
remind that each of the operators $\hat{\mathcal
B}_{j\sigma\mu}(D_y)$ given by~(\ref{eqSystemB}) is the sum of
``local'' operators, which allows us to regard
system~(\ref{eqA2}), (\ref{eqA3}) as an algebraic one. Let us take
for granted that system~(\ref{eqA2}), (\ref{eqA3}) admits a unique
solution for any right-hand side. Denote by $W_{j\alpha}$ a
solution of system~(\ref{eqA2}), (\ref{eqA3}). It is obvious that
$W_{j\alpha}\in W^1(\mathbb R^2)$ and $W_{j\alpha}=0$ for
$|y|\ge2$. By virtue of Lemma~4.17~\cite{KondrTMMO67}, there
exists a linear bounded operator
\begin{equation}\label{eqA4}
 \{W_{j\alpha}\}_{|\alpha|=l+2m-1}\mapsto V_j
\end{equation}
taking a system $\{W_{j\alpha}\}_{|\alpha|=l+2m-1}\in
\prod\limits_{|\alpha|=l+2m-1}W^1(\mathbb R^2)$ to a function
$V_j\in W^{l+2m}(\mathbb R^2)$ such that $V_j(y)=0$ for $|y|\ge
1$,
\begin{equation}\label{eqA5}
  D^\alpha V_j|_{y=0}=0,\quad |\alpha|\le l+2m-2,
\end{equation}
\begin{equation}\label{eqA6}
  D^\alpha V_j-W_{j\alpha}\in H_0^1(\mathbb R^2),\quad |\alpha|=l+2m-1.
\end{equation}

\smallskip

2. Let us show that the function $V=(V_1,\dots, V_N)$ is that we
are seeking for. Inequality~(\ref{eqA''}) follows from
relations~\eqref{eqA5}, Lemma~\ref{lWinHa}, and the boundedness of
the operator~\eqref{eqA4}.

Let us prove~(\ref{eqA'}). Since the functions $W_{j\alpha}$ are
solutions of algebraic system~(\ref{eqA2}), (\ref{eqA3}) and the
functions $V_j$ satisfy~(\ref{eqA6}), it follows that
\begin{equation}\label{eqA8}
 \hat{\mathcal B}_{j\sigma\mu}(D_y)V-\frac{\partial^{l+2m-m_{j\sigma\mu}-1}}
        {\partial\tau_{j\sigma}^{l+2m-m_{j\sigma\mu}-1} }
    \Phi_{j\sigma\mu}\in H_0^1(\mathbb R^2),
\end{equation}
\begin{equation}\label{eqA9}
 D^{l-1}(\mathcal P_j(D_y)V_j-f_j)\in H_0^1(\mathbb R^2).
\end{equation}

Furthermore, from~(\ref{eqA5}) and~(\ref{eqSlN1}), we get
$$
 D^\alpha(\mathcal P_j(D_y)V_j-f_j)|_{y=0}=0,\quad |\alpha|\le l-2.
$$
Combining this with relations~(\ref{eqA9}) and Lemma~\ref{lWinHa},
we see that $\mathcal P_j(D_y)V_j-f_j\in H_0^l(K_j)$.

Now let us show that
\begin{equation}\label{eqA10}
  {\mathcal B}_{j\sigma\mu}(D_y)V|_{\gamma_{j\sigma}}-f_{j\sigma\mu}\in
  H_0^{l+2m-m_{j\sigma\mu}-1/2}(\gamma_{j\sigma}).
\end{equation}
To this end, we pass in~(\ref{eqA8}) from the ``local'' operators
$\hat{\mathcal B}_{j\sigma\mu}(D_y)$ back to the nonlocal ones
$\dfrac{\partial^{l+2m-m_{j\sigma\mu}-1}}
 {\partial \tau_{j\sigma}^{l+2m-m_{j\sigma\mu}-1}}{\mathcal
 B}_{j\sigma\mu}(D_y)$. Then, using Lemma~\ref{lu-uG}, we obtain
 from~(\ref{eqA8}):
\begin{equation}\label{eqA11}
\frac{\partial^{l+2m-m_{j\sigma\mu}-1}}
 {\partial \tau_{j\sigma}^{l+2m-m_{j\sigma\mu}-1}}(\mathcal
 B_{j\sigma\mu}(D_y)V-\Phi_{j\sigma\mu})\in H_0^1(\mathbb R^2).
\end{equation}
Inclusions~(\ref{eqA11}) and Lemma~4.18~\cite{KondrTMMO67} imply
\begin{multline}\label{eqA12}
 \int\limits_0^\infty r^{-1}\left|\frac{\partial^{l+2m-m_{j\sigma\mu}-1}}
 {\partial \tau_{j\sigma}^{l+2m-m_{j\sigma\mu}-1}}(\mathcal
 B_{j\sigma\mu}(D_y)V|_{\gamma_{j\sigma}}-f_{j\sigma\mu})\right|^2dr\\
 \le k_1\left\|\frac{\partial^{l+2m-m_{j\sigma\mu}-1}}
 {\partial \tau_{j\sigma}^{l+2m-m_{j\sigma\mu}-1}}(\mathcal
 B_{j\sigma\mu}(D_y)V-\Phi_{j\sigma\mu})\right\|^2_{H_0^1(K_j)}.
\end{multline}
From inequality~(\ref{eqA12}), relations~(\ref{eqSlN1})
and~(\ref{eqA5}), and Lemma~4.7~\cite{KondrTMMO67}, it follows
that
\begin{equation}\label{eqA13}
  \int\limits_0^\infty r^{1-2(l+2m-m_{j\sigma\mu})}
 |{\mathcal B}_{j\sigma\mu}(D_y)V|_{\gamma_{j\sigma}}-f_{j\sigma\mu}|^2 dr
 \le k_2 \left\|\frac{\partial^{l+2m-m_{j\sigma\mu}-1}}
 {\partial \tau_{j\sigma}^{l+2m-m_{j\sigma\mu}-1}}(\mathcal
 B_{j\sigma\mu}(D_y)V-\Phi_{j\sigma\mu})\right\|^2_{H_0^1(K_j)}.
\end{equation}
Combining this with the relation ${\mathcal
B}_{j\sigma\mu}(D_y)V|_{\gamma_{j\sigma}}-f_{j\sigma\mu}\in
W^{l+2m-m_{j\sigma\mu}-1/2}(\gamma_{j\sigma})$, from~(\ref{eqA13})
and Lemma~4.16~\cite{KondrTMMO67}, we get~(\ref{eqA10}). Using the
boundedness of the operators~(\ref{eqA1}) and~(\ref{eqA4}), one
can easily prove estimate~(\ref{eqA'}) as well.

\smallskip

3. Now it remains to show that system~(\ref{eqA2}), (\ref{eqA3})
admits a unique solution for any right-hand side. Obviously, this
system consists of $(l+2m)N$ equations for $(l+2m)N$ unknowns.
Therefore, it suffices to show that the corresponding homogeneous
system has only a trivial solution. We assume the contrary: there
exists a nontrivial vector of numbers $\{q_{j\alpha}\}$
($j=1,\dots, N$, $|\alpha|=l+2m-1$) such that, after substituting
the numbers $q_{j\alpha}$ instead of $D^\alpha W_j$ into the
left-hand side of system~(\ref{eqA2}), (\ref{eqA3}), its
right-hand side goes to zero. Let us consider the homogeneous
polynomial $Q_j(y)$ of order $l+2m-1$, satisfying $D^\alpha
Q_j(y)\equiv q_{j\alpha}$. Then we have $\mathcal
P_j(D_y)Q_j(y)\equiv 0$ (since $D^\xi\mathcal P_j(D_y)Q_j(y)\equiv
0$ for all $|\xi|=l-1$) and
\begin{equation}\label{eqA14}
 \hat{\mathcal B}_{j\sigma\mu}(D_y)Q(y)\equiv
 \sum\limits_{k,s}\hat B_{j\sigma\mu ks}(D_y)Q_k(y)\equiv 0\quad
 \big(Q=(Q_1,\dots, Q_N)\big).
\end{equation}
Notice that $\hat B_{j\sigma\mu ks}(D_y)Q_k(y)\equiv\const$, while
every operator $\mathcal G_{j\sigma ks}$ of rotation and expansion
takes a constant to itself. Therefore, along with~(\ref{eqA14}),
the following identity holds:
\begin{equation}\label{eqA15}
  \frac{\partial^{l+2m-m_{j\sigma\mu}-1}}
 {\partial \tau_{j\sigma}^{l+2m-m_{j\sigma\mu}-1}}\Big(
 {\mathcal B}_{j\sigma\mu}(D_y)Q(y)\Big)\equiv
 \sum\limits_{k,s}(\hat B_{j\sigma\mu ks}(D_y)Q_k)({\mathcal G}_{j\sigma ks}y)
 \equiv 0.
\end{equation}
Since ${\mathcal B}_{j\sigma\mu}(D_y)Q$ is a homogeneous
polynomial of order $l+2m-m_{j\sigma\mu}-1$, it follows
from~(\ref{eqA15}) that ${\mathcal
B}_{j\sigma\mu}(D_y)Q|_{\gamma_{j\sigma}}\equiv 0$. Thus, we see
that the vector-valued function $Q=(Q_1,\dots, Q_N)$ is a solution
to homogeneous problem~(\ref{eqPinK0}), (\ref{eqBinK0}).
Therefore,
\begin{equation}\label{eqA16}
 \begin{aligned}
& \tilde{\mathcal P}_j(\omega, D_\omega,
rD_r)\bigl(r^{l+2m-1}\tilde
 Q_j(\omega)\bigr)\equiv 0,\\
& \sum\limits_{k,s} (\chi_{j\sigma ks})^{(l+2m-1)-m_{j\sigma\mu}}
 {\tilde B}_{j\sigma\mu ks}(\omega, D_\omega, rD_r)
  \bigl(r^{l+2m-1}\tilde Q_k(\omega+\omega_{j\sigma ks})\bigr)|_{\omega=(-1)^\sigma b_j}
  \equiv 0,
 \end{aligned}
\end{equation}
where $Q_j(y)\equiv r^{l+2m-1}\tilde Q_j(\omega)$. But
identities~(\ref{eqA16}) mean that $\tilde{\mathcal
L}(-i(l+2m-1))\tilde Q(\omega)\equiv 0$, where $\tilde Q=(\tilde
Q_1,\dots, \tilde Q_N)$. This contradicts the assumption that the
line $\Im\lambda=1-l-2m$ contains no eigenvalues of
$\tilde{\mathcal L}(\lambda)$.
\end{proof}

\begin{corollary}\label{corA}
The function $V$ constructed in Lemma~\ref{lA} satisfies the
following inequality:
 \begin{equation}\label{eqA'''}
 \|\mathfrak LV-f\|_{\mathcal H_0^{l,N}(K)}
 \le c\|f\|_{\mathcal W^{l,N}(K,\gamma)}.
\end{equation}
\end{corollary}
\begin{proof}
By virtue of inequality~(\ref{eqA'}), it suffices to estimate the
differences
 $({\bf P}_j(y, D_y)-{\mathcal P}_j(D_y))V_j$ and
 $({\bf B}_{j\sigma\mu}(y, D_y)-
   {\mathcal B}_{j\sigma\mu}(D_y))V|_{\gamma_{j\sigma}}$.
The former contains the terms of the form
$$
 \bigl(a_\alpha(y)-a_\alpha(0)\bigr)D^\alpha V_j\ (|\alpha|=2m),\quad
 a_\beta(y)D^\beta V_j\ (|\beta|\le 2m-1),
$$
where $a_\alpha$ and $a_\beta$ are infinitely differentiable
functions. Fixing some $a$, $0<a<1$, taking into account that
$V=0$ for $|y|\ge1$, and using Lemma~$3.3'$~\cite{KondrTMMO67} and
inequality~(\ref{eqA''}), we obtain
 \begin{multline}\notag
 \|\bigl(a_\alpha(y)-a_\alpha(0)\bigr)D^\alpha V_j\|_{H_0^l(K_j)}\le
 k_1\|\bigl(a_\alpha(y)-a_\alpha(0)\bigr)D^\alpha V_j\|_{H_{a-1}^l(K_j)}\\
 \le k_2\|D^\alpha V_j\|_{H_{a}^l(K_j)}\le k_3\|f\|_{\mathcal W^{l,N}(K,\gamma)}.
\end{multline}
Similarly, from the definition of weighted spaces and
inequality~(\ref{eqA''}), we get
$$
 \|a_\beta(y)D^\beta V_j\|_{H_0^l(K_j)}\le
 k_4\|a_\beta(y)D^\beta V_j\|_{H_{a}^{l+1}(K_j)}
 \le k_5\|V_j\|_{H_{a}^{l+2m}(K_j)}\le k_6\|f\|_{\mathcal W^{l,N}(K,\gamma)}.
$$
The expressions $({\bf B}_{j\sigma\mu}(y, D_y)- {\mathcal
B}_{j\sigma\mu}(D_y))V|_{\gamma_{j\sigma}}$ can be estimated in
the same way.
\end{proof}

\smallskip

Using Lemma~\ref{lA}, we can construct the operator~$\mathfrak R$.

\begin{theorem}\label{thRegL1p}
Let Condition~\ref{condNoEigen} hold. Then, for any $\varepsilon$,
$0<\varepsilon<1$, there exist bounded operators
\begin{align*}
\mathfrak R&:\{f\in \mathcal S^{l,N}(K,\gamma): \supp
f\subset\mathcal O_{\varepsilon}(0)\}\to \{U\in W^{l+2m,N}(K):
\supp U\subset\mathcal O_{\varepsilon_1}(0)\},\\
\mathfrak M,\mathfrak T&:\{f\in \mathcal S^{l,N}(K,\gamma): \supp
f\subset\mathcal O_{\varepsilon}(0)\}\to \{f\in \mathcal
S^{l,N}(K,\gamma): \supp f\subset\mathcal O_{2\varepsilon_1}(0)\}
\end{align*}
with\footnote{We remind that the number $\varepsilon_0$ defines
the diameter for the support of the function~$\zeta$ appearing in
the definition of the nonlocal operator $\mathbf B_{i\mu}^1$ (see
Sec.~\ref{sectStatement}). In other words, the number
$\varepsilon_0$ defines the diameter for the support of the
coefficients of the model operators $B_{j\sigma\mu ks}(y, D_y)$,
$(k,s)\ne(j,0)$
(see~\eqref{eqCoefficientsB}).}\label{pFootnoteEps}\newcounter{qountFootnoteEps}
\setcounter{qountFootnoteEps}{\value{footnote}}
$\varepsilon_1={\rm max\,}\big\{\varepsilon, \varepsilon_0/{\rm
min\,}\{\chi_{j\sigma ks},1\}\big\}$ such that $\|\mathfrak
Mf\|_{\mathcal W^{l,N}(K,\gamma)}\le c\varepsilon_1\|f\|_{\mathcal
W^{l,N}(K,\gamma)}$, where $c>0$ depends only on the coefficients
of the operators $\mathcal P_j(D_y)$ and $B_{j\sigma\mu ks}(D_y)$,
the operator $\mathfrak T$ is compact, and
\begin{equation}\label{eqReg}
 \mathfrak L\mathfrak R f=f+\mathfrak Mf+\mathfrak T f.
\end{equation}
\end{theorem}
\begin{proof}
By virtue of Lemma~\ref{lA}, we have $f-{\mathcal L}\mathcal A
f\in\mathcal
 H_0^{l,N}(K,\gamma)$. Therefore,
$$
 {\mathcal L}_{0}^{-1}(f-{\mathcal L}\mathcal A f)\in H_0^{l+2m,N}(K),
$$
where ${\mathcal L}_{0}:H_0^{l+2m,N}(K)\to \mathcal
H_0^{l,N}(K,\gamma)$ is the operator given by~\eqref{eqL1ap} for
$a=0$. Put
$$
  \mathfrak R f=\psi U,\quad U={\mathcal L}_{0}^{-1}(f-{\mathcal L}\mathcal A f)+\mathcal A
  f.
$$
Here $\psi\in C_0^\infty(\mathbb R^2)$ is such that $\psi(y)=1$
for $|y|\le \varepsilon_1={\rm max\,}\big\{\varepsilon,
\varepsilon_0/{\rm min\,}\{\chi_{j\sigma ks},1\}\big\}$,
$\supp\psi\subset\mathcal O_{2\varepsilon_1}(0)$, and $\psi$ does
not depend on polar angle $\omega$. Let us show that the operator
$\mathfrak R$ is that we are seeking for. Using the continuity of
the embedding $H_0^{l+2m,N}(K)\subset W^{l+2m,N}(K)$, which is
valid for compactly supported functions, inequality~\eqref{eqA'},
and boundedness of the operators $\mathcal A$, we get
$$
 \|\mathfrak R f\|_{W^{l+2m,N}(K)}\le c\|f\|_{W^{l+2m,N}(K)}.
$$

Let us prove relation~\eqref{eqReg}. Since $\mathcal
 P_j(D_y)U_j=f_j$ and $\psi f_j=f_j$, it follows that
\begin{equation}\label{eqPRf-f}
\mathbf P_j(y,D_y)(\psi U_j)-f_j=[\mathbf
P_j(y,D_y),\psi]U_j+\psi(y)\big(\mathbf P_j(y,D_y)-\mathcal
 P_j(D_y)\big)U_j,
\end{equation}
where $[\cdot,\cdot]$ stands for the commutator.

Let $b(y)$ be an arbitrary coefficient of the operator
$B_{j\sigma\mu ks}(y,D_y)$ with $(k,s)\ne(j,0)$. By virtue
of~\eqref{eqCoefficientsB} and the choice of the function $\psi$,
we have
\begin{align*}
b(\mathcal G_{j\sigma ks}y)=0\quad &\text{for }
|y|\ge\varepsilon_0/\chi_{j\sigma ks},\\
(D_y^\alpha\psi)(\mathcal G_{j\sigma ks}y)=D_y^\alpha\psi(y)\quad
&\text{for } |y|\le\varepsilon_0/\chi_{j\sigma ks}
\end{align*}
(the latter expression, for $|y|\le\varepsilon_0/\chi_{j\sigma
ks}$, equals $1$ if $|\alpha|=0$ and equals $0$ if
$|\alpha|\ge1$). Thus, we have
\begin{equation}\label{eqCommutePsiG}
(b v D_y^\alpha\psi)(\mathcal G_{j\sigma ks}y)\equiv
D_y^\alpha\psi(y)(b v)(\mathcal G_{j\sigma ks}y)\quad\text{for any
} v.
\end{equation}
Obviously, if $(k,s)=(j,0)$, then identity~\eqref{eqCommutePsiG}
is also true. Therefore, taking into account that $\mathcal
B_{j\sigma\mu}(D_y)U|_{\gamma_{j\sigma}}=f_{j\sigma\mu}$ and $\psi
f_{j\sigma\mu}=f_{j\sigma\mu}$, we get
\begin{equation}\label{eqBRf-f}
\mathbf B_{j\sigma\mu}(y,D_y)(\psi
U)|_{\gamma_{j\sigma}}-f_{j\sigma\mu}= [\mathbf
B_{j\sigma\mu}(y,D_y),\psi]U|_{\gamma_{j\sigma}}+\psi(y)\big(\mathbf
B_{j\sigma\mu}(y,D_y)-\mathcal
B_{j\sigma\mu}(D_y)\big)U|_{\gamma_{j\sigma}}.
\end{equation}

From~\eqref{eqPRf-f}--\eqref{eqBRf-f} and Leibniz' formula, we
obtain that $\supp(\mathfrak L\mathfrak R f-f)\subset\mathcal
O_{2\varepsilon_1}(0)$ and
\begin{equation}\label{eqReg1}
\|\mathfrak L\mathfrak R f-f\|_{\mathcal W^{l,N}(K,\gamma)}\le
k_1\varepsilon_1\|f\|_{\mathcal
W^{l,N}(K,\gamma)}+k_2(\varepsilon_1)\|\psi_1
U\|_{W^{l+2m-1,N}(K)},
\end{equation}
where $\psi_1\in C_0^\infty(\mathbb R^2)$ is equal to $1$ on the
support of $\psi$. Notice that the function ${\mathcal
L}_{0}^{-1}(f-{\mathcal L}\mathcal A f)$ belongs to
$H_0^{l+2m,N}(K)$ and, therefore, vanishes at $y=0$ together with
all its derivatives of order $\le$$l+2m-2$. By virtue of
Lemma~\ref{lA} (in particular, see~\eqref{eqA5}), the function
$\mathcal A f$ possesses the same property. Hence, we have
$\mathfrak L\mathfrak R f-f\in \mathcal S^{l,N}(K,\gamma)$.

Furthermore, by virtue of Lemma~\ref{lA} and compactness of the
embedding
$$
\{\psi_1 U: U\in W^{l+2m,N}(K)\}\subset W^{l+2m-1,N}(K),
$$
the operator
$$
 f\mapsto \psi_1 U
$$
(see the second norm on the right-hand side of~\eqref{eqReg1})
compactly maps $\{f\in \mathcal S^{l,N}(K,\gamma): \supp
f\subset\mathcal O_{\varepsilon}(0)\}$ into $W^{l+2m-1,N}(K)$.
Combining this with inequality~\eqref{eqReg1} and
Lemma~\ref{lSmallComp}, we complete the proof.
\end{proof}

The operator $\mathfrak R$ has the ``defect'' that the diameter of
the support of $\mathfrak R f$ depends on $\varepsilon_0$ and
cannot be reduced by reducing the diameter of the support of $f$.
However, to construct a right regularizer for
problem~\eqref{eqPinG}, \eqref{eqBinG} in the whole of the domain
$G$, we need, along with $\mathfrak R$, its modification
$\mathfrak R'$ devoid of this defect. In the following theorem, we
construct such a modification $\mathfrak R'$ defined for the
functions $f'=\{f_{j\sigma\mu}\}$.

\begin{theorem}\label{thRegL1p'}
Let condition~\ref{condNoEigen} hold. Then, for any $\varepsilon$,
$0<\varepsilon<1$, there exist bounded operators
\begin{align*}
\mathfrak R'&:\{f': \{0,f'\}\in \mathcal S^{l,N}(K,\gamma),\ \supp
f'\subset\mathcal O_{\varepsilon}(0)\}\to \{U\in W^{l+2m,N}(K):
\supp U\subset\mathcal O_{2\varepsilon}(0)\},\\
\mathfrak M',\mathfrak T'&:\{f': \{0,f'\}\in \mathcal
S^{l,N}(K,\gamma),\ \supp f'\subset\mathcal
O_{\varepsilon}(0)\}\to \{f\in \mathcal S^{l,N}(K,\gamma): \supp
f\subset\mathcal O_{2\varepsilon_2}(0)\}
\end{align*}
with $\varepsilon_2=\varepsilon/{\rm min\,}\{\chi_{j\sigma
ks},1\}$ such that $\|\mathfrak M'f'\|_{\mathcal
W^{l,N}(K,\gamma)}\le c\varepsilon\|\{0,f'\}\|_{\mathcal
W^{l,N}(K,\gamma)}$, where $c>0$ depends only on the coefficients
of the operators $\mathcal P_j(D_y)$ and $B_{j\sigma\mu ks}(D_y)$,
the operator $\mathfrak T'$ is compact, and
$$
 \mathfrak L\mathfrak R' f'=\{0,f'\}+\mathfrak M'f'+\mathfrak T' f'.
$$
\end{theorem}
\begin{proof}
Put
$$
\mathfrak R' f'=\psi U,\quad U={\mathcal
L}_{0}^{-1}\big(\{0,f'\}-{\mathcal L}\mathcal A
\{0,f'\}\big)+\mathcal A
  \{0,f'\},
$$
where $\psi\in C_0^\infty(\mathbb R^2)$ is such that $\psi(y)=1$
for $|y|\le \varepsilon$, $\supp\psi\subset\mathcal
O_{2\varepsilon}(0)$, and $\psi$ does not depend on polar angle
$\omega$.

The subsequent proof coincides with the proof of
Theorem~\ref{thRegL1p} except for the one thing. Namely, in this
case, identity~\eqref{eqCommutePsiG} is not true; therefore,
instead of~\eqref{eqBRf-f}, we have
\begin{multline}\label{eqBRf-f'}
\mathbf B_{j\sigma\mu}(y,D_y)(\psi
U)|_{\gamma_{j\sigma}}-f_{j\sigma\mu}= [\mathbf
B_{j\sigma\mu}(y,D_y),\psi]U|_{\gamma_{j\sigma}}+\psi(y)\big(\mathbf
B_{j\sigma\mu}(y,D_y)-\mathcal
B_{j\sigma\mu}(D_y)\big)U|_{\gamma_{j\sigma}}\\
+\sum\limits_{(k,s)\ne(j,0)}\big(\psi(\mathcal G_{j\sigma
ks}y)-\psi(y)\big)\big(B_{j\sigma\mu
ks}(y,D_y)U_k\big)\big(\mathcal G_{j\sigma
ks}y\big)\big|_{\gamma_{j\sigma}}.
\end{multline}
Thus, to prove the theorem, it suffices to show that each of the
operators
\begin{equation}\label{eqRegL1p1}
 U_k\mapsto
J_{j\sigma\mu ks}=\big(\psi(\mathcal G_{j\sigma
ks}y)-\psi(y)\big)\big(B_{j\sigma\mu
ks}(y,D_y)U_k\big)\big(\mathcal G_{j\sigma
ks}y\big)\big|_{\gamma_{j\sigma}}
\end{equation}
compactly maps $W^{l+2m}(K_k)$ into
$W^{l+2m-m_{j\sigma\mu}-1/2}(\gamma_{j\sigma})$.

Notice that if $(k,s)\ne(j,0)$, the operator $\mathcal G_{j\sigma
ks}$ maps the ray $\gamma_{j\sigma}$ onto the ray
$$
\{y\in{\mathbb R}^2: r>0,\ \omega=(-1)^\sigma
b_{j}+\omega_{j\sigma ks}\}
$$
being strictly inside the angle $K_k$. Therefore, there exists a
function $\xi\in C_0^\infty\big((-b_k,b_k)\big)$ equal to $1$ at
the point $\omega=(-1)^\sigma b_{j}+\omega_{j\sigma ks}$.

Furthermore, notice that the difference $\psi(y)-\psi(\mathcal
G_{j\sigma ks}^{-1}y)$ has a compact support and vanishes near the
origin. Therefore, there exists a function $\psi_1\in
C_0^\infty(K_k)$ vanishing near the origin and equal to $1$ on the
support of the function $\xi(\omega)\big(\psi(y)-\psi(\mathcal
G_{j\sigma ks}^{-1}y)\big)$.

Thus, we have
\begin{multline}\label{eqRegL1p2}
\|J_{j\sigma\mu
ks}\|_{W^{l+2m-m_{j\sigma\mu}-1/2}(\gamma_{j\sigma})}\le k_1
\|\xi(\omega)\big(\psi(y)-\psi(\mathcal G_{j\sigma
ks}^{-1}y)\big)B_{j\sigma\mu
ks}(y,D_y)U_k\|_{W^{l+2m-m_{j\sigma\mu}}(K_k)}\\
\le k_2\|\psi_1U_k\|_{W^{l+2m}(K_k)}.
\end{multline}
Let us estimate the norm on the right-hand side of the last
inequality, applying Theorem~5.1~\cite[Ch.~2]{LM} and taking into
account that {\rm(I)} the function $\psi_1$ is compactly supported
and vanishes both near the origin and near the sides of the angle
$K_k$ and {\rm(II)} $\mathcal P_k(D_y)U_k=0$. As a result, using
Leibniz' formula, we obtain
\begin{equation}\label{eqRegL1p3}
\|J_{j\sigma\mu
ks}\|_{W^{l+2m-m_{j\sigma\mu}-1/2}(\gamma_{j\sigma})} \le k_3
\|\psi_2U_k\|_{W^{l+2m-1}(K_k)},
\end{equation}
where $\psi_2\in C_0^\infty(K_k)$ is equal to $1$ on the support
of $\psi_1$. From estimate~\eqref{eqRegL1p3} and the Rellich
theorem, it follows that the operator~\eqref{eqRegL1p1} is
compact.
\end{proof}

\begin{remark}\label{remUinS^l+2m}
It follows from the proofs of Theorems~\ref{thRegL1p}
and~\ref{thRegL1p'} that
\begin{equation}\notag
D^\alpha \mathfrak R f|_{y=0}=0,\quad  D^\alpha \mathfrak R'
f'|_{y=0}=0,\quad |\alpha|\le l+2m-2.
\end{equation}
\end{remark}

In Sec.~\ref{sectLFredH_a}, we study nonlocal problems in weighted
spaces with small values of the weight exponent $a$. The role of
model operators in weighted spaces is played by the bounded
operator ${\mathfrak L}_{a}: H_a^{l+2m,N}(K)\to \mathcal
H_a^{l,N}(K,\gamma)$ given by
$$
{\mathfrak L}_{a}U=\{{\mathbf P}_j(y,D_y)U_j,\
 {\mathbf B}_{j\sigma\mu}(y,D_y)U|_{\gamma_{j\sigma}}\}.
$$
Let us formulate the analog of Theorem~\ref{thRegL1p'} in weighted
spaces.

\begin{theorem}\label{thRegL1p'H_a}
Let the line $\Im\lambda=a+1-l-2m$ contain no eigenvalues of
$\tilde{\mathcal L}(\lambda)$. Then, for any $\varepsilon$,
$0<\varepsilon<1$, there exist bounded operators
\begin{align*}
\mathfrak R_a'&:\{f': \{0,f'\}\in \mathcal H_a^{l,N}(K,\gamma),\
\supp f'\subset\mathcal O_{\varepsilon}(0)\}\to \{U\in
H_a^{l+2m,N}(K):
\supp U\subset\mathcal O_{2\varepsilon}(0)\},\\
\mathfrak M_a',\mathfrak T_a'&:\{f': \{0,f'\}\in \mathcal
H_a^{l,N}(K,\gamma),\ \supp f'\subset\mathcal
O_{\varepsilon}(0)\}\to \{f\in \mathcal H_a^{l,N}(K,\gamma): \supp
f\subset\mathcal O_{2\varepsilon_2}(0)\}
\end{align*}
with $\varepsilon_2=\varepsilon/{\rm min\,}\{\chi_{j\sigma
ks},1\}$ such that $\|\mathfrak M_a'f'\|_{\mathcal
H_a^{l,N}(K,\gamma)}\le c\varepsilon\|\{0,f'\}\|_{\mathcal
H_a^{l,N}(K,\gamma)}$, where $c>0$ depends only on the
coefficients of the operators $\mathcal P_j(D_y)$ and
$B_{j\sigma\mu ks}(D_y)$, the operator $\mathfrak T_a'$ is
compact, and
$$
 \mathfrak L_a\mathfrak R_a' f'=\{0,f'\}+\mathfrak M_a'f'+\mathfrak T_a' f'.
$$
\end{theorem}
\begin{proof}
From Theorem~2.1~\cite{SkDu90}, it follows that the operator
${\mathcal L}_{a}$ has a bounded inverse. Put
$$
\mathfrak R_a' f'=\psi U,\quad U={\mathcal L}_{a}^{-1}\{0,f'\},
$$
where $\psi$ is the same function as in the proof of
Theorem~\ref{thRegL1p'}. The remaining part of the proof is
analogous to that of Theorem~\ref{thRegL1p'}.
\end{proof}

\section{Nonlocal Problems in Plane Angles in the Case where the Line
$\Im\lambda=1-l-2m$ Contains a Proper Eigenvalue of
$\tilde{\mathcal L}_p(\lambda)$}\label{sectLinKProperEigen}

\subsection{Spaces $\hat{\mathcal S}^{l,N}(K,\gamma)$}
In this section, we keep denoting the operators ${\bf L}_p$,
$\mathcal L_p$, and $\tilde{\mathcal L}_p(\lambda)$ by ${\mathfrak
L}$, $\mathcal L$, and $\tilde{\mathcal L}(\lambda)$ respectively.
Let us consider the situation where the line $\Im\lambda=1-l-2m$
contains eigenvalues of $\tilde{\mathcal L}(\lambda)$. Let
$\lambda=\lambda_0$ be one of such eigenvalues.
\begin{definition}\label{defRegEigVal}
We say that $\lambda=\lambda_0$ is a {\em proper eigenvalue} if
{\rm(I)} neither of the corresponding eigenvectors
$\varphi(\omega)=(\varphi_{j}(\omega),\dots, \varphi_{N}(\omega))$
has associate ones and {\rm(II)} the functions
$r^{i\lambda_0}\varphi_{j}(\omega)$, $j=1, \dots, N$, are
polynomials with respect to $y_1, y_2$.
\end{definition}
\begin{definition}\label{defRegEigValImproper}
An eigenvalue $\lambda=\lambda_0$ which is not proper is said to
be an {\em improper eigenvalue}.
\end{definition}

\begin{remark}
The notion of a proper eigenvalue was originally proposed by
Kondrat'ev~\cite{KondrTMMO67} for ``local'' elliptic
boundary-value problems in angular or conical domains.
\end{remark}

Clearly, if $\lambda_0$ is a proper eigenvalue, then
$\Re\lambda_0=0$. Therefore, the line $\Im\lambda=1-l-2m$ may
contain at most one proper eigenvalue. In this section, we
investigate the case where the following condition holds.

\begin{condition}\label{condProperEigen}
The line $\Im\lambda=1-l-2m$ contains only the eigenvalue
$\lambda_0=i(1-l-2m)$ and it is proper.
\end{condition}

In thats case, the conclusion of Lemma~\ref{lA} is not true, since
algebraic system~(\ref{eqA2}), (\ref{eqA3}) may have no solutions
for some right-hand sides and the system of
operators~(\ref{eqSystemB}), (\ref{eqSystemP}) is not linearly
independent. Indeed, let
$\varphi(\omega)=(\varphi_{1}(\omega),\dots, \varphi_{N}(\omega))$
be an eigenvector corresponding to the proper eigenvalue
$\lambda_0=i(1-l-2m)$. Then, by the definition of a proper
eigenvalue, $Q_j(y)=r^{l+2m-1}\varphi_j(\omega)$ is an $l+2m-1$
order polynomial (obviously, homogeneous) with respect to $y=(y_1,
y_2)$. Repeating the arguments of item~3 in the proof of
Lemma~\ref{lA}, we see that, after substituting
$q_{j\alpha}=D^\alpha Q_j$ instead of $D^\alpha W_j$ into the
left-hand side of system~(\ref{eqA2}), (\ref{eqA3}), its
right-hand side goes to zero. Therefore, system~(\ref{eqSystemB}),
(\ref{eqSystemP}) is linearly dependent. Nevertheless, provided
Condition~\ref{condProperEigen} holds, it turns out to be possible
to construct an operator $\hat{\mathfrak R}$ defined for compactly
supported functions from a certain space $\hat{\mathcal
S}^{l,N}(K,\gamma)$ and being the right inverse for $\mathfrak L$
(see Theorem~\ref{thRegL1p0}). However, in contrast to $\mathcal
S^{l,N}(K,\gamma)$, the set $\hat{\mathcal S}^{l,N}(K,\gamma)$ is
not closed in the topology of the space $\mathcal
W^{l,N}(K,\gamma)$.

\smallskip

We choose from system~(\ref{eqSystemB}) consisting of homogeneous
$l+2m-1$ order operators a maximum number of linearly independent
operators and denote them by
\begin{equation}\label{eqSystemB'}
 \hat{\mathcal B}_{j'\sigma'\mu'}(D_y)U.
\end{equation}
Any operator $\hat{\mathcal B}_{j\sigma\mu}(D_y)$ which is not
included in system~(\ref{eqSystemB'}) can be represented in the
following form:
\begin{equation}\label{eqBviaB'}
 \hat{\mathcal B}_{j\sigma\mu}(D_y)U=\sum\limits_{j',\sigma',\mu'}
 p_{j\sigma\mu}^{j'\sigma'\mu'}\hat{\mathcal B}_{j'\sigma'\mu'}(D_y)U,
\end{equation}
where $p_{j\sigma\mu}^{j'\sigma'\mu'}$ are some constants.

Let us consider the functions $f=\{f_j, f_{j\sigma\mu}\}\in
\mathcal W^{l,N}(K,\gamma)$ satisfying
\begin{equation}\label{eqRelPhi}
 \mathcal T_{j\sigma\mu}f\equiv\frac{\partial^{l+2m-m_{j\sigma\mu}-1}}
 {\partial \tau_{j\sigma}^{l+2m-m_{j\sigma\mu}-1}}\Phi_{j\sigma\mu}-
 \sum\limits_{j',\sigma',\mu'}p_{j\sigma\mu}^{j'\sigma'\mu'}
 \frac{\partial^{l+2m-m_{j'\sigma'\mu'}-1}}
 {\partial \tau_{j'\sigma'}^{l+2m-m_{j'\sigma'\mu'}-1}}\Phi_{j'\sigma'\mu'}\in
 H_0^1(\mathbb R^2).
\end{equation}
Here indices $j',\sigma',\mu'$ correspond to
operators~(\ref{eqSystemB'}) while indices $j,\sigma,\mu$
correspond to the operators from system~(\ref{eqSystemB}) that are
not included in~(\ref{eqSystemB'}); $\Phi_{j\sigma\mu}$ are the
fixed extensions of the functions $f_{j\sigma\mu}$ to $\mathbb
R^2$, defined by the operator~(\ref{eqA1});
$p_{j\sigma\mu}^{j'\sigma'\mu'}$ are the constants appearing in
relation~(\ref{eqBviaB'}). If system~(\ref{eqSystemB}) is linearly
independent, then the set of conditions~(\ref{eqRelPhi}) is empty.

Notice that the fulfilment of conditions~(\ref{eqRelPhi}) does not
depend on the choice of the extension of $f_{j\sigma\mu}$ to
$\mathbb R^2$. Indeed, let $\hat\Phi_{j\sigma\mu}$ be an extension
different from $\Phi_{j\sigma\mu}$. Then we have
$(\Phi_{j\sigma\mu}-\hat\Phi_{j\sigma\mu})|_{\gamma_{j\sigma}}=0$;
therefore, by Theorem~4.8~\cite{KondrTMMO67}, $
 \dfrac{\partial^{l+2m-m_{j\sigma\mu}-1}}
 {\partial \tau_{j\sigma}^{l+2m-m_{j\sigma\mu}-1}}
 (\Phi_{j\sigma\mu}-\hat\Phi_{j\sigma\mu})\in
 H_0^1(\mathbb R^2).
$

\smallskip

Now let us complete system~(\ref{eqSystemB'}) with $l+2m-1$ order
operators from system~(\ref{eqSystemP}) so that the resulting
system consist of linearly independent operators
\begin{equation}\label{eqSystemBP'}
 \hat{\mathcal B}_{j'\sigma'\mu'}(D_y)U,\ D^{\xi'}\mathcal P_{j'}(D_y)U_{j'}
\end{equation}
and any operator $D^{\xi}\mathcal P_{j}(D_y)U_{j}$ not included
in~(\ref{eqSystemBP'}) be represented in the following form:
\begin{equation}\label{eqPviaBP'}
 D^\xi \mathcal P_{j}(D_y)U_j=\sum\limits_{j',\sigma',\mu'}
 p_{j\xi}^{j'\sigma'\mu'}\hat{\mathcal B}_{j'\sigma'\mu'}(D_y)U+
 \sum\limits_{j',\xi'}p_{j\xi}^{j'\xi'}D^{\xi'}\mathcal P_{j'}(D_y)U_{j'},
\end{equation}
where $p_{j\xi}^{j',\sigma',\mu'}$ and $p_{j\xi}^{j',\xi'}$ are
some constants.

Let us extend the components $f_j\in W^l(K_j)$ of the vector $f$
to $\mathbb R^2$. The extended functions are also denoted by
$f_j\in W^l(\mathbb R^2)$. We consider the functions $f$
satisfying
\begin{equation}\label{eqRelPhif}
 \mathcal T_{j\xi} f \equiv D^\xi f_j-\sum\limits_{j',\sigma',\mu'}
 p_{j\xi}^{j'\sigma'\mu'}\frac{\partial^{l+2m-m_{j'\sigma'\mu'}-1}}
 {\partial \tau_{j'\sigma'}^{l+2m-m_{j'\sigma'\mu'}-1}}\Phi_{j'\sigma'\mu'}
 -\sum\limits_{j',\xi'}p_{j\xi}^{j'\xi'}D^{\xi'}f_{j'}\in H_0^1(\mathbb R^2).
\end{equation}
Here indices $j', \sigma', \mu'$ and $j', \xi'$ correspond to the
operators~(\ref{eqSystemBP'}) while indices $j, \xi$ correspond to
the operators from system~(\ref{eqSystemP}) that are not included
in~(\ref{eqSystemBP'}); $p_{j\xi}^{j'\sigma'\mu'}$ and
$p_{j\xi}^{j'\xi'}$ are the constants appearing in
relations~(\ref{eqPviaBP'}). Similarly to the above, one can show
that the fulfilment of conditions~(\ref{eqRelPhif}) does not
depend on the choice of the extension of $f_j$ and
$f_{j\sigma\mu}$ to $\mathbb R^2$. Notice that the set of
conditions~(\ref{eqRelPhif}) is empty if either $l=0$ or $l\ge1$
but system~(\ref{eqSystemBP'}) contains all the operators
from~(\ref{eqSystemP}).

Let us introduce the analog of the set $\mathcal
S^{l,N}(K,\gamma)$ in the case where
Condition~\ref{condProperEigen} holds. We denote by $\hat{\mathcal
S}^{l,N}(K,\gamma)$ the set of functions $f\in\mathcal
W^{l,N}(K,\gamma)$ satisfying conditions~(\ref{eqSlN1}),
(\ref{eqSlN2}), (\ref{eqRelPhi}), and (\ref{eqRelPhif}). Supplying
$\hat{\mathcal S}^{l,N}(K,\gamma)$ with the norm
\begin{equation}\label{eqNormS0lN}
 \|f\|_{\hat{\mathcal S}^{l,N}(K,\gamma)}=\Big(\|f\|^2_{\mathcal W^{l,N}(K,\gamma)}+
 \sum\limits_{j,\sigma,\mu}\|\mathcal T_{j\sigma\mu}f\|^2_{H_0^1(\mathbb R^2)}+
 \sum\limits_{j,\xi}\|\mathcal T_{j\xi} f\|^2_{H_0^1(\mathbb R^2)}\Big)^{1/2}
\end{equation}
makes it a complete space. (In the definition of the
norm~(\ref{eqNormS0lN}), indices $j, \sigma, \mu$ and $j, \xi$
correspond to the operators not included in
system~(\ref{eqSystemBP'}).)

\smallskip

Let us establish some important properties of the space
$\hat{\mathcal S}^{l,N}(K,\gamma)$. The following lemma shows that
if we impose on a compactly supported function $U\in
W^{l+2m,N}(K)$ finitely many orthogonality conditions of the form
 \begin{equation}\label{eqConnectUS0lN}
  D^\alpha U|_{y=0}=0,\quad |\alpha|\le l+2m-2,
 \end{equation}
then the right-hand side of the corresponding nonlocal problem
belongs to $\hat{\mathcal S}^{l,N}(K,\gamma)$.
\begin{lemma}\label{lConnectUS0lN}
Let Condition~\ref{condProperEigen} hold. Suppose that $U\in
W^{l+2m,N}(K)$, $\supp U\subset\mathcal
O_{\varepsilon\min\{\chi_{j\sigma ks}, 1\}}(0)$, and
relations~\eqref{eqConnectUS0lN} hold. Then we have
\begin{equation}\label{eqConnectUS0lN'}
  \|{\mathcal L} U\|_{\hat{\mathcal S}^{l,N}(K,\gamma)}\le
  c\|U\|_{W^{l+2m,N}(K)},\quad
  \|{\mathfrak L} U\|_{\hat{\mathcal S}^{l,N}(K,\gamma)}\le
  c\|U\|_{W^{l+2m,N}(K)}.
\end{equation}
\end{lemma}
\begin{proof}
1. Put $f=\{f_j, f_{j\sigma\mu}\}={\mathcal L} U$. From the
assumptions of the lemma, it follows that $f\in\mathcal
W^{l,N}(K,\gamma)$, $\supp f\subset\mathcal O_{\varepsilon}(0)$,
and the functions $f_j$ and $f_{j\sigma\mu}$ satisfy
relations~(\ref{eqSlN1}) and~(\ref{eqSlN2}) respectively.

We denote by $\Phi_{j\sigma\mu}\in W^{l+2m-m_{j\sigma\mu}}(\mathbb
R^2)$ the extension of $f_{j\sigma\mu}$ defined by the
operator~\eqref{eqA1}. Let us show that
\begin{equation}\label{eqConnectUS0lN1}
  \hat{\mathcal B}_{j\sigma\mu}(D_y)U-\frac{\partial^{l+2m-m_{j\sigma\mu}-1}}
 {\partial \tau_{j\sigma}^{l+2m-m_{j\sigma\mu}-1}}\Phi_{j\sigma\mu}\in
 H_0^1(\mathbb R^2).
\end{equation}
By Lemma~\ref{lu-uG}, we have
 $\hat{\mathcal B}_{j\sigma\mu}(D_y)U-\dfrac{\partial^{l+2m-m_{j\sigma\mu}-1}}
 {\partial \tau_{j\sigma}^{l+2m-m_{j\sigma\mu}-1}}{\mathcal B}_{j\sigma\mu}(D_y)U\in
 H_0^1(\mathbb R^2)$; thus, to prove~(\ref{eqConnectUS0lN1}), it suffices to show that
\begin{equation}\label{eqConnectUS0lN2}
  \frac{\partial^{l+2m-m_{j\sigma\mu}-1}}
 {\partial \tau_{j\sigma}^{l+2m-m_{j\sigma\mu}-1}}
 ({\mathcal B}_{j\sigma\mu}(D_y)U-\Phi_{j\sigma\mu})\in
 H_0^1(\mathbb R^2).
\end{equation}
But $\dfrac{\partial^{l+2m-m_{j\sigma\mu}-1}}
 {\partial \tau_{j\sigma}^{l+2m-m_{j\sigma\mu}-1}}
 ({\mathcal B}_{j\sigma\mu}(D_y)U-\Phi_{j\sigma\mu})\in
 W^1(\mathbb R^2)$ and $\dfrac{\partial^{l+2m-m_{j\sigma\mu}-1}}
 {\partial \tau_{j\sigma}^{l+2m-m_{j\sigma\mu}-1}}
 ({\mathcal
 B}_{j\sigma\mu}(D_y)U-\Phi_{j\sigma\mu})\big|_{\gamma_{j\sigma}}=0$;
hence, relation~(\ref{eqConnectUS0lN2}) follows from
Lemma~4.8~\cite{KondrTMMO67}. Thus,
relation~(\ref{eqConnectUS0lN1}) is also proved.

The operators $\hat{\mathcal B}_{j\sigma\mu}(D_y)U$ satisfy
relations~(\ref{eqBviaB'}); therefore, by virtue
of~(\ref{eqConnectUS0lN1}), the functions $\Phi_{j\sigma\mu}$
satisfy relations~(\ref{eqRelPhi}).

Similarly, from~(\ref{eqConnectUS0lN1}), equalities $\mathcal
P_j(D_y)U_j-f_j=0$, and relations~(\ref{eqPviaBP'}), it follows
that the function $f$ satisfies relations~(\ref{eqRelPhif}).
Therefore, $f\in\hat{\mathcal S}^{l,N}(K,\gamma)$, and it is easy
to check that the first inequality in~(\ref{eqConnectUS0lN'})
holds.

\smallskip

2. Now, to prove that ${\mathfrak L} U\in\hat{\mathcal
S}^{l,N}(K,\gamma)$, it suffices to show that
$$
D^{l-1}\bigl(\mathbf P_j(y, D_y)-\mathcal P_j(D_y)\bigr)U_j\in
H_0^1(\mathbb R^2),\quad
\frac{\partial^{l+2m-m_{j\sigma\mu}-1}}
 {\partial \tau_{j\sigma}^{l+2m-m_{j\sigma\mu}-1}}
 \bigl({\mathbf B}_{j\sigma\mu}(y, D_y)U-{\mathcal B}_{j\sigma\mu}(D_y)U\bigr)\in
 H_0^1(\mathbb R^2),
$$
where $U_j\in W^{l+2m}(\mathbb R^2)$ is an extension of $U_j\in
W^{l+2m}(K_j)$ to $\mathbb R^2$ (which is also denoted by $U_j$).
These expressions consist of the terms
$$
 \big(a_\alpha(y)-a_\alpha(0)\big)D^\alpha U_j\ (|\alpha|=l+2m-1),\quad
 a_\beta(y)D^\beta U_j\ (|\beta|\le l+2m-2),
$$
where $a_\alpha$ and $a_\beta$ are infinitely differentiable
functions.

Since $U_j\in W^{l+2m}(\mathbb R^2)$, it follows that $D^\alpha
U_j\in H_1^1(\mathbb R^2)$. This and
Lemma~$3.3'$~\cite{KondrTMMO67} imply that
$\big(a_\alpha(y)-a_\alpha(0)\big) D^\alpha U_j\in H_0^1(\mathbb
R^2)$.

The function $a_\beta D^\beta U_j$ ($|\beta|\le l+2m-2$) belongs
to $W^2(\mathbb R^2)$. From this,
relations~\eqref{eqConnectUS0lN}, and Lemma~\ref{lWinHa}, it
follows that $a_\beta D^\beta U_j\in H_a^2(\mathbb R^2)\subset
H_{a-1}^1(\mathbb R^2)$, $a>0$. Let us choose $0<a<1$; then, by
virtue of the compactness of supports of $U_j$, we get $a_\beta
D^\beta U_j\in H_0^1(\mathbb R^2)$. Furthermore, it is easy to
show that the second inequality in~(\ref{eqConnectUS0lN'}) also
holds.
\end{proof}

The following lemma shows that the set $\hat{\mathcal
S}^{l,N}(K,\gamma)$ is not closed in the topology of $\mathcal
W^{l,N}(K,\gamma)$.

\begin{lemma}\label{lS0NotClosed}
Let Condition~\ref{condProperEigen} hold. Then there exists a
family of functions $f^\delta\in \hat{\mathcal
S}^{l,N}(K,\gamma)$, $\delta>0$, such that $\supp
f^\delta\subset\mathcal O_{\varepsilon}(0)$ and $f^\delta$
converges in $\mathcal W^{l,N}(K,\gamma)$ to a function
$f^0\notin\hat{\mathcal S}^{l,N}(K,\gamma)$ as
 $\delta\to0$.
\end{lemma}
\begin{proof}
1. As was shown above, if $\lambda_0=i(1-l-2m)$ is a proper
eigenvalue of $\tilde{\mathcal L}(\lambda)$, then
system~(\ref{eqSystemB}), (\ref{eqSystemP}) is linearly dependent.
We consider the two possible cases: {\rm (a)}
system~(\ref{eqSystemB}) is linearly dependent or {\rm(b)}
system~(\ref{eqSystemB}) is linearly independent but
system~(\ref{eqSystemB}), (\ref{eqSystemP}) is linearly dependent.

\smallskip

2. First, let us suppose that system~(\ref{eqSystemB}) is linearly
dependent. Then the set of conditions~(\ref{eqRelPhi}) is not
empty. In this case, for some $j, \sigma, \mu$, the
norm~(\ref{eqNormS0lN}) contains the corresponding term
$\|\mathcal T_{j\sigma\mu}f\|_{H_0^1(\mathbb R^2)}$. We fix such
$j, \sigma, \mu$. Without loss of generality, one may assume that
$\gamma_{j\sigma}$ coincides with the axis $Oy_1$. We introduce
the functions $f^\delta=\{0, f_{j_1\sigma_1\mu_1}^\delta\}$
($0\le\delta\le1$) such that  $f_{j_1\sigma_1\mu_1}^\delta=0$ for
$(j_1, \sigma_1, \mu_1)\ne (j, \sigma, \mu)$ and
$f_{j\sigma\mu}^\delta(y_1)=\psi(y_1)y_1^{l+2m-m_{j\sigma\mu}-1+\delta}$,
where $\psi\in C_0^\infty\big([0, \infty)\big)$, $\psi(y_1)=1$ for
$0\le y_1\le\varepsilon/2$, and $\psi(y_1)=0$ for
$y_1\ge2\varepsilon/3$. Clearly,
$$
\hat\Phi_{j\sigma\mu}^\delta(y)=\psi(r)y_1^{l+2m-m_{j\sigma\mu}-1}r^\delta
$$
is an extension of the function $f_{j\sigma\mu}^\delta$ to
$\mathbb R^2$. Moreover, the extension operator defined for the
functions $f_{j\sigma\mu}^\delta$ ($0\le\delta\le1$) is bounded
from $W^{l+2m-m_{j\sigma\mu}-1/2}(\gamma_{j\sigma})$ into
$W^{l+2m-m_{j\sigma\mu}}(\mathbb R^2)$ (which follows from the
fact that
$\|f_{j\sigma\mu}^\delta\|_{W^{l+2m-m_{j\sigma\mu}-1/2}(\gamma_{j\sigma})}\ge
c_1$ and
$\|\hat\Phi_{j\sigma\mu}^\delta\|_{W^{l+2m-m_{j\sigma\mu}}(\mathbb
R^2)}\le c_2$ with $c_1, c_2>0$ being independent of
$0\le\delta\le1$).

Thus, for $0<\delta\le1$, we have
\begin{equation}\label{eqS0NotClosed1}
 \begin{aligned}
 \|f^\delta\|^2_{\mathcal W^{l,N}(K,\gamma)}&=
 \|f_{j\sigma\mu}^\delta\|^2_{W^{l+2m-m_{j\sigma\mu}-1/2}(\gamma_{j\sigma})},\\
\|f^\delta\|^2_{\hat{\mathcal S}^{l,N}(K,\gamma)}&\approx
 \|f_{j\sigma\mu}^\delta\|^2_{W^{l+2m-m_{j\sigma\mu}-1/2}(\gamma_{j\sigma})}+
 \left\|\frac{\partial^{l+2m-m_{j\sigma\mu}-1}}
 {\partial
 y_{1}^{l+2m-m_{j\sigma\mu}-1}}\hat\Phi_{j\sigma\mu}^\delta\right\|^2_{H_0^1(\mathbb
 R^2)}
  \end{aligned}
\end{equation}
(the finiteness of the norms~(\ref{eqS0NotClosed1}) for each
$\delta>0$ can be verified by straightforward calculations). Here
symbol ``$\approx$'' means that the corresponding norms are
equivalent. Furthermore, one can directly check that
$\hat\Phi_{j\sigma\mu}^\delta\to\hat\Phi_{j\sigma\mu}^0$ in
$W^{l+2m-m_{j\sigma\mu}}(\mathbb R^2)$ as $\delta\to 0$.
Therefore, $f_{j\sigma\mu}^\delta\to f_{j\sigma\mu}^0$ in
$W^{l+2m-m_{j\sigma\mu}-1/2}(\gamma_{j\sigma})$ as $\delta\to 0$.
However, the corresponding function $f^0=\{0, f_{j\sigma\mu}^0\}$
does not belong to $\hat{\mathcal S}^{l,N}(K,\gamma)$. Indeed,
assuming the contrary, by virtue of~(\ref{eqS0NotClosed1}), we
have $\dfrac{\partial^{l+2m-m_{j\sigma\mu}-1}}
 {\partial y_{1}^{l+2m-m_{j\sigma\mu}-1}}\hat\Phi_{j\sigma\mu}^0\in H_0^1(\mathbb
R^2)$, which is not true since the function
$\dfrac{\partial^{l+2m-m_{j\sigma\mu}-1}}
 {\partial y_{1}^{l+2m-m_{j\sigma\mu}-1}}\hat\Phi_{j\sigma\mu}^0$
is equal to a nonzero constant near the origin.

\smallskip

3. Now let system~(\ref{eqSystemB}) be linearly independent; then
system~(\ref{eqSystemB}), (\ref{eqSystemP}) is linearly dependent.
In this case, conditions~(\ref{eqRelPhi}) are absent but the set
of conditions~(\ref{eqRelPhif}) is not empty. Therefore, for some
$j, \xi$, the norm~(\ref{eqNormS0lN}) contains the corresponding
term $\|\mathcal T_{j\xi}f\|_{H_0^1(\mathbb R^2)}$. We fix such
$j, \xi$ and introduce the functions $f^\delta=\{f_{j_1}^\delta,
0\}$ ($0\le\delta\le1$) such that $f_{j_1}^\delta=0$ for $j_1\ne
j$ and $f_j^\delta=\psi(r)y^\xi r^\delta$. One can directly check
that $f_j^\delta\to f_j^0$ in $W^l(\mathbb R^2)$ as $\delta\to 0$,
but $f^0=\{f^0_j, f^0_{j\sigma\mu}\}\notin \hat{\mathcal
S}^{l,N}(K,\gamma)$ since $D^\xi f_j^0\notin H_0^1(\mathbb R^2)$.
\end{proof}

\smallskip

\subsection{Construction of the operator $\hat{\mathfrak R}$}
Let us prove the analog of Lemma~\ref{lA}, which will be used to
construct the operator $\hat{\mathfrak R}$ acting in the space
$\hat{\mathcal S}^{l,N}(K,\gamma)$.
\begin{lemma}\label{lA0}
Let Condition~\ref{condProperEigen} hold. Then, for any
$\varepsilon$, $0<\varepsilon<1$,  there exists a bounded operator
$$
\hat{\mathcal A}: \{f\in \hat{\mathcal
S}^{l,N}(K,\gamma): \supp f\subset\mathcal O_{\varepsilon}(0)\}\to
W^{l+2m,N}(K)
$$
such that, for any $f=\{f_j, f_{j\sigma\mu}\}\in {\rm
Dom}(\hat{\mathcal A})$, the function $V=\hat{\mathcal A} f$
satisfies the following conditions: $V=0$ for $|y|\ge 1$,
 \begin{equation}\label{eqA0'}
\|{\mathcal L}V-f\|_{\mathcal H_0^{l,N}(K)}
 \le c\|f\|_{\hat{\mathcal S}^{l,N}(K,\gamma)},
\end{equation}
and inequality~\eqref{eqA''} holds.
\end{lemma}
\begin{proof}
1. Similarly to the proof of Lemma~\ref{lA}, we consider the
algebraic system for all partial derivatives $D^\alpha W_j$,
$|\alpha|=l+2m-1$, $j=1, \dots, N$:
\begin{equation}\label{eqA01}
 \begin{aligned}
 \hat{\mathcal B}_{j'\sigma'\mu'}(D_y)W&=
   \frac{\partial^{l+2m-m_{j'\sigma'\mu'}-1}}
        {\partial\tau_{j'\sigma'}^{l+2m-m_{j'\sigma'\mu'}-1} }
    \Phi_{j'\sigma'\mu'},\\
  D^{\xi'}\mathcal P_{j'}(D_y)W_{j'}&=D^{\xi'} f_{j'},
  \end{aligned}
\end{equation}
where $\Phi_{j'\sigma'\mu'}$ and $f_{j'}$ are the extensions of
$f_{j'\sigma'\mu'}$ and $f_{j'}$ to $\mathbb R^2$ described in the
proof of Lemma~\ref{lA}. Now the left-hand side of
system~(\ref{eqA01}) contains only the operators included in
system~(\ref{eqSystemBP'}). The matrix of system~(\ref{eqA01})
consists of $(l+2m)N$ columns and $q$, $q<(l+2m)N$, linearly
independent rows. Choosing $q$ linearly independent columns and
putting the unknowns $D^\alpha W_j$ corresponding to the remaining
$(l+2m)N-q$ columns equal to zero, we obtain a system of $q$
equations for $q$ unknowns, which admits a unique solution. Thus,
we defined the linear bounded operator
\begin{equation}\label{eqA02}
 \left\{\frac{\partial^{l+2m-m_{j'\sigma'\mu'}-1}}
        {\partial\tau_{j'\sigma'}^{l+2m-m_{j'\sigma'\mu'}-1} }
    \Phi_{j'\sigma'\mu'},\ D^{\xi'} f_{j'}\right\}\mapsto \{D^\alpha W_j\}\equiv\{W_{j\alpha}\}
\end{equation}
acting from $W^{1,q}(\mathbb R^2)$ into $W^{1,(l+2m)N}(\mathbb
R^2)$ and such that $W_{j\alpha}(y)=0$ for $|y|\ge 2$. Using the
functions $D^\alpha W_j$ and the operator~(\ref{eqA4}), we get
functions $V_j$, $j=1,\dots, N$, satisfying relations~(\ref{eqA5})
and~(\ref{eqA6}). Let us show that $V=(V_1,\dots, V_N)$ is the
function we are seeking for.

\smallskip

2. Analogously to the proof of Lemma~\ref{lA}, one can prove
estimate~(\ref{eqA''}) for the function $V$. Let us prove
inequality~(\ref{eqA0'}). Since $\{W_{j\alpha}\}$ is a solution to
system~(\ref{eqA01}) and the functions $V_j$ satisfy
conditions~(\ref{eqA6}), it follows that
\begin{equation}\label{eqA03}
 \hat{\mathcal B}_{j'\sigma'\mu'}(D_y)V-\frac{\partial^{l+2m-m_{j'\sigma'\mu'}-1}}
        {\partial\tau_{j'\sigma'}^{l+2m-m_{j'\sigma'\mu'}-1} }
    \Phi_{j'\sigma'\mu'}\in H_0^1(\mathbb R^2),
\end{equation}
\begin{equation}\label{eqA04}
 D^{\xi'}(\mathcal P_{j'}(D_y)V_{j'}-f_{j'})\in H_0^1(\mathbb R^2).
\end{equation}
Let us consider an arbitrary operator $\hat{\mathcal
B}_{j\sigma\mu}(D_y)$ that is not included in
system~(\ref{eqSystemBP'}). Using~(\ref{eqBviaB'}), we get
\begin{multline}\label{eqA05}
  \hat{\mathcal B}_{j\sigma\mu}(D_y)V-\frac{\partial^{l+2m-m_{j\sigma\mu}-1}}
        {\partial\tau_{j\sigma}^{l+2m-m_{j\sigma\mu}-1} }\Phi_{j\sigma\mu}=
  \sum\limits_{j',\sigma',\mu'}
 p_{j\sigma\mu}^{j'\sigma'\mu'}\left(\hat{\mathcal B}_{j'\sigma'\mu'}(D_y)V-
\frac{\partial^{l+2m-m_{j'\sigma'\mu'}-1}}
{\partial\tau_{j'\sigma'}^{l+2m-m_{j'\sigma'\mu'}-1}}\Phi_{j'\sigma'\mu'}\right)\\
+\sum\limits_{j',\sigma',\mu'}
 p_{j\sigma\mu}^{j'\sigma'\mu'}\frac{\partial^{l+2m-m_{j'\sigma'\mu'}-1}}
{\partial\tau_{j'\sigma'}^{l+2m-m_{j'\sigma'\mu'}-1}}\Phi_{j'\sigma'\mu'}-\frac{\partial^{l+2m-m_{j\sigma\mu}-1}}
        {\partial\tau_{j\sigma}^{l+2m-m_{j\sigma\mu}-1} }\Phi_{j\sigma\mu}.
\end{multline}
But $f\in \hat{\mathcal S}^{l,N}(K,\gamma)$; therefore
conditions~(\ref{eqRelPhi}) hold. This and relations~(\ref{eqA03})
and~(\ref{eqA05}) imply that, for all $j, \sigma, \mu$, the
following relations hold:
\begin{equation}\label{eqA06}
  \hat{\mathcal B}_{j\sigma\mu}(D_y)V-\frac{\partial^{l+2m-m_{j\sigma\mu}-1}}
        {\partial\tau_{j\sigma}^{l+2m-m_{j\sigma\mu}-1} }\Phi_{j\sigma\mu}\in
        H_0^1(\mathbb R^2).
\end{equation}
Similarly, one can consider the operators $D^\xi P_j(D_y)$ that
are not included in system~(\ref{eqSystemBP'}) and, using
relations~(\ref{eqBviaB'}) and (\ref{eqRelPhi}), (\ref{eqPviaBP'})
and (\ref{eqRelPhif}), as well as~(\ref{eqA03}) and~(\ref{eqA04}),
prove the relations
\begin{equation}\label{eqA07}
 D^{\xi}(\mathcal P_{j}(D_y)V_{j}-f_{j})\in H_0^1(\mathbb R^2)
\end{equation}
for all $j, \xi$.

From~(\ref{eqA06}) and~(\ref{eqA07}), repeating the arguments of
the proof of Lemma~\ref{lA}, we obtain estimate~(\ref{eqA0'}).
\end{proof}

The following corollary of Lemma~\ref{lA0} can be proved in the
same way as Corollary~\ref{corA}.
\begin{corollary}\label{corA0}
The function $V$ constructed in Lemma~\ref{lA0} satisfies the
following inequality:
 \begin{equation}\label{eqA0'''}
 \|\mathfrak LV-f\|_{\mathcal H_0^{l,N}(K)}
 \le c\|f\|_{\hat{\mathcal S}^{l,N}(K,\gamma)}.
\end{equation}
\end{corollary}
With the help of Lemma~\ref{lA0}, we will construct the right
inverse to the operator ${\mathcal L}$, defined for compactly
supported functions $f\in \hat{\mathcal S}^{l,N}(K,\gamma)$, and
prove an analog of Theorem~\ref{thRegL1p}. However, we cannot
formally repeat the arguments of the proof of
Theorem~\ref{thRegL1p} since they are based upon the
invertibility, in weighted spaces, of the operator ${\mathcal
L}_{0}$ given by~(\ref{eqL1ap}). In this case, by
Theorem~2.1~\cite{SkDu90}, the operator ${\mathcal L}_{0}$ is not
invertible, since the line $\Im\lambda=1-l-2m$ contains the
eigenvalue $\lambda_0=i(1-l-2m)$ of $\tilde{\mathcal L}(\lambda)$.
But, as we mentioned before, the spectrum of $\tilde{\mathcal
L}(\lambda)$ is discrete; hence, there is an $a>0$ such that the
line $\Im\lambda=a+1-l-2m$ contains no eigenvalues of
$\tilde{\mathcal L}(\lambda)$, which implies that the operator
${\mathcal L}_{a}$ is invertible. In order to pass from $a>0$ to
$a=0$, we make use of the following result.
\begin{lemma}\label{lato0}
Let $W\in H_a^{l+2m,N}(K)$ for some $a>0$ and $f={\mathcal
L}_{a}W\in \mathcal H_0^{l,N}(K,\gamma)$. Suppose that the closed
strip $1-l-2m\le\Im\lambda\le a+1-l-2m$ contains only the
eigenvalue $\lambda_0=i(1-l-2m)$ of $\tilde{\mathcal L}(\lambda)$
and this eigenvalue is proper. Then we have
  \begin{equation}\label{eqato0}
   \|D^{l+2m}W\|_{H_0^{0,N}(K)}\le c \|f\|_{\mathcal H_0^{l,N}(K,\gamma)}.
  \end{equation}
\end{lemma}
Lemma~\ref{lato0} will be proved in
Sec.~\ref{subsecProofLemmalato0}. Now let us study the solvability
of problems~\eqref{eqPinK0}, \eqref{eqBinK0} and~\eqref{eqPinK},
\eqref{eqBinK} respectively.

We denote $K_j^d=K_j\cap\{y\in\mathbb R^2:\ |y|<d\}$,
$W^{k,N}(K^d)=\prod\limits_{j=1}^N W^k(K_j^d)$, and
$H_a^{k,N}(K^d)=\prod\limits_{j=1}^N H_a^k(K_j^d)$.

\begin{lemma}\label{lInvL1p0}
Let Condition~\ref{condProperEigen} hold. Then, for any
$f\in\hat{\mathcal S}^{l,N}(K,\gamma)$ with $\supp
f\subset\mathcal O_{\varepsilon}(0)$, there exists a solution $U$
to problem~\eqref{eqPinK0}, \eqref{eqBinK0} such that $U\in
W^{l+2m,N}(K^d)$ for any $d>0$ and $U$ satisfies
relations~\eqref{eqConnectUS0lN} and inequalities
 \begin{equation}\label{eqlInvL1p0'}
  \|U\|_{W^{l+2m,N}(K^d)}\le c_d\|f\|_{\hat{\mathcal S}^{l,N}(K,\gamma)},
 \end{equation}
\begin{equation}\label{eqlInvL1p0''}
 \|U\|_{H_0^{l+2m-1,N}(K^d)}\le c_d\|f\|_{\mathcal W^{l,N}(K,\gamma)}.
\end{equation}
\end{lemma}
\begin{proof}
1. Fix an $a$, $0<a<1$, such that the strip $1-l-2m<\Im\lambda\le
a+1-l-2m$ contains no eigenvalues of $\tilde{\mathcal
L}(\lambda)$. (The existence of such an $a$ follows from the
discreteness of the spectrum of~$\tilde{\mathcal L}(\lambda)$.)
From the definition of the space $\hat{\mathcal
S}^{l,N}(K,\gamma)$, it follows that, for the function $f=\{f_j,
f_{j\sigma\mu}\}$ satisfying the assumptions of the lemma,
relations~\eqref{eqSlN1} and~\eqref{eqSlN2} hold. Combining this
with Lemma~\ref{lWinHa}, we get
\begin{equation}\label{eqlInvL1p03}
  \|f\|_{\mathcal H_a^{l,N}(K,\gamma)}\le k_1\|f\|_{\mathcal W^{l,N}(K,\gamma)}.
\end{equation}
Let us consider the function $f-{\mathcal L}V$, where
$V=\hat{\mathcal A} f\in W^{l+2m,N}(K)\cap H_a^{l+2m,N}(K)$ is the
function defined in Lemma~\ref{lA0}. By virtue of
inequalities~(\ref{eqA''}) and~(\ref{eqlInvL1p03}), we have
\begin{equation}\label{eqlInvL1p04}
  \|f-{\mathcal L}V\|_{\mathcal H_a^{l,N}(K,\gamma)}\le
  k_2\|f\|_{\mathcal W^{l,N}(K,\gamma)}.
\end{equation}
Therefore, the function $f-{\mathcal L}V\in \mathcal
H_a^{l,N}(K,\gamma)$ belongs to the domain of the operator
${\mathcal L}_{a}^{-1}$. Denoting $W={\mathcal
L}_{a}^{-1}(f-{\mathcal L}V)$, we see that $U=V+W$ is a solution
to problem~(\ref{eqPinK0}), (\ref{eqBinK0}).

\smallskip

2. Let us prove~(\ref{eqlInvL1p0''}). By virtue of the boundedness
of ${\mathcal L}_{a}^{-1}$ and inequality~(\ref{eqlInvL1p04}), we
have
\begin{equation}\label{eqlInvL1p05}
  \|W\|_{H_a^{l+2m}(K)}\le
  k_3\|f\|_{\mathcal W^{l,N}(K,\gamma)}.
\end{equation}
Now estimate~(\ref{eqlInvL1p0''}) follows from
inequalities~(\ref{eqlInvL1p05}) and~(\ref{eqA''}) and the
boundedness of the embedding $H_a^{l+2m,N}(K)\subset
H_{0}^{l+2m-1,N}(K^d)$.

\smallskip

3. Let us prove~(\ref{eqlInvL1p0'}). By virtue of the boundedness
of the operator $\hat{\mathcal A}:\hat{\mathcal
S}^{l,N}(K,\gamma)\to W^{l+2m,N}(K)$ and
inequality~(\ref{eqlInvL1p05}), it suffices to estimate the
functions $D^{l+2m}W$. From Lemma~\ref{lA0}, it follows that
$f-{\mathcal L}V\in \mathcal H_0^{l,N}(K,\gamma)$ and
estimate~(\ref{eqA0'}) holds. Therefore, applying
Lemma~\ref{lato0} for the function $W={\mathcal
L}_{a}^{-1}(f-{\mathcal L}V)$ and using~(\ref{eqA0'}), we get
$$
 \|D^{l+2m}W\|_{H_0^{0,N}(K)}\le
 k_4\|f-{\mathcal L}V\|_{\mathcal H_0^{l,N}(K,\gamma)}\le
 k_5\|f\|_{\hat{\mathcal S}^{l,N}(K,\gamma)}.
$$
Noticing that $H_0^0(K_j)=L_2(K_j)$ completes the proof
of~(\ref{eqlInvL1p0'}).

4. The fulfilment of relations~(\ref{eqConnectUS0lN}) follows from
the inclusion $U=V+W\in W^{l+2m,N}(K^d)\cap H_a^{l+2m,N}(K)$ for
$a<1$ and Sobolev's embedding theorem.
\end{proof}
Now we can construct the operator $\hat{\mathfrak R}$.
\begin{theorem}\label{thRegL1p0}
Let Condition~\ref{condProperEigen} hold. Then, for any
$\varepsilon$, $0<\varepsilon<1$, there exist bounded operators
\begin{align*}
\hat{\mathfrak R}&:\{f\in \hat{\mathcal S}^{l,N}(K,\gamma): \supp
f\subset\mathcal O_{\varepsilon}(0)\}\to \{U\in W^{l+2m,N}(K):
\supp U\subset\mathcal O_{2\varepsilon_1}(0)\},\\
\hat{\mathfrak M},\hat{\mathfrak T}&:\{f\in \hat{\mathcal
S}^{l,N}(K,\gamma): \supp f\subset\mathcal O_{\varepsilon}(0)\}\to
\{f\in \hat{\mathcal S}^{l,N}(K,\gamma): \supp f\subset\mathcal
O_{2\varepsilon_1}(0)\}
\end{align*}
with\footnote{See footnote~\theqountFootnoteEps~on
p.~\pageref{pFootnoteEps}.} $\varepsilon_1={\rm
max\,}\big\{\varepsilon, \varepsilon_0/{\rm min\,}\{\chi_{j\sigma
ks},1\}\big\}$ such that $\|\hat{\mathfrak M}f\|_{\hat{\mathcal
S}^{l,N}(K,\gamma)}\le c\varepsilon_1\|f\|_{\hat{\mathcal
S}^{l,N}(K,\gamma)}$, where $c>0$ depends only on the coefficients
of the operators $\mathcal P_j(D_y)$ and $B_{j\sigma\mu ks}(D_y)$,
the operator $\hat{\mathfrak T}$ is compact, and
\begin{equation}\label{eqRegL1p0}
 \mathfrak L\hat{\mathfrak R} f=f+\hat{\mathfrak M} f+\hat{\mathfrak T} f.
\end{equation}
\end{theorem}
\begin{proof}
Let us consider a function $\psi\in C_0^\infty(\mathbb R^2)$ such
that $\psi(y)=1$ for $|y|\le \varepsilon_1={\rm
max\,}\big\{\varepsilon, \varepsilon_0/{\rm min\,}\{\chi_{j\sigma
ks},1\}\big\}$, $\supp\psi\subset\mathcal O_{2\varepsilon_1}(0)$,
and $\psi$ does not depend on polar angle $\omega$. We introduce
the operator $\hat{\mathfrak R}$ by the formula
$$
 \hat{\mathfrak R} f=\psi U\quad \big(f\in
\hat{\mathcal S}^{l,N}(K,\gamma),\ \supp f\subset\mathcal
O_{\varepsilon}(0)\big),
$$
where $U\in W^{l+2m,N}(K^{2\varepsilon_1})$ is a solution to
problem~(\ref{eqPinK0}), (\ref{eqBinK0}) with the fight-hand side
$f$ (see Lemma~\ref{lInvL1p0}).

Let us prove~\eqref{eqRegL1p0}. Relation~\eqref{eqCommutePsiG} and
Leibniz' formula imply that $\supp(\mathfrak L\hat{\mathfrak R}
f-f)\subset\mathcal O_{2\varepsilon_1}(0)$ and
\begin{equation}\label{eqReg10}
\|\mathfrak L\hat{\mathfrak R} f-f\|_{\hat{\mathcal
S}^{l,N}(K,\gamma)}\le k_1\varepsilon_1\|f\|_{\hat{\mathcal
S}^{l,N}(K,\gamma)}+k_2(\varepsilon_1)\|\psi_1
U\|_{H_0^{l+2m-1,N}(K)}),
\end{equation}
where $\psi_1\in C_0^\infty(\mathbb R^2)$ is equal to $1$ on the
support of $\psi$. From the proof of Lemma~\ref{lA0}, it follows
that the operator $f\mapsto\psi U$ acting from $\{f\in
\hat{\mathcal S}^{l,N}(K,\gamma): \supp f\subset\mathcal
O_{\varepsilon}(0)\}$ into $H_a^{l+2m,N}(K)$, $0<a<1$ is bounded.
From this and the compactness of the embedding
$$
 \{\psi_1V : V\in H_a^{l+2m,N}(K)\}\subset
H_0^{l+2m-1,N}(K),\ a<1
$$
(see Lemma~3.5~\cite{KondrTMMO67}), it follow that the operator
$f\mapsto\psi_1 U$ compactly maps $\{f\in \hat{\mathcal
S}^{l,N}(K,\gamma): \supp f\subset\mathcal O_{\varepsilon}(0)\}$
into $H_0^{l+2m-1,N}(K)$. Thus, using Lemma~\ref{lSmallComp} and
estimate~\eqref{eqReg10}, we complete the proof.
\end{proof}

Let us formulate the analog of Theorem~\ref{thRegL1p'}.

\begin{theorem}\label{thRegL1p0'}
Let Condition~\ref{condProperEigen} hold. Then, for any
$\varepsilon$, $0<\varepsilon<1$, there exist bounded operators
\begin{align*}
\hat{\mathfrak R}'&:\{f': \{0,f'\}\in \hat{\mathcal
S}^{l,N}(K,\gamma),\ \supp f'\subset\mathcal
O_{\varepsilon}(0)\}\to \{U\in
W^{l+2m,N}(K): \supp U\subset\mathcal O_{2\varepsilon}(0)\},\\
\hat{\mathfrak M}',\hat{\mathfrak T}'&:\{f': \{0,f'\}\in
\hat{\mathcal S}^{l,N}(K,\gamma),\ \supp f'\subset\mathcal
O_{\varepsilon}(0)\}\to \{f\in \hat{\mathcal S}^{l,N}(K,\gamma):
\supp f\subset\mathcal O_{2\varepsilon_2}(0)\}
\end{align*}
with $\varepsilon_2=\varepsilon/{\rm min\,}\{\chi_{j\sigma
ks},1\}$ such that $\|\hat{\mathfrak M}'f'\|_{\hat{\mathcal
S}^{l,N}(K,\gamma)}\le c\varepsilon\|\{0,f'\}\|_{\hat{\mathcal
S}^{l,N}(K,\gamma)}$, where $c>0$ depends only on the coefficients
of the operators $\mathcal P_j(D_y)$ and $B_{j\sigma\mu ks}(D_y)$,
the operator $\hat{\mathfrak T}'$ is compact, and
$$
 \mathfrak L\hat{\mathfrak R}' f'=\{0,f'\}+\hat{\mathfrak M}'f'+\hat{\mathfrak T}' f'.
$$
\end{theorem}

{\it Proof} of Theorem~\ref{thRegL1p0'} is analogous to that of
Theorem~\ref{thRegL1p'}.\qed

\subsection{Proof of Lemma~\ref{lato0}}\label{subsecProofLemmalato0} First,
we assume that $W\in\prod\limits_{j=1}^N C_0^\infty(\bar
K_j\setminus\{0\})$; then $f_j\in C_0^\infty(\bar
K_j\setminus\{0\})$ and $f_{j\sigma\mu}\in
C_0^\infty(\gamma_{j\sigma})$, where $f=\{f_j,
f_{j\sigma\mu}\}={\mathcal L} W$. We denote by $W_j(\omega, r)$
and $f_j(\omega, r)$ the functions $W_j(y)$ and $f_j(y)$
respectively, written in polar coordinates. Let $\tilde
W_j(\omega, \lambda)$, $\tilde f_j(\omega, \lambda)$, and $\tilde
f_{j\sigma\mu}(\lambda)$ be the Fourier transforms of $W_j(\omega,
e^\tau)$, $e^{2m\tau}f_j(\omega, e^\tau)$, and
$e^{m_{j\sigma\mu}\tau} f_{j\sigma\mu}(e^\tau)$ with respect to
$\tau$. Denote $\tilde f=\{\tilde f_j, \tilde f_{j\sigma\mu}\}$.
Under our assumptions, the function $\lambda\mapsto \tilde
f(\lambda)$ is analytic in the whole of the complex plane;
moreover, for $|\Im\lambda|\le\const$, this function tends to
zero, uniformly with respect to $\omega$ and $\lambda$, at a rate
higher than $|\lambda|$ to any power as $|\Re\lambda|\to\infty$ .

By virtue of Lemma~2.1~\cite{SkDu90}, there exists a
finite-meromorphic operator-valued function $\tilde{\mathcal
R}(\lambda)$ such that $\tilde{\mathcal
R}(\lambda)=\big(\tilde{\mathcal L}(\lambda)\big)^{-1}$ for any
$\lambda$ which is not an eigenvalue of $\tilde{\mathcal
L}(\lambda)$. Furthermore, if the line $\Im\lambda=a+1-l-2m$
contains no eigenvalues of $\tilde{\mathcal L}(\lambda)$, then, by
virtue of the proof of Theorem~2.1~\cite{SkDu90}, the solution $W$
is given by
\begin{equation}\label{eqato01}
  W(\omega, e^\tau)=\int\limits_{-\infty+i(a+1-l-2m)}^{+\infty+i(a+1-l-2m)}
  e^{i\lambda\tau}\tilde{\mathcal R}(\lambda){\tilde f}(\lambda)\,d\lambda.
\end{equation}

Let us consider an arbitrary $l+2m$ order derivative
$D^{l+2m}W(y)$ of the function $W$ with respect to $y_1, y_2$. Let
the operator $D^{l+2m}$ be represented in polar coordinates as
$r^{-(l+2m)}\tilde M(\omega, D_\omega, rD_r)$. After the
substitution $r=e^\tau$, the operator $D^{l+2m}$ assumes the form
$e^{-(l+2m)\tau}\tilde M(\omega, D_\omega, D_\tau)$, where
$D_\tau=-i{\partial}/{\partial\tau}$. Combining this
with~(\ref{eqato01}), we see that the function $D^{l+2m}W(y)$ can
be obtained from the function
\begin{equation}\label{eqato02}
  e^{-(l+2m)\tau}\int\limits_{-\infty+i(a+1-l-2m)}^{+\infty+i(a+1-l-2m)}
  e^{i\lambda\tau}\tilde M(\omega, D_\omega, \lambda)
  \tilde{\mathcal R}(\lambda)\tilde f(\lambda)\,d\lambda
\end{equation}
by substituting $\tau=\ln r$, followed by passing from polar
coordinates to Cartesian coordinates. Let us show that the
operator-valued function $\tilde M(\omega, D_\omega,
\lambda)\tilde{\mathcal R}(\lambda)$ is analytic near the point
$\lambda_0=i(1-l-2m)$. Since $\lambda_0$ is an eigenvalue of
$\tilde{\mathcal L}(\lambda)$, it follows from~\cite{GS} that
$$
  \tilde{\mathcal R}(\lambda)=\frac{A_{-1}}{\lambda-\lambda_0}+\Gamma(\lambda),
$$
where $\Gamma(\lambda)$ is an analytic operator-valued function
near $\lambda_0$ and the image of $A_{-1}$ coincides with the
linear span of eigenvectors corresponding to $\lambda_0$.
Therefore, for any $\tilde f\in\mathcal W^{l,N}[-b, b]$, we have
$$
 \tilde M(\omega, D_\omega, \lambda)
  \tilde{\mathcal R}(\lambda)\tilde f=
  \frac{\tilde M(\omega, D_\omega, \lambda)A_{-1}\tilde f}{\lambda-\lambda_0}+
  \tilde M(\omega, D_\omega, \lambda)\Gamma(\lambda)\tilde f.
$$
By the definition of a proper eigenvalue, the function
$r^{l+2m-1}A_{-1}\tilde f$ is a vector $Q(y)=(Q_{1}(y),\dots,
Q_{N}(y))$, where $Q_{j}(y)$ are some $l+2m-1$ order polynomials
with respect to $y_1, y_2$. Hence,
$$
 \tilde M(\omega, D_\omega, \lambda)A_{-1}\tilde f=
 r^{1-l-2m}\tilde M(\omega, D_\omega, rD_r)(r^{l+2m-1}A_{-1}\tilde f)=
 rD^{l+2m}Q(y)=0.
$$
Thus, the operator-valued function $\tilde M(\omega, D_\omega,
\lambda)\tilde{\mathcal R}(\lambda)$ is analytic near
$\lambda_0=i(1-l-2m)$ and, therefore, in the closed strip
$1-l-2m\le\Im\lambda\le a+1-l-2m$.

Furthermore, for $|\Im\lambda|\le\const$, the norm $\|\tilde
M(\omega, D_\omega, \lambda)\tilde{\mathcal
R}(\lambda)\|_{\mathcal W^{l,N}[-b, b]\to W^{0,N}(-b, b)}$ grows
at most as $|\lambda|$ to some power (see Lemma~2.1~\cite{SkDu90})
while $\|\tilde f(\lambda)\|_{\mathcal W^{l,N}[-b, b]}$ tends to
zero at a rate higher than $|\lambda|$ to any power as
$|\Re\lambda|\to\infty$. Therefore, in~(\ref{eqato02}), we can
replace the integration line $\Im\lambda=a+1-l-2m$ by the line
$\Im\lambda=1-l-2m$. Thus, the function $D^{l+2m}W(y)$ can be
obtained from the function
\begin{equation}\label{eqato03}
  e^{-(l+2m)\tau}\int\limits_{-\infty+i(1-l-2m)}^{+\infty+i(1-l-2m)}
  e^{i\lambda\tau}\tilde M(\omega, D_\omega, \lambda)
  \tilde{\mathcal R}(\lambda)\tilde f(\lambda)\,d\lambda
\end{equation}
by substituting $\tau=\ln r$, followed by passing from polar
coordinates to Cartesian coordinates. Let us estimate the norm of
$D^{l+2m}W$:
 \begin{multline}\notag
  \|D^{l+2m}W\|^2_{H_0^{0,N}(K)}=\sum\limits_j\int\limits_{K_j}|D^{l+2m}W_j|^2dy\\
  =\sum\limits_j\int\limits_{-b_j}^{b_j}d\omega\int\limits_{-\infty}^{+\infty}e^{-2(l+2m-1)\tau}
  \left|\int\limits_{-\infty+i(1-l-2m)}^{+\infty+i(1-l-2m)}
  e^{i\lambda\tau}\tilde M(\omega, D_\omega, \lambda)
  \tilde{\mathcal R}(\lambda)\tilde f(\lambda)\,d\lambda\right|^2d\tau.
 \end{multline}
Combining this with the complex analog of Parseval's equality, we
get
\begin{equation}\label{eqato04}
 \|D^{l+2m}W\|^2_{H_0^{0,N}(K)}=
 \int\limits_{-\infty+i(1-l-2m)}^{+\infty+i(1-l-2m)}
 \|\tilde M(\omega, D_\omega, \lambda)
 \tilde{\mathcal R}(\lambda)\tilde f(\lambda)\|^2_{W^{0,N}(-b, b)}d\lambda.
\end{equation}
Let us estimate the norm which is the integrand on the right-hand
side. To this end, we introduce the equivalent norms depending on
parameter $\lambda\ne0$ as follows:
 \begin{align*}
 |||\tilde U_j|||^2_{W^k(-b_j, b_j)}&=\|\tilde U_j\|^2_{W^k(-b_j, b_j)}+
 |\lambda|^{2k}\|\tilde U_j\|^2_{L_2(-b_j, b_j)},\\
 |||\tilde f|||^2_{\mathcal W^{l,N}[-b, b]}&=\sum\limits_j
 \{|||\tilde f_j|||^2_{W^l(-b_j, b_j)}+
 \sum\limits_{\sigma,\mu}|\lambda|^{2(l+2m-m_{j\sigma\mu}-1/2)}
 |\tilde f_{j\sigma\mu}|^2\}.
 \end{align*}
By virtue of the interpolation inequality
$$
 |\lambda|^{l+2m-k}\|\tilde U_j\|_{W^k(-b_j, b_j)}\le
 c_k|||\tilde U_j|||_{W^{l+2m}(-b_j, b_j)},\quad 0<k<l+2m
$$
(see.~\cite[Ch.~1]{AV}), and Lemma~2.1~\cite{SkDu90}, there exists
$C>0$ such that the following estimate holds for all
$\lambda\in\mathbb C$ satisfying $\Im\lambda=1-l-2m$ and
$|\Re\lambda|>C$:
\begin{equation}\label{eqato05}
\|\tilde M(\omega, D_\omega, \lambda)
 \tilde{\mathcal R}(\lambda)\tilde f(\lambda)\|^2_{W^{0,N}(-b, b)}\le
 k_1|||\tilde f(\lambda)|||^2_{\mathcal W^{l,N}[-b, b]}.
\end{equation}
Since the operator-valued function $\tilde M(\omega, D_\omega,
\lambda)\tilde{\mathcal R}(\lambda):\mathcal W^{l,N}[-b, b]\to
W^{0,N}(-b, b)$ is analytic on the segment $\{\lambda\in\mathbb C:
\Im\lambda=1-l-2m,\ |\Re\lambda|\le C\}$,
inequality~(\ref{eqato05}) holds on the whole line
$\Im\lambda=1-l-2m$. From~(\ref{eqato04}) and~(\ref{eqato05}), it
follows that
$$
  \|D^{l+2m}W\|^2_{H_0^{0,N}(K)}\le k_1
 \int\limits_{-\infty+i(1-l-2m)}^{+\infty+i(1-l-2m)}
 |||\tilde f(\lambda)|||^2_{\mathcal W^{l,N}[-b, b]}d\lambda.
$$
Combining this with inequalities~(1.9), (1.10)~\cite{KondrTMMO67}
yields estimate~(\ref{eqato0}). Since $C_0^\infty(\bar
K_j\setminus\{0\})$ is everywhere dense in $H_a^k(K_j)$ for any
$a$ and $k$, it follows that estimate~(\ref{eqato0}) holds for all
$W\in H_a^{l+2m,N}(K)$ and~$f\in\mathcal H_0^{l,N}(K,\gamma)$.

\section{Nonlocal Problems in Bounded Domains
in the Case where the Line $\Im\lambda=1-l-2m$ Contains no
Eigenvalues of $\tilde{\mathcal L}_p(\lambda)$}\label{sectLFred}

In this section, based on the results of
Sec.~\ref{sectLinKNoEigen}, we construct a right regularizer for
the operator $\mathbf L$ corresponding to problem~(\ref{eqPinG}),
(\ref{eqBinG}). From the existence of a right regularizer, it
follows that the image of $\mathbf L$ is closed and of finite
codimension. To prove that the kernel of $\mathbf L$ is of finite
dimension, we reduce the operator $\mathbf L$ to the operator
acting in weighted spaces and having a finite-dimensional kernel.

\smallskip

We denote $\mathbf B^k=\{\mathbf B^k_{i\mu}\}_{i,\mu}$, $k=0,
\dots, 2$; $\mathbf B=\mathbf B^0+\mathbf B^1+\mathbf B^2$,
$\mathbf C=\mathbf B^0+\mathbf B^1$. Along with the nonlocal
operator $\mathbf L=\{\mathbf P,\ \mathbf B\}$ introduced in
Sec.~\ref{sectStatement}, we also consider the bounded operators
$$
 \mathbf L^1=\{\mathbf P,\ \mathbf C\}:
W^{l+2m}(G)\to  \mathcal W^l(G,\Upsilon)\quad\text{and}\quad
\mathbf L^0=\{\mathbf P,\ \mathbf B^0\}:
 W^{l+2m}(G)\to  \mathcal W^l(G,\Upsilon).
$$

First, we will consider the operator $\mathbf L^1$ (i.e. assume
that $\mathbf B_{i\mu}^2=0$); then we will study the operator
$\mathbf L$ in the general case where $\mathbf B_{i\mu}^2\ne0$.
Throughout this section, we assume that the following condition
holds.

\begin{condition}\label{condNoEigenG}
For each orbit $\Orb_p$, $p=1, \dots, N_1$, the line
$\Im\lambda=1-l-2m$ contains no eigenvalues of the corresponding
operator $\tilde{\mathcal L}_p(\lambda)$.
\end{condition}

\subsection{Construction of the right regularizer in the case where $\mathbf
B_{i\mu}^2=0$}\label{subsectRegB^2=0}

In this subsection, we deal with the situation where $\mathbf
B_{i\mu}^2=0$, i.e., the support of nonlocal terms is concentrated
near the set $\mathcal K$.

\smallskip

For each curve $\Upsilon_i$ ($i=1, \dots, N_0$), we denote by
$g_{i1}$ and $g_{i2}$ its end points. We remind that, in some
neighborhood of the point $g_{i1}\ (g_{i2})$, the domain $G$
coincides with a plane angle while the curve $\Upsilon_i$
coincides with a segment $I_{i1}\ (I_{i2})$. Let $\tau_{i1}\
(\tau_{i2})$ be the unit vector parallel to the segment $I_{i1}\
(I_{i2})$.

We introduce the set $\mathcal S_1^l(G,\Upsilon)$ that consists of
the functions $f=\{f_0, f_{i\mu}\}\in\mathcal W^l(G,\Upsilon)$
satisfying
\begin{gather}
 D^\alpha f_0(y)=0\ (y\in\mathcal K),\quad |\alpha|\le
 l-2,\label{eqS1l1}\\
 \left.\frac{\partial^\beta f_{i\mu}}
 {\partial\tau_{i1}^\beta}\right|_{y=g_{i1}}=0,\quad
\left.\frac{\partial^\beta f_{i\mu}}
 {\partial\tau_{i2}^\beta}\right|_{y=g_{i2}}=0,\quad
 \beta\le l+2m-m_{i\mu}-2.\label{eqS1l2}
\end{gather}
From Sobolev's embedding theorem and Riesz' theorem on a general
form of linear continuous functionals in Hilbert spaces, it
follows that $\mathcal S_1^l(G,\Upsilon)$ is a closed subset of
finite codimension in the space $\mathcal W^l(G,\Upsilon)$.

\begin{lemma}\label{lRegL1}
Let Condition~\ref{condNoEigenG} hold. Then, for sufficiently
small $\varepsilon_0$, there exist a bounded operator $\mathbf
R_1:\mathcal S_1^l(G,\Upsilon)\to W^{l+2m}(G)$ and a compact
operator $\mathbf T_1: \mathcal S_1^l(G,\Upsilon)\to \mathcal
S_1^l(G,\Upsilon)$ such that
\begin{equation}\label{eqL1R1}
 \mathbf L^1\mathbf R_1=\mathbf I_1+\mathbf T_1,
\end{equation}
where $\mathbf I_1$ denotes the identity operator in $\mathcal
S_1^l(G,\Upsilon)$.
\end{lemma}
\begin{proof}
1. By virtue of Theorem~\ref{thRegL1p}, there exist bounded
operators
\begin{align*}
\mathbf R_{\mathcal K}&:\{f\in \mathcal S_1^l(G,\Upsilon): \supp
f\subset\mathcal O_{2\varepsilon_0}(\mathcal K)\}\to W^{l+2m}(G),\\
\mathbf M_{\mathcal K},\mathbf T_{\mathcal K}&:\{f\in \mathcal
S_1^l(G,\Upsilon): \supp f\subset\mathcal
O_{2\varepsilon_0}(\mathcal K)\}\to \mathcal S_1^l(G,\Upsilon)
\end{align*}
such that $\|\mathbf M_{\mathcal K}f\|_{\mathcal
W^{l}(G,\Upsilon)}\le c\varepsilon_0\|f\|_{\mathcal
W^{l}(G,\Upsilon)}$, where $c>0$ is independent of
$\varepsilon_0$, the operator $\mathbf T_{\mathcal K}$ is compact,
and
\begin{equation}\label{eqRegL1_1}
\mathbf L^1\mathbf R_{\mathcal K} f=f+\mathbf M_{\mathcal
K}f+\mathbf T_{\mathcal K} f.
\end{equation}

\smallskip

2. For each point $g\in\bar G\backslash\mathcal
O_{2\varepsilon_0}(\mathcal K)$, we consider its
$\varepsilon_0/2$-neighborhood $\mathcal O_{\varepsilon_0/2}(g)$.
All such neighborhoods, together with the set $\mathcal
O_{2\varepsilon_0}(\mathcal K)$, cover $\bar G$. Let us choose a
finite subcovering $\mathcal O_{2\varepsilon_0}(\mathcal K)$,
$\mathcal O_{\varepsilon_0/2}(g_j)$,
$j=1,\dots,J=J(\varepsilon_0)$. Let $\psi,\psi_j\in
C_0^\infty(\mathbb R^2)$, $j=1,\dots,J$, be a unity partition
corresponding to the covering $\mathcal
O_{2\varepsilon_0}(\mathcal K)$, $\mathcal
O_{\varepsilon_0/2}(g_j)$, $j=1,\dots,J$.

According to the general theory of elliptic boundary-value
problems in smooth domains (see, e.g.,~\cite{Volevich}), there
exist bounded operators
\begin{equation}\label{eqR_0j}
 \mathbf R_{0j}:\{f\in \mathcal
W^l(G,\Upsilon): \supp f\subset\mathcal
O_{\varepsilon_0/2}(g_j)\}\to \{u\in W^{l+2m}(G): \supp
u\subset\mathcal O_{\varepsilon_0}(g_j)\}
\end{equation}
and compact operators
$$
\mathbf T_{0j}:\{f\in \mathcal W^l(G,\Upsilon): \supp
f\subset\mathcal O_{\varepsilon_0/2}(g_j)\}\to \{f\in\mathcal
W^l(G,\Upsilon): \supp f\subset\mathcal O_{\varepsilon_0}(g_j)\}
$$
such that
\begin{equation}\label{eqRegL1_2}
\mathbf L^0\mathbf R_{0j} f=f+\mathbf T_{0j} f.
\end{equation}

\smallskip

3. For any $f\in \mathcal S_1^l(G,\Upsilon)$, we put $\mathbf
R_0f=\sum\limits_{j=1}^J \mathbf R_{0j}(\psi_j f)$ and
$\hat{\mathbf R}_1 f=\mathbf R_{\mathcal K}(\psi f)+\mathbf R_0f$.

Then, we have
\begin{equation}\label{eqRegL1_3}
\mathbf P\hat{\mathbf R}_1 f=\mathbf P\mathbf R_{\mathcal K}(\psi
f)+\mathbf P\mathbf R_0 f.
\end{equation}

Since $\supp\mathbf R_0f\subset\bar G\setminus\overline{\mathcal
O_{\varepsilon_0}(\mathcal K)}$, it follows from the definition of
the operator $\mathbf B^1$ that $\mathbf B^1\mathbf R_0f=0$.
Therefore,
\begin{equation}\label{eqRegL1_4}
\mathbf C\hat{\mathbf R}_1 f=\mathbf C\mathbf R_{\mathcal K}(\psi
f)+\mathbf B^0\mathbf R_0 f.
\end{equation}

From relations~\eqref{eqRegL1_3} and~\eqref{eqRegL1_4}, taking
into account~\eqref{eqRegL1_1} and~\eqref{eqRegL1_2}, we get
\begin{equation}\label{eqRegL1_5}
\mathbf L^1\hat{\mathbf R}_1 f=f+\mathbf M_{\mathcal K}(\psi
f)+\mathbf T_{\mathcal K}(\psi f)+\mathbf T_0 f,
\end{equation}
where $\mathbf T_0f=\sum\limits_{j=1}^J \mathbf T_{0j}(\psi_j f)$.

4. Let us estimate the norm of $\mathbf M_{\mathcal K}(\psi f)$:
\begin{multline}\label{eqRegL1_6}
\|\mathbf M_{\mathcal K}(\psi f)\|_{\mathcal W^l(G,\Upsilon)}\le
k_1\varepsilon_0\|\psi f\|_{\mathcal W^l(G,\Upsilon)}\\
\le k_2\varepsilon_0\|f\|_{\mathcal W^l(G,\Upsilon)}+
k_3(\varepsilon_0)\left(\|f_0\|_{W^{l-1}(G)}+\sum\limits_{i,\mu}\|\Phi_{i\mu}\|_{
W^{l+2m-m_{i\mu}-1}(G)}\right),
\end{multline}
where $\Phi_{i\mu}\in W^{l+2m-m_{i\mu}}(G)$ is an extension of
$f_{i\mu}\in W^{l+2m-m_{i\mu}-1/2}(\Upsilon_i)$ to the domain $G$
(if $l=0$, the term $\|f_0\|_{W^{l-1}(G)}$ on the right-hand side
of~\eqref{eqRegL1_6} is absent).

From~\eqref{eqRegL1_6}, the Rellich theorem, and
Lemma~\ref{lSmallComp}, it follows that
$$
\mathbf M_{\mathcal K}(\psi f)=\hat{\mathbf M}_1 f+{\mathbf T}_2f,
$$
where $\hat{\mathbf M}_1, {\mathbf T}_2:\mathcal
S_1^l(G,\Upsilon)\to \mathcal S_1^l(G,\Upsilon)$ are such that
$\|\hat{\mathbf M}_1\|\le c\varepsilon_0$ ($c>0$ is independent of
$\varepsilon_0$) and ${\mathbf T}_2$ is compact. Combining this
with relation~\eqref{eqRegL1_5}, we obtain
$$
\mathbf L^1\hat{\mathbf R}_1=\mathbf I_1+\hat{\mathbf
M}_1+\hat{\mathbf T}_1,
$$
where $\hat{\mathbf T}_1f={\mathbf T}_2f+\mathbf T_{\mathcal
K}(\psi f)+\mathbf T_0 f$.

For $\varepsilon_0\le \frac{1}{2c}$, the operator $\mathbf
I_1+\hat{\mathbf M}_1:\mathcal S_1^l(G,\Upsilon)\to \mathcal
S_1^l(G,\Upsilon)$ has a bounded inverse. Denoting $\mathbf
R_1=\hat{\mathbf R}_1(\mathbf I_1+\hat{\mathbf M}_1)^{-1}$ and
$\mathbf T_1=\hat{\mathbf T}_1(\mathbf I_1+\hat{\mathbf
M}_1)^{-1}$ yields~\eqref{eqL1R1}.
\end{proof}

\subsection{Construction of the right regularizer in the case where $\mathbf
B_{i\mu}^2\ne0$}\label{subsectRegB^2ne0}

In this subsection, we assume that $\varepsilon_0$ is fixed and
consider the operator $\mathbf L$ with $\mathbf B_{i\mu}^2\ne0$.
In other words, we suppose that there are nonlocal terms with the
support both near the set $\mathcal K$ and outside a neighborhood
of $\mathcal K$.

By virtue of Theorem~\ref{thRegL1p'}, for any sufficiently small
$\varepsilon>0$, there exist bounded operators
\begin{align*}
\mathbf R_{\mathcal K}'&:\{f': \{0,f'\}\in \mathcal
S_1^l(G,\Upsilon),\ \supp f'\subset\mathcal
O_{2\varepsilon}(\mathcal K)\}\to \{u\in W^{l+2m}(G): \supp
f'\subset\mathcal
O_{4\varepsilon}(\mathcal K)\},\\
\mathbf M_{\mathcal K}',\mathbf T_{\mathcal K}'&:\{f': \{0,f'\}\in
\mathcal S_1^l(G,\Upsilon),\ \supp f'\subset\mathcal
O_{2\varepsilon}(\mathcal K)\}\to \mathcal S_1^l(G,\Upsilon)
\end{align*}
such that $\|\mathbf M_{\mathcal K}'f'\|_{\mathcal
W^{l}(G,\Upsilon)}\le c\varepsilon\|\{0,f'\}\|_{\mathcal
W^{l}(G,\Upsilon)}$, where $c>0$ is independent of $\varepsilon$,
the operator $\mathbf T_{\mathcal K}'$ is compact, and
$$
\mathbf L^1\mathbf R_{\mathcal K}' f'=\{0,f'\}+\mathbf M_{\mathcal
K}'f'+\mathbf T_{\mathcal K}' f'.
$$
Notice that the diameter of the support of $\mathbf R_{\mathcal
K}' f'$ depends on $\varepsilon$ but is independent of
$\varepsilon_0$.

\smallskip

Similarly to the proof of Lemma~\ref{lRegL1}, we construct a
covering $\mathcal O_{2\varepsilon}(\mathcal K)$, $\mathcal
O_{\varepsilon/2}(g_j)$ ($g_j\in\partial G$, $j=1,\dots,J$,
$J=J(\varepsilon)$) of the boundary $\partial G$. Let
$\psi',\psi'_j\in C_0^\infty(\mathbb R^2)$, $j=1,\dots,J$, be a
unity partition corresponding to this covering.

According to the general theory of elliptic boundary-value
problems in smooth domains (see, e.g.,~\cite{Volevich}), there
exist bounded operators
$$
\mathbf R'_{0j}:\{f': \{0,f'\}\in \mathcal W^l(G,\Upsilon),\ \supp
f\subset\mathcal O_{\varepsilon/2}(g_j)\}\to \{u\in W^{l+2m}(G):
\supp u\subset\mathcal O_{\varepsilon}(g_j)\}
$$
and compact operators
$$
\mathbf T'_{0j}:\{f': \{0,f'\}\in \mathcal W^l(G,\Upsilon): \supp
f\subset\mathcal O_{\varepsilon/2}(g_j)\}\to \{f\in\mathcal
W^l(G,\Upsilon): \supp f\subset\mathcal O_{\varepsilon}(g_j)\}
$$
such that
\begin{equation}\notag
\mathbf L^0\mathbf R'_{0j} f'=\{0,f'\}+\mathbf T'_{0j} f'.
\end{equation}

\smallskip

For any $f'$ satisfying $\{0,f'\}\in \mathcal S_1^l(G,\Upsilon)$,
we put
\begin{equation}\label{eqR_1'}
\mathbf R'_1f'=\mathbf R_{\mathcal K}'(\psi'
f')+\sum\limits_{j=1}^J \mathbf R'_{0j}(\psi'_j f').
\end{equation}

Analogously to the proof of Lemma~\ref{lRegL1}, one can show that
\begin{equation}\label{eqL1R1'}
\mathbf L^1\mathbf R_1' f'=\{0,f'\}+\mathbf M_1'f'+\mathbf T_1'
f'.
\end{equation}
Here $\mathbf M_1',\mathbf T_1':\{f': \{0,f'\}\in \mathcal
S_1^{l}(G,\Upsilon)\}\to \mathcal S_1^{l}(G,\Upsilon)$ are bounded
operators such that $\|\mathbf M_1'f'\|_{\mathcal
W^{l}(G,\Upsilon)}\le c\varepsilon\|\{0,f'\}\|_{\mathcal
W^{l}(G,\Upsilon)}$, where $c>0$ is independent of $\varepsilon$,
and $\mathbf T_1'$ is compact.

With the help of the operators $\mathbf R_1$ (see
Lemma~\ref{lRegL1}) and $\mathbf R_1'$, we will construct a right
regularizer for the operator $\mathbf L$ with $\mathbf
B_{i\mu}^2\ne0$.

\medskip

Let us introduce the set
$$
\mathcal S^l(G,\Upsilon)=\left\{f\in\mathcal S_1^l(G,\Upsilon):
\text{the functions }\Phi=\mathbf B^2\mathbf R_1f\ \text{and }
\mathbf B^2\mathbf R_1'\Phi\ \text{satisfy
relations~\eqref{eqS1l2}}\right\}.
$$

\smallskip

From Sobolev's embedding theorem and Riesz' theorem on a general
form of linear continuous functionals in Hilbert spaces, it
follows that $\mathcal S^l(G,\Upsilon)$ is a closed subset of
finite codimension in $\mathcal W^l(G,\Upsilon)$. It is also clear
that $\mathcal S^l(G,\Upsilon)\subset\mathcal S_1^l(G,\Upsilon)$.

\begin{lemma}\label{lRegL}
Let Condition~\ref{condNoEigenG} hold. Then there exist a bounded
operator $\mathbf R:\mathcal W^l(G,\Upsilon)\to W^{l+2m}(G)$ and a
compact operator $\mathbf T: \mathcal W^l(G,\Upsilon)\to \mathcal
W^l(G,\Upsilon)$ such that
\begin{equation}\label{eqRegL}
\mathbf L\mathbf R=\mathbf I+\mathbf T,
\end{equation}
where $\mathbf I$ denotes the identity operator in $\mathcal
W^l(G,\Upsilon)$.
\end{lemma}
\begin{proof}
1. We put $\Phi=\mathbf B^2\mathbf R_1f$, where
$f=\{f_0,f'\}\in\mathcal S^l(G,\Upsilon)$. Then, by virtue of the
definition of the space $\mathcal S^l(G,\Upsilon)$, the functions
$\Phi$ and $\mathbf B^2\mathbf R_1'\Phi$ belong to the domain of
the operator $\mathbf R'_1$. Therefore, we can introduce the
bounded operator $\mathbf R_{\mathcal S}:\mathcal
S^l(G,\Upsilon)\to W^{l+2m}(G)$ by the formula
$$
 {\mathbf R}_{\mathcal S}f=\mathbf R_1f-\mathbf R_1'\Phi+\mathbf R_1'\mathbf
 B^2\mathbf R_1'\Phi.
$$

Let us show that the operator ${\mathbf R}_{\mathcal S}$ is the
right inverse to $\mathbf L$, up to the sum of small and compact
perturbations. For simplicity, we will denote by the same letter
$M$ different operators (acting in corresponding spaces) with the
norms majorized by $c\varepsilon$ and by the same letter $T$
different compact operators.

By virtue of~\eqref{eqL1R1} and \eqref{eqL1R1'}, we have
\begin{multline}\label{eqRegL_1}
\mathbf P{\mathbf R}_{\mathcal S}f=\mathbf P\mathbf R_1f-\mathbf
P\mathbf R_1'(\Phi-\mathbf B^2\mathbf R_1'\Phi)\\
=f_0+Tf_0-M(\Phi-\mathbf B^2\mathbf R_1'\Phi)-T(\Phi-\mathbf
B^2\mathbf R_1'\Phi)=f_0+Mf+Tf,
\end{multline}
\begin{multline}\label{eqRegL_2}
\mathbf C{\mathbf R}_{\mathcal S}f=\mathbf C\mathbf R_1f-\mathbf
C\mathbf R_1'\Phi+\mathbf
C\mathbf R_1'\mathbf B^2\mathbf R_1'\Phi\\
=(f'+Tf')-(\Phi+M\Phi+T\Phi)+(\mathbf B^2\mathbf R_1'\Phi+M\mathbf
B^2\mathbf R_1'\Phi+T\mathbf B^2\mathbf R_1'\Phi)=f'-\Phi+\mathbf
B^2\mathbf R_1'\Phi+Mf+Tf.
\end{multline}
Applying the operator $\mathbf B^2$ to the function ${\mathbf
R}_{\mathcal S}f$, we obtain
\begin{equation}\label{eqRegL_3}
\mathbf B^2{\mathbf R}_{\mathcal S}f=\Phi-\mathbf B^2\mathbf
R_1'\Phi+\mathbf B^2\mathbf R_1'\mathbf B^2\mathbf R_1'\Phi.
\end{equation}
Summing up equalities~\eqref{eqRegL_2} and~\eqref{eqRegL_3}, we
get
\begin{equation}\label{eqRegL_4}
\mathbf B{\mathbf R}_{\mathcal S}f=f'+Mf+Tf+\mathbf B^2\mathbf
R_1'\mathbf B^2\mathbf R_1'\Phi.
\end{equation}

Let us show that
\begin{equation}\label{eqRegL_5}
\mathbf B^2\mathbf R_1'\mathbf B^2\mathbf R_1'\Phi=0
\end{equation}
for sufficiently small
$\varepsilon=\varepsilon(\varkappa_1,\varkappa_2,\rho)$, where
$\varkappa_1,\varkappa_2,\rho$ are the constants appearing in
Condition~\ref{condSeparK23}. (Notice that $\varepsilon$ does not
depend on $\varepsilon_0$.)

By virtue of~\eqref{eqR_1'}, we have $\supp\mathbf
R_1'\Phi\subset\bar G\setminus\bar G_{4\varepsilon}$. Let
$\varepsilon$ be so small that $4\varepsilon<\rho$. Then
estimate~\eqref{eqSeparK23''} implies that $\supp\mathbf
B^2\mathbf R_1'\Phi\subset\mathcal O_{\varkappa_2}(\mathcal K)$.

Furthermore, let $\varepsilon$ be so small that
$4\varepsilon<\varkappa_1$ and
$\varkappa_2+3\varepsilon/2<\varkappa_1$. Then,
using~\eqref{eqR_1'} once more, we see that $\supp\mathbf
R_1'\mathbf B^2\mathbf R_1'\Phi\subset\mathcal
O_{\varkappa_1}(\mathcal K)$. Combining this with
inequality~\eqref{eqSeparK23'}, we get~\eqref{eqRegL_5}.

From relations~\eqref{eqRegL_1}, \eqref{eqRegL_4},
and~\eqref{eqRegL_5}, it follows that
$$
\mathbf L{\mathbf R}_{\mathcal S}=\mathbf I_{\mathcal S}+M+T,
$$
where $\mathbf I_{\mathcal S},M,T:\mathcal
S^l(G,\Upsilon)\to\mathcal W^l(G,\Upsilon)$ are bounded operators
such that $\mathbf I_{\mathcal S}f=f$, $\|M\|\le c\varepsilon$
($c>0$ is independent of $\varepsilon$), and $T$ is compact.

\smallskip

3. Since the subspace $\mathcal S^l(G,\Upsilon)$ is of finite
codimension in $\mathcal W^l(G,\Upsilon)$, the operator $\mathbf
I_{\mathcal S}$ is Fredholm. Therefore, by Theorems~16.2
and~16.4~\cite{Kr}, the operator $\mathbf I_{\mathcal S}+M+T$ is
also Fredholm, provided that $\varepsilon$ is small enough. Now
Theorem~15.2~\cite{Kr} implies the existence of a bounded operator
$\tilde{\mathbf R}$ and a compact operator $\mathbf T$ acting from
$\mathcal W^l(G,\Upsilon)$ into $\mathcal S^l(G,\Upsilon)$ and
$\mathcal W^l(G,\Upsilon)$ respectively and such that $(\mathbf
I_{\mathcal S}+M+T)\tilde{\mathbf R}=\mathbf I+\mathbf T$.
Denoting $\mathbf R=\mathbf R_{\mathcal S}\tilde{\mathbf
R}:\mathcal W^l(G,\Upsilon)\to W^{l+2m}(G)$ yields~\eqref{eqRegL}.
\end{proof}

\begin{remark}
We underline that, in the proof of Lemma~\ref{lRegL}, the numbers
$\varepsilon_0,\varkappa_1,\varkappa_2,\rho$ are fixed.
\end{remark}
\begin{remark}
The construction of the operator ${\mathbf R}$ is close to that
in~\cite{KovSk}, where the authors study nonlocal problems in
weighted spaces in the case where $\mathbf B^1=0$ (i.e., the
support of nonlocal terms does not intersect with the set
$\mathcal K$).
\end{remark}

\subsection{Fredholm solvability of nonlocal problems}\label{subsectLFred}

In this subsection, we prove the following result on the
solvability of problem~\eqref{eqPinG}, \eqref{eqBinG} in the
bounded domain in Sobolev spaces.

\begin{theorem}\label{thLFred}
Let Condition~\ref{condNoEigenG} hold; then the operator ${\bf
L}:W^{l+2m}(G)\to \mathcal W^l(G,\Upsilon)$ is Fredholm and
$\ind\mathbf L=\ind\mathbf L^1$.

Conversely, let the operator ${\bf L}:W^{l+2m}(G)\to \mathcal
W^l(G,\Upsilon)$ be Fredholm; then Condition~\ref{condNoEigenG}
holds.
\end{theorem}

We will show below that if Condition~\ref{condNoEigenG} fails,
then the image of ${\bf L}$ is not closed
(Lemma~\ref{lNon-ClosedRange}). Combining this with
Theorem~\ref{thLFred} and Theorem~7.1~\cite{Kr} implies the
following corollary.
\begin{corollary}
Condition~\ref{condNoEigenG} holds if and only if the following a
priori estimate holds:
$$
 \|u\|_{W^{l+2m}(G)}\le c(\|{\bf L}u\|_{{\mathcal W}^l(G,\Upsilon)}+\|u\|_{L_2(G)}),
$$
where $c>0$ is independent of $u$.
\end{corollary}

\subsubsection{Proof of Theorem~\ref{thLFred}. Sufficiency}
Let us show that the kernel of $\mathbf L$ is of finite dimension.
To this end, we consider problem~(\ref{eqPinG}), (\ref{eqBinG}) in
weighted spaces. We denote by $H_a^k(G)$ the completion of the set
$C_0^\infty(\bar G\setminus\mathcal K)$ with respect to the norm
$$
 \|u\|_{H_a^k(G)}=
 \left(\sum\limits_{|\alpha|\le k}\int\limits_G \rho^{2(a-k+|\alpha|)}
 |D^\alpha u|^2\right)^{1/2}.
$$
Here $k\ge 0$ is an integer; $a\in\mathbb R;$ $\rho=\rho(y)={\rm
dist}(y, {\mathcal K})$. For $k\ge1$, we denote by
$H_a^{k-1/2}(\Upsilon)$ the space of traces on a smooth curve
$\Upsilon\subset\bar G$ with the norm
$$
\|\psi\|_{H_a^{k-1/2}(\Upsilon)}=\inf\|u\|_{H_a^k(G)}\quad (u\in
H_a^k(G): u|_\Upsilon=\psi).
$$

Let us introduce the operator corresponding to
problem~(\ref{eqPinG}), (\ref{eqBinG}) in weighted spaces:
$$
 \mathbf L_a=\{{\bf P},\ {\bf B}\}: H_a^{l+2m}(G)\to  \mathcal
 H_a^l(G,\Upsilon),\quad a>l+2m-1,
$$
where $\mathcal H_a^l(G,\Upsilon)=
H_a^l(G)\times\prod\limits_{i=1}^{N_0}\prod\limits_{\mu=1}^m
H_a^{l+2m-m_{i\mu}-1/2}(\Upsilon_i)$. Notice that, by virtue
of~\eqref{eqSeparK23'} and Lemma~5.2~\cite{KovSk}, we have for
$a>l+2m-1$:
$$
\mathbf B^2_{i\mu}u\in W^{l+2m-m_{i\mu}-1/2}(\Upsilon_i)\subset
H_a^{l+2m-m_{i\mu}-1/2}(\Upsilon_i)\quad \text{for all } u\in
H_a^{l+2m}(G)\subset W^{l+2m}(G\setminus\overline{\mathcal
O_{\varkappa_1}(\mathcal
  K)}).
$$
Since the functions $\mathbf B^0_{i\mu}u$ and $\mathbf
B^1_{i\mu}u$ also belong to $H_a^{l+2m-m_{i\mu}-1/2}(\Upsilon_i)$,
it follows that the operator $\mathbf L_a$ is well defined.

Thus, the operators $\mathbf L$ and $\mathbf L_a$ correspond to
the same nonlocal problem~(\ref{eqPinG}), (\ref{eqBinG})
considered in Sobolev spaces and weighted spaces respectively.
\begin{lemma}\label{lLFiniteKernel}
The kernel of the operator $\mathbf L$ is of finite dimension.
\end{lemma}
\begin{proof}
From Lemma~2.1~\cite{SkDu90} and
Theorem\footnote{Theorem~3.2~\cite{SkDu91} is formulated for the
case where the operators $\mathbf B_{i\mu}^2$ have the same
particular form as in Example~\ref{exGeneralProblem}. However, the
proof of Theorem~3.2~\cite{SkDu91} is based on
inequalities~\eqref{eqSeparK23'} and~\eqref{eqSeparK23''} and does
not depend on any explicit form of the operators $\mathbf
B_{i\mu}^2$.}\label{pageLaFred}~3.2~\cite{SkDu91}, it follows that
the operator $\mathbf L_a$ is Fredholm for almost all $a>l+2m-1$.
Fix some $a>l+2m-1$ for which the operator $\mathbf L_a$ is
Fredholm. By Lemma~5.2~\cite{KovSk}, we have $W^{l+2m}(G)\subset
H_a^{l+2m}(G)$; therefore, $\ker\mathbf L\subset\ker\mathbf L_a$.
Since $\ker\mathbf L_a$ is of finite dimension for $a$ fixed
above, it follows that $\ker\mathbf L$ is also of finite
dimension.
\end{proof}
\begin{remark}
We underline that the kernel of the operator $\mathbf L$ is
finite-dimensional irrespective of the location of eigenvalues for
the operators $\tilde{\mathcal L}_p(\lambda)$, $p=1, \dots, N_1$.
\end{remark}

By virtue of Theorem~15.2~\cite{Kr} and Lemma~\ref{lRegL}, the
image of the operator $\mathbf L$ is closed and of finite
codimension. Combining this with Lemma~\ref{lLFiniteKernel}
implies that $\mathbf L$ is Fredholm.

\medskip

Let us show that $\ind\mathbf L=\ind\mathbf L^1$. We introduce the
operator
$$
 {\bf L}_tu=\{{\bf P}u,\ {\bf C}u+(1-t){\bf B}^2u\}.
$$
Clearly, we have ${\bf L}_0={\bf L}$, ${\bf L}_1={\bf L}^1$.

From what was proved, it follows that the operators ${\bf L}_t$
are Fredholm for all $t$. Furthermore, for any $t_0$ and $t$, the
following estimate holds:
$$
\| {\bf L}_tu- {\bf L}_{t_0}u\|_{{\cal W}^l(G,\Upsilon)}\le
k_{t_0}|t-t_0|\cdot \|u\|_{W^{l+2m}(G)},
$$
where $k_{t_0}>0$ is independent of $t$. Therefore, by
Theorem~16.2~\cite{Kr}, we have $\ind{\bf L}_t=\ind{\bf L}_{t_0}$
for all $t$ from a sufficiently small neighborhood of the point
$t_0$. Such neighborhoods cover the segment $[0, 1]$. Choosing a
finite subcovering, we obtain $\ind{\bf L}=\ind{\bf L}_0=\ind{\bf
L}_1=\ind{\bf L}^1$.

\subsubsection{Proof of Theorem~\ref{thLFred}. Necessity}

Let, to the orbit $\Orb_p$, there correspond model
problem~(\ref{eqPinK0}), (\ref{eqBinK0}) in the plane angles
$K_j=K_j^p$ with the sides $\gamma_{j\sigma}=\gamma_{j\sigma}^p$,
$j=1, \dots, N=N_{1p}$, $\sigma=1, 2$.

For any $d>0$, we consider the sets $K_j^d=K_j\cap\{y\in\mathbb
R^2:\ |y|<d\}$ and
$\gamma_{j\sigma}^d=\gamma_{j\sigma}\cap\{y\in\mathbb R^2:
|y|<d\}$ and the spaces $H_a^{l,N}(K^d)=\prod\limits_{j=1}^N
H_a^l(K_j^d)$ and
 \begin{gather*}
 W^{l,N}(K^d)=\prod\limits_{j=1}^N W^l(K_j^d),\\
 \mathcal W^{l}(K_j^d,\gamma_j^d)=
 W^l(K_j^d)\times\prod\limits_{\sigma=1,2}\prod\limits_{\mu=1}^m
W^{l+2m-m_{j\sigma\mu}-1/2}
   (\gamma_{j\sigma}^d),\\
\mathcal W^{l,N}(K^d,\gamma^d)=\prod\limits_{j=1}^{N} {\mathcal
W}^l(K_j^d,\gamma_j^d).
\end{gather*}

Put $d_1=\min\{\chi_{j\sigma ks}, 1\}/2$,
$d_2=2\max\{\chi_{j\sigma ks}, 1\}$, and
$d=d(\varepsilon)=2d_2\varepsilon$.

\begin{lemma}\label{lClosedRange}
Let the image of the operator ${\mathbf L}$ be closed. Then, for
each orbit $\Orb_p$, sufficiently small $\varepsilon$, and all
functions $U\in W^{l+2m,N}(K^{d})$, the following estimate holds:
\begin{equation}\label{eqClosedRange}
  \|U\|_{W^{l+2m,N}(K^\varepsilon)}\le c\left(\|\mathcal L_p U\|_{\mathcal
  W^{l,N}(K^{2\varepsilon},\gamma^{2\varepsilon})}+\sum\limits_{j=1}^N
  \|\mathcal P_j(D_y)U_j\|_{W^l(K_j^d)}+\|U\|_{W^{l+2m-1,N}(K^{d})}\right).
\end{equation}
\end{lemma}
\begin{proof}
1. Since the image of $\mathbf L$ is closed, it follows from
Lemma~\ref{lLFiniteKernel}, compactness of the embedding
$W^{l+2m}(G)\subset W^{l+2m-1}(G)$, and Theorem~7.1~\cite{Kr} that
\begin{equation}\label{eqClosedRange1}
  \|u\|_{W^{l+2m}(G)}\le c(\|\mathbf L u\|_{\mathcal
  W^l(G,\Upsilon)}+\|u\|_{W^{l+2m-1}(G)}).
\end{equation}
Let us substitute functions $u\in W^{l+2m}(G)$ such that $\supp
u\subset\bigcup\limits_{j=1}^{N_{1p}}\mathcal
O_{2\varepsilon}(g_j^p)$,
$2\varepsilon<\min\{\varepsilon_0,\varkappa_1\}$,
into~(\ref{eqClosedRange1}). By virtue of~\eqref{eqSeparK23'}, for
such functions, we have $\mathbf B^2u=0$. Therefore, using
Lemma~3.2~\cite[Ch.~2]{LM}, for a sufficiently small
$\varepsilon$, we obtain the following estimate:
\begin{equation}\label{eqClosedRange2}
  \|U\|_{W^{l+2m,N}(K)}\le c(\|\mathcal L_p U\|_{\mathcal
  W^{l,N}(K,\gamma)}+\|U\|_{W^{l+2m-1,N}(K)}),
\end{equation}
which holds for all $U\in W^{l+2m,N}(K)$ with $\supp
U\subset\mathcal O_{2\varepsilon}(0)$.

\smallskip

2. Now let us refuse the assumption $\supp U\subset\mathcal
O_{2\varepsilon}(0)$ and show that, for any $U\in
W^{l+2m,N}(K^{d})$, estimate~\eqref{eqClosedRange} holds.

We introduce a function $\psi\in C_0^\infty(\mathbb R^2)$ such
that $\psi(y)=1$ for $|y|\le \varepsilon$,
$\supp\psi\subset\mathcal O_{2\varepsilon}(0)$, and $\psi$ does
not depend on polar angle $\omega$. Using
inequality~(\ref{eqClosedRange2}) and Leibniz' formula, for $U\in
W^{l+2m,N}(K^{d})$, we obtain
\begin{multline}\label{eqClosedRange4}
  \|U\|_{W^{l+2m,N}(K^\varepsilon)}\le \|\psi U\|_{W^{l+2m,N}(K)}\le k_1
  (\|\mathcal L_p (\psi U)\|_{\mathcal
  W^{l,N}(K,\gamma)}
  +\|\psi U\|_{W^{l+2m-1,N}(K)})\\
  \le
  k_2(\|\psi \mathcal L_p U\|_{\mathcal
  W^{l,N}(K,\gamma)}+\sum\limits_{j,\sigma,\mu}\,\sum\limits_{(k,s)\ne(j,0)}\|J_{j\sigma\mu
ks}\|_{W^{l+2m-m_{j\sigma\mu}-1/2}(\gamma_{j\sigma})}
  +\|U\|_{W^{l+2m-1,N}(K^{2\varepsilon})}),
\end{multline}
where
 \begin{equation}\notag 
J_{j\sigma\mu ks}=\big(\psi(\mathcal G_{j\sigma
ks}y)-\psi(y)\big)\big(B_{j\sigma\mu ks}(D_y)U_k\big)\big(\mathcal
G_{j\sigma ks}y\big)\big|_{\gamma_{j\sigma}}.
 \end{equation}
Let us estimate the norm of $J_{j\sigma\mu ks}$. Notice that, for
$(k,s)\ne(j,0)$, the operator $\mathcal G_{j\sigma ks}$ maps the
ray $\gamma_{j\sigma}$ onto the ray
$$
\{y\in{\mathbb R}^2: r>0,\ \omega=(-1)^\sigma
b_{j}+\omega_{j\sigma ks}\},
$$
which is strictly inside the angle $K_k$. Therefore, there exists
a function $\xi_{j\sigma ks}\in C_0^\infty(-b_k,b_k)$ equal to $1$
at the point $\omega=(-1)^\sigma b_{j}+\omega_{j\sigma ks}$.

Furthermore, the support of the function $\psi(y)-\psi(\mathcal
G_{j\sigma ks}^{-1}y)$ is contained in the set
$\{d_1\varepsilon<|y|<d_2\varepsilon\}$. Therefore, there exists a
function $\psi_1\in C_0^\infty(K_k)$ equal to $1$ on the support
of the function $\xi(\omega)\big(\psi(y)-\psi(\mathcal G_{j\sigma
ks}^{-1}y)\big)$ and such that
$\supp\psi_1\subset\{d_1\varepsilon<|y|<d_2\varepsilon\}$. Then,
similarly to~\eqref{eqRegL1p2}, we obtain
\begin{equation}\notag
\|J_{j\sigma\mu
ks}\|_{W^{l+2m-m_{j\sigma\mu}-1/2}(\gamma_{j\sigma})} \le k_3
\|\psi_1U_k\|_{W^{l+2m}(K_k)}.
\end{equation}
Let us estimate the norm on the right-hand side of this inequality
by using Theorem~5.1~\cite[Ch.~2]{LM} and Leibniz's formula. Then,
taking into account that $\psi_1$ is compactly supported and
vanishes both near the origin and near the sides of $K_k$, we get
\begin{equation}\label{eqClosedRange6}
\|J_{j\sigma\mu
ks}\|_{W^{l+2m-m_{j\sigma\mu}-1/2}(\gamma_{j\sigma})} \le
k_4(\|\mathcal
P_k(D_y)U_k\|_{W^l(\{d_1\varepsilon/2<|y|<2d_2\varepsilon\})}+
\|U_k\|_{W^{l+2m-1}(\{d_1\varepsilon/2<|y|<2d_2\varepsilon\})}).
\end{equation}

Now estimate~(\ref{eqClosedRange}) follows
from~(\ref{eqClosedRange4}) and~(\ref{eqClosedRange6}).
\end{proof}

\begin{lemma}\label{lNon-ClosedRange}
Let the line $\Im\lambda=1-l-2m$ contain an eigenvalue of the
operator $\tilde{\mathcal L}_p(\lambda)$ for some $p$. Then the
image of the operator $\mathbf L$ is not closed.
\end{lemma}
\begin{proof}
1. Suppose that the image of $\mathbf L$ is closed. The following
two cases are possible: either {\rm (a)}~the line
$\Im\lambda=1-l-2m$ contains an improper eigenvalue or {\rm(b)}
the line $\Im\lambda=1-l-2m$ contains only the eigenvalue
$\lambda_0=i(1-l-2m)$, which is proper (see
Definitions~\ref{defRegEigVal} and~\ref{defRegEigValImproper}).

\smallskip

2. First we assume that there is an improper eigenvalue
$\lambda=\lambda_0$. Let us show that, in this case,
estimate~(\ref{eqClosedRange}) does not hold. Denote by
$\varphi^{(0)}(\omega), \dots, \varphi^{(\varkappa-1)}(\omega)$ an
eigenvector and associate vectors (the {\it Jordan chain} of
length $\varkappa\ge1$) corresponding to the eigenvalue
$\lambda_0$ (see~\cite{GS}). According to
Remark~2.1~\cite{GurAsympAngle}, the vectors
$\varphi^{(k)}(\omega)$ belong to $W^{l+2m,N}(-b, b)$, and, by
Lemma~2.1~\cite{GurAsympAngle}, we have
\begin{equation}\label{eqNon-ClosedRange8}
 \mathcal L_p V^k=0,
\end{equation}
where $V^k=r^{i\lambda_0}\sum\limits_{s=0}^k\frac{1}{s!}(i\ln
r)^k\varphi^{(k-s)}(\omega)$,~~$k=0, \dots, \varkappa-1$. Since
$\lambda_0$ is not a proper eigenvalue, it follows that, for some
$k\ge0$, the function $V^k(y)$ is not a vector-polynomial. For
simplicity, we suppose that
$V^0=r^{i\lambda_0}\varphi^{(0)}(\omega)$ is not a
vector-polynomial (the case where $k>0$ can be considered
analogously).

We introduce the sequence $U^\delta=r^\delta V^0/\|r^\delta
V^0\|_{W^{l+2m,N}(K^\varepsilon)}$. For any $\delta>0$, the
denominator is finite, but $\|r^\delta
V^0\|_{W^{l+2m,N}(K^\varepsilon)}\to\infty$ as $\delta\to 0$ since
$V^0$ is not a vector-polynomial. However, $\|r^\delta
V^0\|_{W^{l+2m-1,N}(K^{d})}\le c$, where $c>0$ is independent of
$\delta\ge0$; therefore,
\begin{equation}\label{eqNon-ClosedRange9}
 \|U^\delta\|_{W^{l+2m-1,N}(K^{d})}\to 0\quad \mbox{as } \delta\to 0.
\end{equation}
Moreover, relation~(\ref{eqNon-ClosedRange8}) implies
$$
 \mathcal P_j(D_y)U^\delta=\frac{r^\delta\mathcal P_j(D_y)V^0+
 \sum\limits_{|\alpha|+|\beta|=2m,|\alpha|\ge1}
 p_{j\alpha\beta}D^\alpha r^\delta\cdot D^\beta V_j^0}{\|r^\delta
V^0\|_{W^{l+2m,N}(K^\varepsilon)}}
=\frac{\sum\limits_{|\alpha|+|\beta|=2m,|\alpha|\ge1}
 p_{j\alpha\beta}D^\alpha r^\delta\cdot D^\beta V_j^0}{\|r^\delta
V^0\|_{W^{l+2m,N}(K^\varepsilon)}},
$$
where $p_{j\alpha\beta}$ are some complex constants. Hence,
$|D^\xi\mathcal P_j(D_y)U^\delta|\le c_{j\xi}\delta
r^{l-1-|\xi|+\delta}/\|r^\delta V^0\|_{W^{l+2m,N}(K^\varepsilon)}$
($|\xi|\le l$), which implies that
\begin{equation}\label{eqNon-ClosedRange10}
 \|\mathcal P_j(D_y)U^\delta\|_{W^{l}(K_j^{d})}\to 0\quad \mbox{as } \delta\to 0.
\end{equation}
Similarly, by using~(\ref{eqNon-ClosedRange8}), one can prove that
\begin{equation}\label{eqNon-ClosedRange11}
 \|\mathcal B_{j\sigma\mu}(D_y)U^\delta|_{\gamma_{j\sigma}}\|_{
 W^{l+2m-m_{j\sigma\mu}-1/2}(\gamma_{j\sigma}^{2\varepsilon})}\to 0\quad \mbox{as } \delta\to
 0.
\end{equation}
(Here one must additionally estimate the expression
$$
 \frac{\sum\limits_{(k,s)\ne(j,0)}\|(\chi_{j\sigma ks}^\delta-1)r^\delta
 (B_{j\sigma\mu ks}(y, D_y)V^0)(\mathcal G_{j\sigma
ks}y)|_{\gamma_{j\sigma}}\|_{W^{l+2m-m_{j\sigma\mu}-1/2}(\gamma_{j\sigma}^{2\varepsilon})}}{\|r^\delta
V^0\|_{W^{l+2m,N}(K^\varepsilon)}},
$$
which also tends to zero as $\delta\to 0$ because of the
inequality $|\chi_{j\sigma ks}^\delta-1|\le k_6\delta$.)

However,
assertions~(\ref{eqNon-ClosedRange9})--(\ref{eqNon-ClosedRange11})
contradict estimate~(\ref{eqClosedRange}), since
$\|U^\delta\|_{W^{l+2m,N}(K^\varepsilon)}=1$.

\smallskip

3. It remains to consider the case where the line
$\Im\lambda=1-l-2m$ contains only the eigenvalue
$\lambda_0=i(1-l-2m)$ of $\tilde{\mathcal L}_p(\lambda)$ and it is
proper. In this case, we cannot repeat the arguments above, since
$V^0$ is a vector-polynomial and the norm $\|r^\delta
V^0\|_{W^{l+2m,N}(K^\varepsilon)}$ is uniformly bounded as
$\delta\to 0$.

Let us make use of the results of Sec.~\ref{sectLinKProperEigen}.
By Lemma~\ref{lS0NotClosed}, there exists a sequence $f^\delta\in
\hat{\mathcal S}^{l,N}(K,\gamma)$, $\delta>0$, such that $\supp
f^\delta\subset\mathcal O_{\varepsilon}(0)$ and $f^\delta$
converges in $\mathcal W^{l,N}(K,\gamma)$ to
$f^0\notin\hat{\mathcal S}^{l,N}(K,\gamma)$ as $\delta\to0$. By
Lemma~\ref{lInvL1p0}, for each $f^\delta$, there exists a function
$U^\delta\in W^{l+2m,N}(K^d)$ such that
\begin{equation}\label{eqNon-ClosedRange12}
 \mathcal L_p U^\delta = f^\delta,
\end{equation}
\begin{equation}\label{eqNon-ClosedRange13}
\|U^\delta\|_{W^{l+2m-1,N}(K^{d})}\le c\|f^\delta\|_{\mathcal
W^{l,N}(K,\gamma)}
\end{equation}
($c>0$ is independent of $\delta$) and $U^\delta$ satisfies
relations~(\ref{eqConnectUS0lN}). From
inequalities~(\ref{eqClosedRange})
and~(\ref{eqNon-ClosedRange13}),
relation~(\ref{eqNon-ClosedRange12}), and convergence of
$f^\delta$ in $\mathcal W^{l,N}(K,\gamma)$, it follows that the
sequence $U^\delta$ is fundamental in
$W^{l+2m,N}(K^{\varepsilon})$. Therefore, $U^\delta$ converges in
$W^{l+2m,N}(K^{\varepsilon})$ to some function $U$ as $\delta\to
0$. Moreover, the limit function $U$ also satisfies
relations~(\ref{eqConnectUS0lN}), and, by virtue of the
boundedness of the operator $\mathcal L_p:
W^{l+2m,N}(K^{\varepsilon})\to \mathcal
W^{l,N}(K^{2d_1\varepsilon},\gamma^{2d_1\varepsilon})$, the
following relation holds:
$$
 \mathcal L_p U = f^0\quad \mbox{for } y\in\mathcal
O_{2d_1\varepsilon}(0).
$$
We consider a function $\psi\in C_0^\infty(\mathbb R^2)$ such that
$\psi(y)=1$ for $|y|\le d_1^2\varepsilon$ and
$\supp\psi\subset\mathcal O_{2d_1^2\varepsilon}(0)$. Clearly,
$\psi U\in W^{l+2m,N}(K)$, $\psi U$ satisfies
relations~(\ref{eqConnectUS0lN}), and $\supp\mathcal L_p (\psi
U)\subset\mathcal O_{2d_1\varepsilon}(0)$. Therefore,
$$
 \mathcal L_p (\psi U) = \psi f^0+\hat f,
$$
where $\hat f\in\mathcal W^{l,N}(K,\gamma)$ and the support of
$\hat f$ is compact and does not intersect with the origin. Hence,
the function $\psi f^0+\hat f$, together with $f^0$, does not
belong to $\hat{\mathcal S}^{l,N}(K,\gamma)$. This contradicts
Lemma~\ref{lConnectUS0lN}.
\end{proof}

Now the necessity of the conditions in Theorem~\ref{thLFred}
follows from Lemma~\ref{lNon-ClosedRange}.

\section{Asymptotics of Solutions to Nonlocal Problems
in Sobolev Spaces}\label{sectAsymp}

\subsection{Smoothness of solutions outside the set $\mathcal K$}

In this subsection, we prove the following result on smoothness of
solutions to problem~(\ref{eqPinG}), (\ref{eqBinG}) inside the
domain and near a smooth part of the boundary.

\begin{lemma}\label{lSmoothOutsideK}
Let $u\in W^{l+2m}(G)$ be a solution to problem~\eqref{eqPinG},
\eqref{eqBinG}. Suppose that the right-hand side $f=\{f_0,
f_{i\mu}\}$ belongs to $\mathcal W^{l_1}(G,\Upsilon)$ with $l_1>l$
and Condition~\ref{condSeparK23} holds for $l_1$ substituted for
$l$. Then
\begin{equation}\label{eqSmoothOutsideK}
u\in W^{l_1+2m}\bigl(G\setminus\overline{\mathcal
O_\delta(\mathcal K)}\bigr)\quad \text{for any } \delta>0.
\end{equation}
\end{lemma}
\begin{proof}
1. We denote by $W^l_\loc(G)$ the space of distributions $v$ in
$G$ such that $\psi v\in W^l(G)$ for all $\psi\in C_0^\infty(G)$.
By Theorem~3.2~\cite[Ch.~2]{LM}, we have
\begin{equation}\label{eqSmoothLocG}
 u\in W^{l_1+2m}_\loc(G).
\end{equation}
Combining this with estimate~\eqref{eqSeparK23''} implies that
\begin{equation}\label{eqSmoothOutsideK_2}
\mathbf B_{i\mu}^2 u\in
W^{l_1+2m-m_{i\mu}-1/2}(\Upsilon_i\setminus\overline{\mathcal
O_{\varkappa_2}(\mathcal K)}).
\end{equation}

We fix an arbitrary point
$g\in\Upsilon_i\setminus\overline{\mathcal
O_{\varkappa_2}(\mathcal K)}$ and choose a $\delta>0$ such that
\begin{equation}\label{eqSmoothOutsideK_3}
\begin{aligned}
 &\overline{\mathcal
O_{\delta}(g)\cap\Upsilon_i}\subset\Upsilon_i\setminus\overline{\mathcal
O_{\varkappa_2}(\mathcal K)}\quad \text{and}\\
&\text{if } g\in\mathcal O_{\varepsilon_0}(\mathcal K),\
\text{then } \overline{\Omega_{is}\bigl(\mathcal
O_{\delta}(g)\cap\mathcal O_{\varepsilon_0}(\mathcal
K)\bigr)}\subset G.
\end{aligned}
\end{equation}
Then, in the neighborhood $\mathcal O_{\delta}(g)$, the function
$u$ is a solution to the following problem:
\begin{align}
\mathbf P(y, D_y)u&=f_0(y)\quad (y\in\mathcal
O_{\delta}(g)\cap G),\label{eqSmoothOutsideK_4}\\
B_{i\mu0}(y, D_y)u&=f^2_{i\mu}(y)\quad (y\in\mathcal
O_{\delta}(g)\cap \Upsilon_i;\ \mu=1, \dots,
m),\label{eqSmoothOutsideK_5}
\end{align}
where
$f^2_{i\mu}(y)=f_{i\mu}(y)-\sum\limits_{s=1}^{S_i}\bigl(B_{i\mu
s}(y, D_y)(\zeta u)\bigr)\bigl(\Omega_{is}(y)\bigr)-\mathbf
B_{i\mu}^2u(y)$ ($y\in\mathcal O_{\delta}(g)\cap \Upsilon_i$).
From relations~\eqref{eqSmoothLocG}, \eqref{eqSmoothOutsideK_2},
and~\eqref{eqSmoothOutsideK_3}, it follows that $f^2_{i\mu}\in
W^{l_1+2m-m_{i\mu}-1/2}(\mathcal O_{\delta}(g)\cap \Upsilon_i)$.

Applying Theorem\footnote{In Theorem~5.1~\cite[Ch.~2]{LM}, the
operators $B_{i\mu0}(y, D_y)$ are additionally supposed to be
normal on $\Upsilon_i$ while their orders are supposed not to
exceed $2m-1$. However, it is easy to check that
Theorem~5.1~\cite[Ch.~2]{LM} remains true without these
assumptions (see~\cite[Ch.~2, \S~8.3]{LM}).}~5.1~\cite[Ch.~2]{LM}
to problem~(\ref{eqSmoothOutsideK_4}), (\ref{eqSmoothOutsideK_5}),
we see that
\begin{equation}\label{eqSmoothOutsideK_6}
 u\in W^{l_1+2m}(\mathcal O_{\delta/2}(g)\cap G).
\end{equation}

By using the unity partition method, we derive
from~(\ref{eqSmoothLocG}) and~(\ref{eqSmoothOutsideK_6}) that
\begin{equation}\label{eqSmoothOutsideK_7}
u\in W^{l_1+2m}\bigl(G\setminus\overline{\mathcal
O_{\varkappa_1}(\mathcal K)}\bigr).
\end{equation}

\smallskip

2. From inclusion~\eqref{eqSmoothOutsideK_7} and
inequality~\eqref{eqSeparK23'}, it follows that
\begin{equation}\label{eqSmoothOutsideK_8}
\mathbf B_{i\mu}^2 u\in W^{l_1+2m-m_{i\mu}-1/2}(\Upsilon_i).
\end{equation}

Taking into account~\eqref{eqSmoothOutsideK_8}, we can repeat the
arguments of item~1 of this proof for an arbitrary point
$g\in\Upsilon_i$ and for $\delta>0$ such that
\begin{equation}\notag
\begin{aligned}
&\overline{\mathcal
O_{\delta}(g)\cap\Upsilon_i}\subset\Upsilon_i\quad \text{and}\\
&\text{if } g\in\mathcal O_{\varepsilon_0}(\mathcal K),\
\text{then } \overline{\Omega_{is}\bigl(\mathcal
O_{\delta}(g)\cap\mathcal O_{\varepsilon_0}(\mathcal
K)\bigr)}\subset G.
\end{aligned}
\end{equation}
As a result, we arrive at~\eqref{eqSmoothOutsideK_6}, which is now
true for an arbitrary $g\in\Upsilon_i$. Combining this
with~(\ref{eqSmoothLocG}) and using the unity partition method, we
obtain~\eqref{eqSmoothOutsideK}.
\end{proof}

\subsection{Asymptotics of solutions near the set $\mathcal K$}

In this subsection, we obtain an asymptotic formula for the
solution $u$ near an arbitrary orbit $\Orb_p\subset\mathcal K$,
provided that the line $\Im\lambda=1-l_1-2m$ contains no
eigenvalues of the operator $\tilde{\mathcal L}_p(\lambda)$.

\smallskip

Thus, let us fix some orbit $\Orb_p\subset\mathcal K$, which
consists of the points $g_j^p$, $j=1, \dots, N=N_{1p}$, and choose
an $\varepsilon>0$ such that $\mathcal
O_{\varepsilon}(g_j^p)\subset\mathcal V(g_j^p)$. Then, in the
neighborhood $\bigcup\limits_{j=1}^{N}\mathcal
O_{\varepsilon}(g_j^p)$ of the orbit $\Orb_p$, the function $u$ is
a solution to the following problem:
\begin{gather}
  {\bf P}(y, D_y) u_j=f(y) \quad (y\in\mathcal O_{\varepsilon}(g_j)\cap
  G),\label{eqPinOrbp}\\
\begin{aligned}
 B_{i\mu 0}(y, D_y)u_j(y)|_{\Upsilon_i}+\sum\limits_{s=1}^{S_i}
 \bigl(B_{i\mu s}(y, D_y)(\zeta u_k)\bigr)\bigl(\Omega_{is}(y)\bigr)|_{\Upsilon_i}=
   f'_{i\mu}(y)& \\
   \big(y\in \mathcal O_{\varepsilon}(g_j^p)\cap\Upsilon_i;\
    i\in\{1\le i\le N_0: g_j\in\bar\Upsilon_i\};\
   j=1, \dots, N;\ \mu=1, \dots, m\big)&.
\end{aligned}\label{eqBinOrbp}
\end{gather}
Here $u_1(y), \dots, u_N(y)$ denote the same as in
Sec.~\ref{subsectStatementNearK} and
$f'_{i\mu}(y)=f_{i\mu}(y)-\mathbf B_{i\mu}^2u(y)$ for $y\in
\mathcal O_{\varepsilon}(g_j^p)\cap\Upsilon_i$. By virtue
of~\eqref{eqSmoothOutsideK_8}, we have $f'_{i\mu}\in
W^{l_1+2m-m_{i\mu}-1/2}(\mathcal
O_{\varepsilon}(g_j^p)\cap\Upsilon_i)$.

Let $y\mapsto y'(g^{p}_j)$ be the argument transformation
described in Sec.~\ref{subsectStatement}. Analogously to
Sec.~\ref{subsectStatementNearK}, we introduce the function
$U_j(y')=u_j(y(y'))$ and denote $y'$ again by $y$. For $p$ being
fixed above, we put $b_j=b^{p}_j$, $K_j=K^p_j$, and
$\gamma_{j\sigma}=\{y\in{\mathbb R}^2: r>0,\ \omega=(-1)^\sigma
b_{j}\}$ ($\sigma=1, 2$). Then problem~(\ref{eqPinOrbp}),
(\ref{eqBinOrbp}) assumes the following form (cf.~\eqref{eqPinK},
\eqref{eqBinK}):
\begin{gather}
{\bf P}_{j}(y, D_y)U_j=f_{j}(y) \quad (y\in K_j^{\varepsilon}),\label{eqPinKEpsilon}\\
{\bf B}_{j\sigma\mu}(y,
D_y)U|_{\gamma_{j\sigma}^{\varepsilon}}\equiv \sum\limits_{k,s}
       (B_{j\sigma\mu ks}(y, D_y)U_k)({\mathcal G}_{j\sigma ks}y)|_{\gamma_{j\sigma}^{\varepsilon}}
    =f_{j\sigma\mu}(y) \quad (y\in\gamma_{j\sigma}^{\varepsilon}).\label{eqBinKEpsilon}
\end{gather}
Here $f=\{f_j, f_{j\sigma\mu}\}\in\mathcal
W^{l_1,N}(K^{\varepsilon},\gamma^{\varepsilon})$ and $U\in
W^{l+2m,N}(K^{{d}})$, where $d=\varepsilon\max\{\chi_{j\sigma
ks},1\}$ ($\chi_{j\sigma ks}$ are the argument expansion
coefficients corresponding to the orbit $\Orb_p$).

To obtain the asymptotics of the solution $u$ to
problem~(\ref{eqPinG}), (\ref{eqBinG}) near the orbit $\Orb_p$, we
preliminarily investigate the asymptotics of the solution $U$ to
problem~(\ref{eqPinKEpsilon}), (\ref{eqBinKEpsilon}) near the
origin.

\medskip

By virtue of Lemma~4.11~\cite{KondrTMMO67}, the function $U_j\in
W^{l+2m}(K_j^{{d}})$ can be represented in the following form:
\begin{equation}\label{eqQ+U1}
  U_j(y)=Q_j(y)+U_j^1(y),
\end{equation}
where $Q_j(y)$ is an $l+2m-2$ order polynomial and $U_j^1\in
W^{l+2m}(K_j^{{d}})\cap H_a^{l+2m}(K_j^{{d}})$ for any $a>0$.
Putting $Q=(Q_1, \dots, Q_N)$, we see that the function
$U^1=(U^1_1, \dots, U^1_N)$ is a solution to the problem
\begin{align}
  {\bf P}_{j}(y, D_y)U_j^1=f_{j}(y)-{\bf P}_{j}(y, D_y)Q_j(y)&\equiv
  f_j^1(y)\quad
    (y\in K_j^{\varepsilon}),\label{eqPinKEpsilon1}\\
 {\bf B}_{j\sigma\mu}(y, D_y)U^1|_{\gamma_{j\sigma}^{\varepsilon}}
    =f_{j\sigma\mu}(y)-
    {\bf B}_{j\sigma\mu}(y, D_y)Q|_{\gamma_{j\sigma}^{\varepsilon}}&\equiv
    f_{j\sigma\mu}^1(y)\quad
      (y\in\gamma_{j\sigma}^{\varepsilon}),\label{eqBinKEpsilon1}
\end{align}
where $f^1=\{f_j^1, f_{j\sigma\mu}^1\}\in\mathcal
W^{l_1,N}(K^{\varepsilon},\gamma^{\varepsilon})$.

Now, using Lemma~4.11~\cite{KondrTMMO67}, we represent the
function $f_j^1\in W^{l_1}(K_j^{\varepsilon})$ in the following
form:
\begin{equation}\label{eqP+f2}
 f_j^1(y)=P_j(y)+f_j^2(y),
\end{equation}
where $P_j(y)$ is an $l_1-2$ order polynomial (if $l_1\ge2$) and
$f_j^2\in W^{l_1}(K_j^{\varepsilon})\cap
H_a^{l_1}(K_j^{\varepsilon})$ for any $a>0$. If $l_1\le1$, then we
put $P_j(y)\equiv 0$, in which case $f_j^1=f_j^2\in
H_a^{l_1}(K_j^{\varepsilon})$ by Lemma~\ref{lWinHa}. We notice
that, on the one hand, the inclusion $U_j^1\in
H_a^{l+2m}(K_j^{{d}})$ implies the inclusion $f_j^1\in
H_a^{l}(K_j^{\varepsilon})$ and, on the other hand, $f_j^2\in
H_a^{l_1}(K_j^{\varepsilon})\subset H_a^{l}(K_j^{\varepsilon})$.
Thus, $P_j\in H_a^{l}(K_j^{\varepsilon})$ and, therefore, $P_j$
consists of monomials of order $\ge$$l-1$.

Similarly, we have
\begin{equation}\label{eqP+f2'}
 f_{j\sigma\mu}^1(y)=P_{j\sigma\mu}(y)+f^2_{j\sigma\mu}(y),
\end{equation}
where $P_{j\sigma\mu}(y)$ is an $l_1+2m-m_{j\sigma\mu}-2$ order
polynomial (if $l_1+2m-m_{j\sigma\mu}\ge 2$) consisting of
monomials of order $\ge$$l+2m-m_{j\sigma\mu}-1$ and
$f^2_{j\sigma\mu}\in
W^{l_1+2m-m_{j\sigma\mu}-1/2}(\gamma_{j\sigma}^{\varepsilon})\cap
H_a^{l_1+2m-m_{j\sigma\mu}-1/2}(\gamma_{j\sigma}^{\varepsilon})$
for any $a>0$. If $l_1+2m-m_{j\sigma\mu}\le1$, then
$P_{j\sigma\mu}(y)\equiv 0$.

By virtue of Lemma\footnote{In Lemma~3.1~\cite{SkMs86} (as well as
in Lemma~3.2~\cite{SkMs86} below), nonlocal terms are supposed to
contain only rotation operators but not expansion ones. However,
the corresponding results also remain true in our case
(see~\cite{GurAsympAngle}).}~3.1~\cite{SkMs86}, there exist the
functions
$$
 W_j=\sum\limits_{s=l+2m-1}^{l_1+2m-1}\sum\limits_{q=0}^{q_j}
     r^s(i\ln r)^q\varphi_{jsq}(\omega)\in
     H_a^{l+2m}(K_j^{\varepsilon}),\quad a>0,
$$
with $\varphi_{jsk}\in C^\infty\big([-b_j, b_j]\big)$ such that
the vector $W=(W_1, \dots, W_N)$ satisfies the following relation:
\begin{equation}\label{eqPW-P}
  \mathbf P_j(y, D_y)W_j-P_j\in H_0^{l_1}(K_j^{\varepsilon}),
\end{equation}
\begin{equation}\label{eqBW-P}
  \mathbf B_{j\sigma\mu}(y, D_y)W-P_{j\sigma\mu}\in
H_0^{l_1+2m-m_{j\sigma\mu}-1/2}(\gamma_{j\sigma}^{\varepsilon}).
\end{equation}

Further, since
$$
f_j^2\in W^{l_1}(K_j^{\varepsilon})\cap
H_a^{l_1}(K_j^{\varepsilon}),\quad f_{j\sigma\mu}^2\in
W^{l_1+2m-m_{j\sigma\mu}-1/2}(\gamma_{j\sigma}^{\varepsilon})\cap
H_a^{l_1+2m-m_{j\sigma\mu}-1/2}(\gamma_{j\sigma}^{\varepsilon})
$$
for all $a>0$, it follows that the functions $f_j^2$ and
$f_{j\sigma\mu}^2$ satisfy the relations
\begin{align}
 D^\alpha f_j^2|_{y=0}&=0,\quad |\alpha|\le
 l_1-2,\label{eqSl_1N1}\\
 \left.\frac{\partial^\beta f_{j\sigma\mu}^2}
 {\partial\tau_{j\sigma}^\beta}\right|_{y=0}&=0,\quad \beta\le
 l_1+2m-m_{j\sigma\mu}-2.\label{eqSl_1N2}
\end{align}
Therefore, by virtue of Lemma~\ref{lA} and Corollary~\ref{corA},
there exist functions $V_j\in W^{l_1+2m}(K_j^{{d}})\cap
H_a^{l_1+2m}(K_j^{{d}})$ for any $a>0$ such that the vector
$V=(V_1, \dots, V_N)$ satisfies the relation
\begin{equation}\label{eqPV-f2}
  \mathbf P_j(y, D_y)V_j-f_j^2\in H_0^{l_1}(K_j^{\varepsilon}),
\end{equation}
\begin{equation}\label{eqBV-f2}
  \mathbf B_{j\sigma\mu}(y, D_y)V-f^2_{j\sigma\mu}\in
H_0^{l_1+2m-m_{j\sigma\mu}-1/2}(\gamma_{j\sigma}^{\varepsilon}).
\end{equation}

From~(\ref{eqPinKEpsilon1})--(\ref{eqBV-f2}), it follows that the
vector
\begin{equation}\label{eqU2}
  U^2=U^1-V-W\in H_a^{l+2m,N}(K^{{d}})
\end{equation}
is a solution to the problem
\begin{equation}\label{eqPU2}
  \mathbf P_j(y, D_y)U_j^2=(P_j-\mathbf P_j(y, D_y)W_j)+
  (f_j^2-\mathbf P_j(y, D_y)V_j)\in H_0^{l_1}(K_j^{\varepsilon}),
\end{equation}
\begin{multline}\label{eqBU2}
\mathbf B_{j\sigma\mu}(y,
D_y)U^2|_{\gamma_{j\sigma}^{\varepsilon}}=(P_{j\sigma\mu}-\mathbf
B_{j\sigma\mu}(y, D_y)W)|_{\gamma_{j\sigma}^{\varepsilon}}\\
+(f_{j\sigma\mu}^2-\mathbf B_{j\sigma\mu}(y,
D_y)V|_{\gamma_{j\sigma}^{\varepsilon}})\in
H_0^{l_1+2m-m_{j\sigma\mu}-1/2}(\gamma_{j\sigma}^{\varepsilon}).
\end{multline}
Let us choose $a>0$ so small that the strip $1-l-2m<\Im\lambda\le
a+1-l-2m$ contains no eigenvalues of $\tilde{\mathcal
L}_p(\lambda)$, which is possible since the spectrum of
$\tilde{\mathcal L}_p(\lambda)$ is discrete. Then
equalities~(\ref{eqPU2}) and~(\ref{eqBU2}) and
Lemma~3.2~\cite{SkMs86} imply the following asymptotic formula for
$U_j^2\in H_a^{l+2m}(K_j^{{d}})$:
\begin{equation}\label{eqAsympU2}
 U_j^2=\sum\limits_{1-l_1-2m<\Im\lambda_n\le 1-l-2m}\,\sum\limits_{s,q}
 r^{i\lambda_n+s}(i\ln r)^q\psi_{jnsq}(\omega)+U_j^3\quad (y\in
 K_j^{\varepsilon}).
\end{equation}
Here $U_j^3\in H_0^{l_1+2m}(K_j^{\varepsilon})$, $\lambda_n$ are
eigenvalues of $\tilde{\mathcal L}_p(\lambda)$, $\psi_{jnsq}\in
C^\infty\big([-b_j, b_j]\big)$, $s=0, \dots, s_n$,
$s_n=[l_1+2m-1+\Im\lambda_n]$, and $q=0, \dots, q_{jn}\ge0$.

Formula~(\ref{eqAsympU2}) and relations~(\ref{eqQ+U1})
and~(\ref{eqU2}) imply
\begin{equation}\label{eqAsympUinK}
 U_j=\sum\limits_{n}\sum\limits_{s,q}
 r^{i\lambda_n+s}(i\ln r)^q\psi_{jnsq}(\omega)
 +\sum\limits_{s,q}
     r^s(i\ln r)^q\varphi_{jsk}(\omega)+U_j^4\quad (y\in
 K_j^{\varepsilon}),
\end{equation}
where $U_j^4=U_j^3+V_j+Q_j\in W^{l_1+2m}(K^{{d}})$.

Notice that the function
 \begin{equation}\notag
  J_j=\sum\limits_{\Im\lambda_n=1-l-2m}\sum\limits_{q=0}^{q_{jn}}
 r^{i\lambda_n}(i\ln r)^q\psi_{j,n0q}(\omega)
 +\sum\limits_{q=0}^{q_j}
     r^{l+2m-1}(i\ln r)^q\varphi_{j,l+2m-1,k}(\omega)
 \end{equation}
is a homogeneous $l+2m-1$ order polynomial with respect to $y_1,
y_2$ (otherwise, Lemma~4.20~\cite{KondrTMMO67} implies that
$J_j\notin W^{l+2m}(K_j^{{d}})$, while the other terms
in~(\ref{eqAsympUinK}) do belong to $W^{l+2m}(K_j^{{d}}$)). Thus,
we finally obtain
\begin{multline}\label{eqAsympUinKFin}
 U_j=\sum\limits_{1-l_1-2m<\Im\lambda_n\le 1-l-2m}\,\sum\limits_{s,q}
 r^{i\lambda_n+s}(i\ln r)^q\psi_{jnsq}(\omega)\\
 +\sum\limits_{s=l+2m}^{l_1+2m-1}\sum\limits_{q=0}^{q_j}
     r^s(i\ln r)^q\varphi_{jsk}(\omega)+U_j^5\quad (y\in
 K_j^{\varepsilon}),
\end{multline}
where $U_j^5=U_j^4+J_j\in W^{l_1+2m}(K^{{d}})$ and, in the first
interior sum, indices range as follows: $s=1, \dots, s_n$ if
$\Im\lambda_n=1-l-2m$, $s=0, \dots, s_n$ if $\Im\lambda_n<1-l-2m$,
and $q=0, \dots, q_{jn}\ge0$.

\medskip

From Lemma~\ref{lSmoothOutsideK} and
representation~\eqref{eqAsympUinKFin}, we derive the main result
of this section.

\begin{theorem}\label{thAsympinG}
Let $u\in W^{l+2m}(G)$ be a solution to problem~\eqref{eqPinG},
\eqref{eqBinG}, and let the conditions of
Lemma~\ref{lSmoothOutsideK} hold. Then the solution $u$ satisfies
relations~\eqref{eqSmoothOutsideK}. If we additionally suppose
that the line $\Im\lambda=1-l_1-2m$ contains no eigenvalues of the
operator $\tilde{\mathcal L}_p(\lambda)$ for some $p\in\{1, \dots,
N_1\}$, then, in the neighborhood $\mathcal
O_{\varepsilon}(g_j^p)$ ($j=1, \dots, N_{1p}$), the following
representation holds:
\begin{equation}\label{eqAsympinG}
 u=\sum\limits_{n}\sum\limits_{s,q}
 r^{i\lambda_n+s}(i\ln r)^q\psi'_{jnsq}(\omega)
 +\sum\limits_{s,q}
     r^s(i\ln r)^q\varphi'_{jsk}(\omega)+u'\quad \big(y\in
 \mathcal O_{\varepsilon}(g_j^p)\cap G\big).
\end{equation}
Here $(\omega, r)$ are polar coordinates with pole at $g_j^p$,
$\psi'_{jnsq}$ and $\varphi'_{jsk}$ are the functions infinitely
differentiable with respect to $\omega$ and turning into the
functions $\psi_{jnsq}$ and $\varphi_{jsk}$ respectively after the
change of variables $y\mapsto y'(g_j^p)$; further, $u'\in
W^{l_1+2m}\big(\mathcal O_{\varepsilon}(g_j^p)\cap G\big)$ and all
the indices in~\eqref{eqAsympinG} range as
in~\eqref{eqAsympUinKFin}.
\end{theorem}

Theorem~\ref{thAsympinG}, in particular, means that if $u\in
W^{l+2m}(G)$ is a solution to problem~(\ref{eqPinG}),
(\ref{eqBinG}) with the right-hand side $f=\{f_0,
f_{i\mu}\}\in\mathcal W^{l_1}(G,\Upsilon)$ ($l_1>l$) and the
closed strip $1-l_1-2m\le\Im\lambda\le 1-l-2m$ contains no
eigenvalues of the operators $\tilde{\mathcal L}_p(\lambda)$,
$p=1, \dots, N_1$, then $u\in W^{l_1+2m}(G)$.

\section{Nonlocal Problems in Bounded Domains in Weighted Spaces with Small
 Weight Exponents}\label{sectLFredH_a}

\subsection{Formulation of the main result}

In Sec.~\ref{subsectLFred}, we introduced the operator
\begin{equation}\label{eqLa>l+2m-1}
 \mathbf L_a=\{{\bf P},\ {\bf B}\}: H_a^{l+2m}(G)\to  \mathcal
 H_a^l(G,\Upsilon),\quad a>l+2m-1.
\end{equation}
As we mentioned in the proof of Lemma~\ref{lLFiniteKernel}, the
operator $\mathbf L_a$ is Fredholm for almost all $a>l+2m-1$ due
to Lemma~2.1~\cite{SkDu90} and Theorem~3.2~\cite{SkDu91}.

In this subsection, we consider problem~(\ref{eqPinG}),
(\ref{eqBinG}) in weighted spaces with the weight exponent $a>0$.
In that case, as before, $\mathbf B^2_{i\mu}u\in
W^{l+2m-m_{i\mu}-1/2}(\Upsilon_i)$ for all $u\in
H_a^{l+2m}(G)\subset W^{l+2m}(G\setminus\overline{\mathcal
O_{\varkappa_1}(\mathcal K)})$. However, the function $\mathbf
B^2_{i\mu}u$ may now not belong to the space
$H_a^{l+2m-m_{i\mu}-1/2}(\Upsilon_i)$, which implies that the
operator $\mathbf L_a$ given by~\eqref{eqLa>l+2m-1} is not well
defined.

\smallskip

Let us introduce the set
$$
 S_a^{l+2m}(G)=\left\{u\in H_a^{l+2m}(G): \text{the functions }\mathbf
 B^2_{i\mu}u\ \text{satisfy conditions~\eqref{eqS1l2}}\right\}.
$$
Using inequality~\eqref{eqSeparK23'}, for $\beta\le
l+2m-m_{i\mu}-2$ and $k=1,2$, we obtain
$$
\|\mathbf
 B^2_{i\mu}u\|_{W^{l+2m-m_{i\mu}-1/2}(\Upsilon_i)} \le k_1
\|u\|_{W^{l+2m}(G\setminus\overline{\mathcal
O_{\varkappa_1}(\mathcal
  K)})}\le k_2 \|u\|_{H_a^{l+2m}(G)}.
$$
Combining this with Sobolev's embedding theorem and Riesz' theorem
on a general form of linear continuous functionals in Hilbert
spaces, we see that $S_a^{l+2m}(G)$ is a closed subspace of finite
codimension in $H_a^{l+2m}(G)$.

On the other hand, by Lemma~\ref{lWinHa}, for any $u\in
S_a^{l+2m}(G)$, we have $\mathbf B^2_{i\mu}u\in
H_a^{l+2m-m_{i\mu}-1/2}(\Upsilon_i)$ for $a>0$. Since the
functions $\mathbf B^0_{i\mu}u$ and $\mathbf B^1_{i\mu}u$ also
belong to $H_a^{l+2m-m_{i\mu}-1/2}(\Upsilon_i)$ for all
$a\in\mathbb R$ and $u\in S_a^{l+2m}(G)$ (and even for $u\in
H_a^{l+2m}(G)$), it follows that
$$
 \{{\bf P}u,\ {\bf B}u\}\in \mathcal
 H_a^l(G,\Upsilon)\quad \text{for all } u\in S_a^{l+2m}(G),\ a>0.
$$

Thus, there exists a \textit{finite-dimensional} space $\mathcal
R_a^l(G,\Upsilon)$ (which is embedded into
$H_a^l(G)\times\prod\limits_{i,\mu}
H_{a'}^{l+2m-m_{i\mu}-1/2}(\Upsilon_i)$, $a'>l+2m-1$) such that
$\mathcal H_a^l(G,\Upsilon)\cap\mathcal R_a^l(G,\Upsilon)=\{0\}$
and
$$
\{{\bf P}u,\ {\bf B}u\}\in \mathcal
H_a^l(G,\Upsilon)\oplus\mathcal R_a^l(G,\Upsilon) \quad \text{for
all } u\in H_a^{l+2m}(G),\ a>0.
$$
Therefore, we can define the bounded operator
\begin{equation}\notag
 \mathbf L_a=\{{\bf P},\ {\bf B}\}: H_a^{l+2m}(G)\to  \mathcal
 H_a^l(G,\Upsilon)\oplus\mathcal R_a^l(G,\Upsilon),\quad a>0.
\end{equation}
Clearly, one can put $\mathcal R_a^l(G,\Upsilon)=\{0\}$ if
$a>l+2m-1$.

\begin{theorem}\label{thLFredH_a}
Let $a>0$ and the line $\Im\lambda=a+1-l-2m$ contain no
eigenvalues of the operators $\tilde{\mathcal L}_p(\lambda)$,
$p=1,\dots,N_1$; then the operator ${\bf L}_a:H_a^{l+2m}(G)\to
\mathcal H_a^l(G,\Upsilon)\oplus\mathcal R_a^l(G,\Upsilon)$ is
Fredholm.

Conversely, let the operator ${\bf L}_a:H_a^{l+2m}(G)\to \mathcal
H_a^l(G,\Upsilon)\oplus\mathcal R_a^l(G,\Upsilon)$ be Fredholm;
then the line $\Im\lambda=a+1-l-2m$ contains no eigenvalues of
either of the operators $\tilde{\mathcal L}_p(\lambda)$,
$p=1,\dots,N_1$.
\end{theorem}

Notice that if $f\in \mathcal H_a^l(G,\Upsilon)$, then $
\|f\|_{\mathcal H_a^l(G,\Upsilon)\oplus\mathcal
R_a^l(G,\Upsilon)}=\|f\|_{\mathcal H_a^l(G,\Upsilon)}.$ Combining
this with Theorem~\ref{thLFredH_a} and Riesz' theorem on a general
form of linear continuous functionals in Hilbert spaces, we obtain
the following result.

\begin{corollary}\label{corLFredH_a}
Let $a>0$ and the line $\Im\lambda=a+1-l-2m$ contain no
eigenvalues of the operators $\tilde{\mathcal L}_p(\lambda)$,
$p=1,\dots,N_1$. Then there exist functions $f^q\in \mathcal
H_a^l(G,\Upsilon)$, $q=1,\dots,q_1$, such that if the right-hand
side $f$ of problem~\eqref{eqPinG}, \eqref{eqBinG} belongs to
$\mathcal H_a^l(G,\Upsilon)$ and
$$
 (f,\ f^q)_{\mathcal H_a^l(G,\Upsilon)}=0,\quad q=1,\dots,q_1,
$$
then problem~\eqref{eqPinG}, \eqref{eqBinG} admits a solution
$u\in H_a^{l+2m}(G)$.
\end{corollary}

Corollary~\ref{corLFredH_a} shows: \textit{in spite of the fact
that the inclusion $u\in H_a^{l+2m}(G)$ for $0<a\le l+2m-1$ does
not, generally speaking, implies the inclusion $\mathbf L_a u\in
\mathcal H_a^l(G,\Upsilon)$, if we impose on the right-hand side
$f\in \mathcal H_a^l(G,\Upsilon)$ finitely many orthogonality
conditions, then problem~\eqref{eqPinG}, \eqref{eqBinG} yet admits
a solution $u\in H_a^{l+2m}(G)$.}

\subsection{Proof of the main result}

\subsubsection{Proof of Theorem~\ref{thLFredH_a}. Sufficiency}

\begin{lemma}\label{lLFiniteKernelH_a}
The kernel of the operator $\mathbf L_a$ is of finite dimension.
\end{lemma}
\begin{proof}
Notice that $H_a^{l+2m}(G)\subset H_{a'}^{l+2m}(G)$ for $a\le a'$.
Thus, the lemma can be proved in the same way as
Lemma~\ref{lLFiniteKernel}.
\end{proof}

Let us proceed to construct the right regularizer for the operator
$\mathbf L_a$.

As we mentioned before, the functions $\mathbf B^0_{i\mu}u$ and
$\mathbf B^1_{i\mu}u$ belong to
$H_a^{l+2m-m_{i\mu}-1/2}(\Upsilon_i)$ for all $u\in H_a^{l+2m}(G)$
and $a\in\mathbb R$. Therefore, we can introduce the bounded
operator
$$
 \mathbf L_a^1=\{{\bf P},\ {\bf C}\}: H_a^{l+2m}(G)\to  \mathcal
 H_a^l(G,\Upsilon).
$$

In~\cite[\S~3]{SkDu91}, it is proved that there exist a bounded
operator $\mathbf R_{a,1}:\mathcal H_a^l(G,\Upsilon)\to
H_a^{l+2m}(G)$ and a compact operator $\mathbf T_{a,1}: \mathcal
H_a^l(G,\Upsilon)\to \mathcal H_a^l(G,\Upsilon)$ such that
\begin{equation}\label{eqL1R1H_a}
 \mathbf L_a^1\mathbf R_{a,1}=\mathbf I_a+\mathbf T_a,
\end{equation}
where $\mathbf I_a$ denotes the identity operator in $\mathcal
H_a^l(G,\Upsilon)$.

\smallskip

Further, from Theorem~\ref{thRegL1p'H_a}, it follows that, for any
sufficiently small $\varepsilon>0$, there exist bounded operators
\begin{align*}
{\mathbf R}_{a,\mathcal K}'&:\{f': \{0,f'\}\in {\mathcal
H}_a^l(G,\Upsilon),\ \supp f'\subset\mathcal
O_{2\varepsilon}(\mathcal K)\}\to \{u\in H_a^{l+2m}(G): \supp
f'\subset\mathcal
O_{4\varepsilon}(\mathcal K)\},\\
{\mathbf M}_{a,\mathcal K}',{\mathbf T}_{a,\mathcal K}'&:\{f':
\{0,f'\}\in {\mathcal H}_a^l(G,\Upsilon),\ \supp f'\subset\mathcal
O_{2\varepsilon}(\mathcal K)\}\to {\mathcal H}_a^l(G,\Upsilon)
\end{align*}
such that $\|{\mathbf M}_{a,\mathcal K}'f'\|_{{\mathcal
H}_a^{l}(G,\Upsilon)}\le c\varepsilon\|\{0,f'\}\|_{{\mathcal
H}_a^{l}(G,\Upsilon)}$, where $c>0$ is independent of
$\varepsilon$, the operator ${\mathbf T}_{a,\mathcal K}'$ is
compact, and
$$
{\mathbf L}_a^1{\mathbf R}_{a,\mathcal K}' f'=\{0,f'\}+{\mathbf
M}_{a,\mathcal K}'f'+{\mathbf T}_{a,\mathcal K}' f'.
$$

\smallskip

For any $f'$ such that $\{0,f'\}\in {\mathcal H}_a^l(G,\Upsilon)$,
we put
\begin{equation}\label{eqR_1'H_a}
{\mathbf R}'_{a,1}f'={\mathbf R}_{a,\mathcal K}'(\psi'
f')+\sum\limits_{j=1}^J \mathbf R'_{0j}(\psi'_j f'),
\end{equation}
where the functions $\psi',\psi_j'$ and the operators $\mathbf
R'_{0j}$ are the same as in Sec.~\ref{subsectRegB^2ne0}.

By using Theorem~\ref{thRegL1p'H_a}, one can easily show that
\begin{equation}\label{eqL1R1'inSH_a}
{\mathbf L}_a^1{\mathbf R}_{a,1}' f'=\{0,f'\}+{\mathbf
M}_{a,1}'f'+{\mathbf T}_{a,1}' f'.
\end{equation}
Here ${\mathbf M}_{a,1}',{\mathbf T}_{a,1}':\{f': \{0,f'\}\in
{\mathcal H}_a^{l}(G,\Upsilon)\}\to {\mathcal
H}_a^{l}(G,\Upsilon)$ are bounded operators such that $\|{\mathbf
M}_{a,1}'f'\|_{{\mathcal H}_a^{l}(G,\Upsilon)}\le
c\varepsilon\|\{0,f'\}\|_{{\mathcal H}_a^{l}(G,\Upsilon)}$, where
$c>0$ is independent of $\varepsilon$, and the operator ${\mathbf
T}_{a,1}'$ is compact.

With the help of the operators ${\mathbf R}_{a,1}$ and ${\mathbf
R}_{a,1}'$, we will construct the right regularizer for
problem~(\ref{eqPinG}), (\ref{eqBinG}) with $\mathbf
B_{i\mu}^2\ne0$ in weighted spaces.

\medskip

For $a>0$, we introduce the set
$$
\mathcal S_a^l(G,\Upsilon)=\left\{f\in\mathcal H_a^l(G,\Upsilon):
\text{the functions }\Phi=\mathbf B^2\mathbf R_{a,1}f\ \text{and }
\mathbf B^2\mathbf R_{a,1}'\Phi\ \text{satisfy
conditions~\eqref{eqS1l2}} \right\}.
$$

\smallskip

First, let us show that $\mathcal S_a^l(G,\Upsilon)$ is a closed
subspace of finite codimension in $\mathcal H_a^l(G,\Upsilon)$.
Indeed, by using inequality~\eqref{eqSeparK23'}, for $\beta\le
l+2m-m_{i\mu}-2$ and $k=1,2$, we obtain
\begin{equation}\label{eqH_a1_1}
\|\Phi_{i\mu}\|_{W^{l+2m-m_{i\mu}-1/2}(\Upsilon_i)} \le k_1
\|\mathbf R_{a,1}f\|_{W^{l+2m}(G\setminus\overline{\mathcal
O_{\varkappa_1}(\mathcal
  K)})}
  \le k_2 \|\mathbf
R_{a,1}f\|_{H_a^{l+2m}(G)}\le k_3\|f\|_{\mathcal
H_a^l(G,\Upsilon)}.
\end{equation}
Since the function $\Phi_{i\mu}$ satisfies
conditions~\eqref{eqS1l2}, it follows from~\eqref{eqH_a1_1} and
Lemma~\ref{lWinHa} that $\Phi_{i\mu}\in
H_a^{l+2m-m_{i\mu}-1/2}(\Upsilon_i)$ and
\begin{equation}\label{eqH_a1_1'}
\|\Phi_{i\mu}\|_{H_a^{l+2m-m_{i\mu}-1/2}(\Upsilon_i)}\le
k_4\|f\|_{\mathcal H_a^l(G,\Upsilon)}.
\end{equation}
Therefore, the expression $\mathbf B^2\mathbf R_{a,1}'\Phi$ is
well defined. Similarly, using~\eqref{eqH_a1_1'}
and~\eqref{eqS1l2}, we get
\begin{equation}\label{eqH_a1_1''}
\|[\mathbf B^2\mathbf
R_{a,1}'\Phi]_{i\mu}\|_{W^{l+2m-m_{i\mu}-1/2}(\Upsilon_i)}\le
k_5\|f\|_{\mathcal H_a^l(G,\Upsilon)},
\end{equation}
and
\begin{equation}\label{eqH_a1_1'''}
\|[\mathbf B^2\mathbf
R_{a,1}'\Phi]_{i\mu}\|_{H_a^{l+2m-m_{i\mu}-1/2}(\Upsilon_i)}\le
k_6\|f\|_{\mathcal H_a^l(G,\Upsilon)},
\end{equation}
where $[\,\cdot\,]_{i\mu}$ stands for the corresponding vector's
component.

From~\eqref{eqH_a1_1}, \eqref{eqH_a1_1''}, Sobolev's embedding
theorem, and Riesz' theorem on a general form of linear continuous
functionals in Hilbert spaces, it follows that $\mathcal
S_a^l(G,\Upsilon)$ is a closed subspace of finite codimension in
$\mathcal H_a^l(G,\Upsilon)$. Hence,
\begin{equation}\label{eqHR=SHatR}
\mathcal H_a^l(G,\Upsilon)\oplus\mathcal
R_a^l(G,\Upsilon)=\mathcal S_a^l(G,\Upsilon)\oplus \hat{\mathcal
R}_a^l(G,\Upsilon),
\end{equation}
where $\hat{\mathcal R}_a^l(G,\Upsilon)$ is some
\textit{finite-dimensional} space. Now we are in a position to
prove the following result.

\begin{lemma}\label{lRegLH_a}
Let $a>0$ and the line $\Im\lambda=a+1-l-2m$ contain no
eigenvalues of the operators $\tilde{\mathcal L}_p(\lambda)$,
$p=1,\dots,N_1$. Then there exist a bounded operator $\mathbf R_a:
\mathcal H_a^l(G,\Upsilon)\oplus\mathcal R_a^l(G,\Upsilon)\to
H_a^{l+2m}(G)$ and a compact operator $\mathbf T_a: \mathcal
H_a^l(G,\Upsilon)\oplus\mathcal R_a^l(G,\Upsilon)\to \mathcal
H_a^l(G,\Upsilon)\oplus\mathcal R_a^l(G,\Upsilon)$ such that
\begin{equation}\label{eqRegLH_a}
\mathbf L\mathbf R=\hat{\mathbf I}_a+\mathbf T_a,
\end{equation}
where $\hat{\mathbf I}_a$ denotes the identity operator in
$\mathcal H_a^l(G,\Upsilon)\oplus\mathcal R_a^l(G,\Upsilon)$.
\end{lemma}
\begin{proof}
1. Put $\Phi=\mathbf B^2\mathbf R_{a,1}f$, where $f\in\mathcal
S_a^l(G,\Upsilon)$. Then, from~\eqref{eqH_a1_1'}
and~\eqref{eqH_a1_1'''}, it follows that the functions
$\{0,\Phi\}$ and $\{0,\mathbf B^2\mathbf R_{a,1}'\Phi\}$ belong to
$\mathcal H_a^l(G,\Upsilon)$. Therefore, the functions $\Phi$ and
$\mathbf B^2\mathbf R_{a,1}'\Phi$ belong to the domain of the
operator $\mathbf R'_{a,1}$, and we may introduce the bounded
operator $\mathbf R_{a,\mathcal S}:\mathcal S_a^l(G,\Upsilon)\to
H_a^{l+2m}(G)$ by the formula
$$
{\mathbf R}_{a,\mathcal S}f=\mathbf R_{a,1}f-\mathbf
R_{a,1}'\Phi+\mathbf R_{a,1}'\mathbf B^2\mathbf R_{a,1}'\Phi.
$$

Analogously to the proof of Lemma~\ref{lRegL}, using
equalities~\eqref{eqL1R1H_a} and~\eqref{eqL1R1'inSH_a}, one can
show that
$$
\mathbf L_a{\mathbf R}_{a,\mathcal S}=\mathbf I_{a,\mathcal
S}+M+T,
$$
where $\mathbf I_{a,\mathcal S},M,T:\mathcal
S_a^l(G,\Upsilon)\to\mathcal \mathcal
H_a^l(G,\Upsilon)\oplus\mathcal R_a^l(G,\Upsilon)$ are bounded
operators such that $\mathbf I_{a,\mathcal S}f=f$, $\|M\|\le
c\varepsilon$ ($c>0$ is independent of $\varepsilon$), and $T$ is
compact.

\smallskip

2. Due to~\eqref{eqHR=SHatR}, the subspace $\mathcal
S_a^l(G,\Upsilon)$ is of finite codimension in $\mathcal
H_a^l(G,\Upsilon)\oplus\mathcal R_a^l(G,\Upsilon)$. Therefore the
operator $\mathbf I_{a,\mathcal S}$ is Fredholm. By Theorem~16.2
and~16.4~\cite{Kr}, the operator $\mathbf I_{a,\mathcal S}+M+T$ is
also Frdholm, provided that $\varepsilon$ is small enough. From
Theorem~15.2~\cite{Kr}, it follows that there exist a bounded
operator $\tilde{\mathbf R}_a$ and a compact operator $\mathbf
T_a$ acting from $\mathcal H_a^l(G,\Upsilon)\oplus\mathcal
R_a^l(G,\Upsilon)$ into $\mathcal S_a^l(G,\Upsilon)$ and $\mathcal
H_a^l(G,\Upsilon)\oplus\mathcal R_a^l(G,\Upsilon)$ respectively
and such that $(\mathbf I_{a,\mathcal S}+M+T)\tilde{\mathbf
R}_a=\hat{\mathbf I}_a+\mathbf T_a$. Denoting $\mathbf R_a=\mathbf
R_{a,\mathcal S}\tilde{\mathbf R}_a:\mathcal
H_a^l(G,\Upsilon)\oplus\mathcal R_a^l(G,\Upsilon)\to
H_a^{l+2m}(G)$ yields~\eqref{eqRegLH_a}.
\end{proof}

By virtue of Theorem~15.2~\cite{Kr} and Lemma~\ref{lRegLH_a}, the
image of the operator $\mathbf L_a$, $a>0$, is closed and of
finite codimension. Combining this with
Lemma~\ref{lLFiniteKernelH_a} proves the sufficiency of the
conditions in Theorem~\ref{thLFredH_a}.

\subsubsection{Proof of Theorem~\ref{thLFredH_a}. Necessity}

\begin{lemma}\label{lNon-ClosedRangeH_a}
Let $a>0$ and the line $\Im\lambda=a+1-l-2m$ contain an eigenvalue
of the operator $\tilde{\mathcal L}_p(\lambda)$ for some $p$. Then
the image of $\mathbf L_a$ is not closed.
\end{lemma}
\begin{proof}
1. Let, to the orbit $\Orb_p$, there correspond model
problem~(\ref{eqPinK0}), (\ref{eqBinK0}) in the angles $K_j=K_j^p$
with the sides $\gamma_{j\sigma}=\gamma_{j\sigma}^p$, $j=1, \dots,
N=N_{1p}$, $\sigma=1, 2$.

For any $d>0$, we introduce the spaces
 \begin{gather*}
 \mathcal H_a^{l}(K_j^d,\gamma_j^d)=
 H_a^l(K_j^d)\times\prod\limits_{\sigma=1,2}\prod\limits_{\mu=1}^m
H_a^{l+2m-m_{j\sigma\mu}-1/2}
   (\gamma_{j\sigma}^d),\\
\mathcal H_a^{l,N}(K^d,\gamma^d)=\prod\limits_{j=1}^{N} {\mathcal
H}_a^l(K_j^d,\gamma_j^d).
\end{gather*}

Put $d_1=\min\{\chi_{j\sigma ks}, 1\}/2$,
$d_2=2\max\{\chi_{j\sigma ks}, 1\}$,
$d=d(\varepsilon)=2d_2\varepsilon$.

Assume that the image of $\mathbf L_a$ is closed. Then, similarly
to the proof of Lemma~\ref{lClosedRange}, by using
Lemma~\ref{lLFiniteKernelH_a}, compactness of the embedding
$H_a^{l+2m}(G)\subset H_a^{l+2m-1}(G)$, and Theorem~7.1~\cite{Kr},
one can show that
\begin{equation}\label{eqClosedRangeH_a}
  \|U\|_{H_a^{l+2m,N}(K^\varepsilon)}\le c\left(\|\mathcal L_p U\|_{\mathcal
  H_a^{l,N}(K^{2\varepsilon},\gamma^{2\varepsilon})}+\sum\limits_{j=1}^N
  \|\mathcal P_j(D_y)U_j\|_{H_a^l(K_j^d)}+\|U\|_{H_a^{l+2m-1,N}(K^{d})}\right)
\end{equation}
for all $U\in H_a^{l+2m,N}(K^{d})$ and sufficiently small
$\varepsilon$.

\smallskip

2. Let $\lambda_0$ be an eigenvalue of ~$\tilde{\mathcal
L}_p(\lambda)$, lying on the line $\Im\lambda=a+1-l-2m$, and
$\varphi^{(0)}(\omega)$ an eigenvector corresponding to the
eigenvalue $\lambda_0$. According to
Remark~2.1~\cite{GurAsympAngle}, the vector
$\varphi^{(0)}(\omega)$ belongs to the space $W^{l+2m,N}(-b, b)$,
and, by Lemma~2.1~\cite{GurAsympAngle}, we have
\begin{equation}\label{eqNon-ClosedRange8H_a}
 \mathcal L_p V^0=0,
\end{equation}
where $V^0=r^{i\lambda_0}\varphi^{(0)}(\omega)$.

We substitute the sequence $U^\delta=r^\delta V^0/\|r^\delta
V^0\|_{H_a^{l+2m,N}(K^\varepsilon)}$, $\delta>0$,
into~\eqref{eqClosedRangeH_a} and let $\delta$ tend to zero.
Analogously to the proof of Lemma~\ref{lNon-ClosedRange}, by using
relation~\eqref{eqNon-ClosedRange8H_a}, one can check that the
right-hand side of inequality~\eqref{eqClosedRangeH_a} tends to
zero while its left-hand side is equal to one.
\end{proof}

Now the necessity of the conditions in Theorem~\ref{thLFredH_a}
follows from Lemma~\ref{lNon-ClosedRangeH_a}.

\section{Nonlocal Problems in Bounded Domains
in the Case where the Line $\Im\lambda=1-l-2m$ Contains an
Eigenvalue of $\tilde{\mathcal L}_p(\lambda)$}\label{sectHatLFred}

In the previous sections, we proved the Fredholm solvability and
obtained the asymptotics of solutions to problem~(\ref{eqPinG}),
(\ref{eqBinG}) in the case where the corresponding line in the
complex plane contains no eigenvalues of the operators
$\tilde{\mathcal L}_p(\lambda)$, $p=1, \dots, N_1$. In this
section, by using the results of Sec.~\ref{sectLinKProperEigen},
we thoroughly study the case where the line $\Im\lambda=1-l-2m$
contains only the proper eigenvalue $\lambda_0=i(1-l-2m)$ of the
operators $\tilde{\mathcal L}_p(\lambda)$ for some $p\in\{1,
\dots, N_1\}$. In this case, the operator $\mathbf
L:W^{l+2m}(G)\to {\mathcal W}^l(G,\Upsilon)$ is not Fredholm due
to Theorem~\ref{thLFred} (its image is not closed). Thus, we
associate to problem~(\ref{eqPinG}), (\ref{eqBinG}) an operator
acting in other spaces and prove that it is Fredholm.

\subsection{Construction of the right regularizer in the case where $\mathbf
B_{i\mu}^2=0$}\label{subsectHatLFredB^2=0}

We study nonlocal elliptic problem~(\ref{eqPinG}), (\ref{eqBinG})
under the following condition.
\begin{condition}\label{condProperEigenG}
The eigenvalue $\lambda_0=i(1-l-2m)$ is a proper eigenvalue of the
operators $\tilde{\mathcal L}_p(\lambda)$, $p\in\Pi$, where $\Pi$
is a nonempty subset of the set $\{1, \dots, N_1\}$. Neither of
the operators $\tilde{\mathcal L}_p(\lambda)$, $p=1, \dots, N_1$,
contains any other eigenvalues on the line $\Im\lambda=1-l-2m$.
\end{condition}

We introduce functions $\psi^p\in C_0^\infty(\mathbb R^2)$ such
that $\psi^p(y)=1$ for $y\in\bigcup\limits_{j=1}^{N_{1p}}\mathcal
O_{\varepsilon/2}(g_j^p)$ and
$\supp\psi^p\subset\bigcup\limits_{j=1}^{N_{1p}}\mathcal
O_\varepsilon(g_j^p)$. Here $\varepsilon>0$ is so small that
$\mathcal O_{\varepsilon}(g_j^p)\subset\mathcal V(g^{p}_j)$. We
also denote $\psi=1-\sum\limits_{p=1}^{N_{1}}\psi^p$. Let, to the
vector $\psi^p f= \{\psi^p f_0, \psi^p f_{i\mu}\}$ of the
right-hand sides in problem~(\ref{eqPinG}), (\ref{eqBinG}), there
correspond the vector $f^p=\{f_j^p, f_{j\sigma\mu}^p\}$ of the
right-hand sides in problem~(\ref{eqPinK}), (\ref{eqBinK}).
Clearly, $\supp f^p\subset\mathcal O_\varepsilon(0)$.

We introduce the space $\hat{\mathcal S}^l(G,\Upsilon)$ with the
norm
\begin{equation}\label{eqHatSNorm}
 \|f\|_{\hat{\mathcal S}^l(G,\Upsilon)}=\Bigl(\|\psi f\|^2_{\mathcal
 W^l(G,\Upsilon)}+\sum\limits_{p\in\Pi}
 \|f^p\|^2_{\hat{\mathcal S}^l(K^p,\gamma^p)}+\sum\limits_{p\notin\Pi}
 \|f^p\|^2_{{\mathcal S}^l(K^p,\gamma^p)}\Bigr)^{1/2}.
\end{equation}
According to Condition~\ref{condProperEigenG}, the set of indices
$\Pi$ is not empty; therefore, by Lemma~\ref{lS0NotClosed}, the
set $\hat{\mathcal S}^l(G,\Upsilon)$ is not closed in the topology
of $\mathcal W^l(G,\Upsilon)$.

On the other hand, it follows from Lemma~\ref{lConnectUS0lN} that,
provided $u\in W^{l+2m}(G)$ satisfies the relations
 \begin{equation}\label{eqConnectuS0l}
  D^\alpha u|_{y=g_j^p}=0,\quad |\alpha|\le l+2m-2;\ p=1,\dots,N_1;\
  j=1, \dots, N_{1p},
 \end{equation}
we have $\{\mathbf P u,\ \mathbf C u\}\in \hat{\mathcal
S}^l(G,\Upsilon)$ (the operator $\mathbf C=\mathbf B^0+\mathbf
B^1$ is defined in Sec.~\ref{sectLFred}). Let us introduce the
space
$$
S^{l+2m}(G)=\left\{u\in W^{l+2m}(G): u\ \text{satisfy
relations~\eqref{eqConnectuS0l}}\right\}
$$
and consider the operator
$$
 \hat{\mathbf L}^1=\{\mathbf P,\ \mathbf C\}:
 S^{l+2m}(G)\to  \hat{\mathcal S}^l(G,\Upsilon).
$$
Lemma~\ref{lConnectUS0lN} implies that the operator $\hat{\mathbf
L}^1$ is bounded.

\begin{lemma}\label{lRegL1inS}
Let Condition~\ref{condProperEigenG} hold. Then there exist a
bounded operator $\hat{\mathbf R}_1:\hat{\mathcal
S}^l(G,\Upsilon)\to S^{l+2m}(G)$ and a compact operator
$\hat{\mathbf T}_1: \hat{\mathcal S}^l(G,\Upsilon)\to
\hat{\mathcal S}^l(G,\Upsilon)$ such that
\begin{equation}\label{eqRegL1inS}
 \hat{\mathbf L}^1\hat{\mathbf R}_1=\hat{\mathbf I}+\hat{\mathbf
 T}_1,
\end{equation}
where $\hat{\mathbf I}$ denotes the unity operator in the space
$\hat{\mathcal S}^l(G,\Upsilon)$.
\end{lemma}
\textit{Proof} is analogous to that of Lemma~\ref{lRegL1} with the
following modifications: {\rm(I)} Theorem~\ref{thRegL1p} (which we
now apply to the orbits $\Orb_p$, $p\notin\Pi$) should be
supplemented with Theorem~\ref{thRegL1p0} (which we apply to the
orbits $\Orb_p$, $p\in\Pi$) and {\rm(II)}
Remark~\ref{remUinS^l+2m} should be taken into account.\qed

\subsection{Construction of the right regularizer in the case where $\mathbf
B_{i\mu}^2\ne0$} Theorem~\ref{thRegL1p'},
Remark~\ref{remUinS^l+2m}, and Theorem~\ref{thRegL1p0'} imply
that, for any sufficiently small $\varepsilon>0$, there exist
bounded operators
\begin{align*}
\hat{\mathbf R}_{\mathcal K}'&:\{f': \{0,f'\}\in \hat{\mathcal
S}^l(G,\Upsilon),\ \supp f'\subset\mathcal
O_{2\varepsilon}(\mathcal K)\}\to \{u\in S^{l+2m}(G): \supp
f'\subset\mathcal
O_{4\varepsilon}(\mathcal K)\},\\
\hat{\mathbf M}_{\mathcal K}',\hat{\mathbf T}_{\mathcal K}'&:\{f':
\{0,f'\}\in \hat{\mathcal S}^l(G,\Upsilon),\ \supp
f'\subset\mathcal O_{2\varepsilon}(\mathcal K)\}\to \hat{\mathcal
S}^l(G,\Upsilon)
\end{align*}
such that $\|\hat{\mathbf M}_{\mathcal K}'f'\|_{\hat{\mathcal
S}^{l}(G,\Upsilon)}\le c\varepsilon\|\{0,f'\}\|_{\hat{\mathcal
S}^{l}(G,\Upsilon)}$, where $c>0$ is independent of $\varepsilon$,
the operator $\hat{\mathbf T}_{\mathcal K}'$ is compact, and
$$
\hat{\mathbf L}^1\hat{\mathbf R}_{\mathcal K}'
f'=\{0,f'\}+\hat{\mathbf M}_{\mathcal K}'f'+\hat{\mathbf
T}_{\mathcal K}' f'.
$$

\smallskip

For any $f'$ such that $\{0,f'\}\in \hat{\mathcal
S}^l(G,\Upsilon)$, we put
\begin{equation}\notag
\hat{\mathbf R}'_1f'=\hat{\mathbf R}_{\mathcal K}'(\psi'
f')+\sum\limits_{j=1}^J \mathbf R'_{0j}(\psi'_j f'),
\end{equation}
where the functions $\psi',\psi_j'$ and the operators $\mathbf
R'_{0j}$ are the same as in Sec.~\ref{subsectRegB^2ne0}.

By using Theorems~\ref{thRegL1p'} and~\ref{thRegL1p0'}, one can
easily show that
\begin{equation}\label{eqL1R1'inS}
\hat{\mathbf L}^1\hat{\mathbf R}_1' f'=\{0,f'\}+\hat{\mathbf
M}_1'f'+\hat{\mathbf T}_1' f'.
\end{equation}
Here $\hat{\mathbf M}_1',\hat{\mathbf T}_1':\{f': \{0,f'\}\in
\hat{\mathcal S}^{l}(G,\Upsilon)\}\to \hat{\mathcal
S}^{l}(G,\Upsilon)$ are bounded operators such that
$\|\hat{\mathbf M}_1'f'\|_{\hat{\mathcal S}^{l}(G,\Upsilon)}\le
c\varepsilon\|\{0,f'\}\|_{\hat{\mathcal S}^{l}(G,\Upsilon)}$,
where $c>0$ is independent of $\varepsilon$, and the operator
$\hat{\mathbf T}_1'$ is compact.

With the help of the operators $\hat{\mathbf R}_1$ and
$\hat{\mathbf R}_1'$, we will construct the right regularizer for
problem~(\ref{eqPinG}), (\ref{eqBinG}) with $\mathbf
B_{i\mu}^2\ne0$. To this end, we will need the following
consistency condition to hold.
\begin{condition}\label{condS_1=S}
For any $u\in S^{l+2m}(G)$, we have $\{0, \mathbf
B^2u\}\in\hat{\mathcal S}^l(G,\Upsilon)$ and
$$
\|\{0, \mathbf B^2u\}\|_{\hat{\mathcal S}^l(G,\Upsilon)}\le
c\|u\|_{W^{l+2m}(G)}.
$$
\end{condition}

\begin{remark}
According to~\eqref{eqSeparK23'}, the operator $\mathbf B^2$
corresponds to nonlocal terms with the support outside the set
$\mathcal K$. Therefore, if Condition~\ref{condS_1=S} holds for
functions $u\in S^{l+2m}(G)$, it also holds for functions $u\in
W^{l+2m}(G\setminus\overline{\mathcal O_{\varkappa_1}(\mathcal
K)})$
\end{remark}

\begin{remark}\label{remS_1=S}
Let us illustrate with Example~\ref{exGeneralProblem} how to
achieve that Condition~\ref{condS_1=S} hold.

We consider problem~\eqref{eqPinGEx}, \eqref{eqBinGEx} and
additionally assume that the transformations $\Omega_{is}$ in this
problem satisfy condition~\eqref{eqOmega} (which is a restriction
on the geometrical structure of the transformations
$\Omega_{is}$). Then, by virtue of the continuity of
$\Omega_{is}$, we have $\Omega_{is}\big(\mathcal
O_\delta(g)\big)\subset\mathcal O_{\varepsilon_0/2}(\mathcal K)$
for any $g\in\bar\Upsilon_i\cap\mathcal K$, provided that
$\delta>0$ is small enough. Therefore, for any $u\in
W^{l+2m}(G\setminus\overline{\mathcal O_{\varkappa_1}(\mathcal
K)})$, we have
\begin{equation}\label{eqS_1=S}
 \mathbf B^2_{i\mu}u(y)=0\quad \text{for } y\in\mathcal O_\delta(\mathcal
 K)
\end{equation}
since $1-\zeta\big(\Omega_{is}(y)\big)=0$ for $y\in\mathcal
O_\delta(\mathcal K)$. In this case, Condition~\ref{condS_1=S}
obviously holds.

One may refuse condition~\eqref{eqOmega} but assume the following:
if $\Omega_{is}(g)\notin\mathcal K$ (where
$g\in\bar\Upsilon_i\cap\mathcal K$), then the coefficients of
$B_{i\mu s}(y,D_y)$ have zeros of certain orders at the points
$\Omega_{is}(g)$, which also guarantees that $\{0,\mathbf
B^2u\}\in\hat{\mathcal S}^l(G,\Upsilon)$ for any $u\in
W^{l+2m}(G\setminus\overline{\mathcal O_{\varkappa_1}(\mathcal
K)})$. However, in this paper, we do not study this issue in
detail.
\end{remark}

By virtue of Lemma~\ref{lConnectUS0lN} and
Condition~\ref{condS_1=S}, we have
$$
\{\mathbf Pu,\ \mathbf Bu\}\in\hat{\mathcal S}^l(G,\Upsilon)\quad
\text{for all } u\in S^{l+2m}(G).
$$
Therefore, the operator
$$
 \hat{\mathbf L}_S=\{\mathbf P,\ \mathbf B\}:
 S^{l+2m}(G)\to \hat{\mathcal S}^l(G,\Upsilon)
$$
is well defined and bounded by virtue of Lemma~\ref{lConnectUS0lN}
and condition~\ref{condS_1=S}.

\begin{lemma}\label{lRegLinS}
Let Conditions~\ref{condProperEigenG} and~\ref{condS_1=S} hold.
Then there exist a bounded operator $\hat{\mathbf R}:\hat{\mathcal
S}^l(G,\Upsilon)\to S^{l+2m}(G)$ and a compact operator
$\hat{\mathbf T}: \hat{\mathcal S}^l(G,\Upsilon)\to \hat{\mathcal
S}^l(G,\Upsilon)$ such that
\begin{equation}\label{eqRegLinS}
\hat{\mathbf L}_S\hat{\mathbf R}=\hat{\mathbf I}+\hat{\mathbf T}.
\end{equation}
\end{lemma}
\begin{proof}
We put $\Phi=\mathbf B^2\hat{\mathbf R}_1f$, where
$f=\{f_0,f'\}\in\hat{\mathcal S}^l(G,\Upsilon)$. Then, according
to Condition~\ref{condS_1=S}, the functions $\Phi$ and $\mathbf
B^2\hat{\mathbf R}_1'\Phi$ belong to the domain of the operator
$\hat{\mathbf R}'_1$. Therefore, we can define the bounded
operator $\hat{\mathbf R}_{\mathcal S}:\hat{\mathcal
S}^l(G,\Upsilon)\to S^{l+2m}(G)$ by the formula
$$
\hat{\mathbf R}_{\mathcal S}f=\hat{\mathbf R}_1f-\hat{\mathbf
R}_1'\Phi+\hat{\mathbf R}_1'\mathbf
 B^2\hat{\mathbf R}_1'\Phi.
$$
Analogously to the proof of Lemma~\ref{lRegL}, by using
equalities~\eqref{eqRegL1inS} and~\eqref{eqL1R1'inS}, one can show
that
$$
\hat{\mathbf L}_S\hat{\mathbf R}_{\mathcal S}=\hat{\mathbf I}+M+T,
$$
where $M,T:\hat{\mathcal S}^l(G,\Upsilon)\to\hat{\mathcal
S}^l(G,\Upsilon)$ are bounded operators such that $\|M\|\le
c\varepsilon$ ($c>0$ is independent of $\varepsilon$) and $T$ is
compact.

For $\varepsilon\le \frac{1}{2c}$, the operator $\hat{\mathbf
I}+M:\hat{\mathcal S}^l(G,\Upsilon)\to \hat{\mathcal
S}^l(G,\Upsilon)$ is invertible. Denoting $\hat{\mathbf
R}=\hat{\mathbf R}_{\mathcal S}(\hat{\mathbf I}+M)^{-1}$, $\mathbf
T=T(\hat{\mathbf I}+M)^{-1}$ yields~\eqref{eqRegLinS}.
\end{proof}

\subsection{Fredholm solvability of nonlocal problems}
Since the subspace $S^{l+2m}(G)$ is of finite codimension in
$W^{l+2m}(G)$, there exists a \textit{finite-dimensional} subspace
$\mathcal R^l(G,\Upsilon)$ in $\mathcal W^l(G,\Upsilon)$ such that
$$
\{\mathbf Pu,\ \mathbf Bu\}\in \hat{\mathcal
S}^l(G,\Upsilon)\oplus\mathcal R^l(G,\Upsilon)\quad \text{for all
} u\in W^{l+2m}(G).
$$
Therefore, we can define the bounded operator
$$
 \hat{\mathbf L}=\{\mathbf P,\ \mathbf B\}:W^{l+2m}(G)\to \hat{\mathcal
S}^l(G,\Upsilon)\oplus\mathcal R^l(G,\Upsilon).
$$

\begin{theorem}\label{thLFredS0}
Let Conditions~\ref{condProperEigenG} and~\ref{condS_1=S} hold.
Then the operator $\hat{\mathbf L}$ is Fredholm.
\end{theorem}
\begin{proof}
Lemmas~\ref{lLFiniteKernel} and~\ref{lRegLinS} and
Theorem~15.2~\cite{Kr} imply that the operator $\hat{\mathbf
L}_S:S^{l+2m}(G)\to \hat{\mathcal S}^l(G,\Upsilon)$ is Fredholm.
Since the domain $W^{l+2m}(G)$ of the operator $\hat{\mathbf L}$
is an extension of the domain $S^{l+2m}(G)$ of the operator
$\hat{\mathbf L}_S$ by a finite-dimensional subspace and
$\hat{\mathbf L}$ coincides with $\hat{\mathbf L}_S$ on
$S^{l+2m}(G)$, it follows that $\hat{\mathbf L}$ is also Fredholm.
\end{proof}

\section{Elliptic Problems with Homogeneous Nonlocal Conditions}\label{secL_B}
In this section, we study the operator corresponding to
problem~(\ref{eqPinG}), (\ref{eqBinG}) with the homogeneous
nonlocal conditions. By using the results of
Sec.~\ref{sectHatLFred}, we show that if the line
$\Im\lambda=1-l-2m$ contains only a proper eigenvalue, then the
operator under consideration, unlike the operator ${\mathbf L}$,
may be Fredholm. This turns out to depend on whether some
algebraic relations between the operators $\mathbf P$, $\mathbf
B^0$, and $\mathbf B^1$ hold at the points of the set $\mathcal
K$.

\subsection{The case where the line $\Im\lambda=1-l-2m$ contains no
eigenvalues of~$\tilde{\mathcal L}_p(\lambda)$}

Let us introduce the space
$$
 W_B^{l+2m}(G)=\{u\in W^{l+2m}(G): \mathbf Bu=0\}.
$$
Clearly, $W_B^{l+2m}(G)$ is a closed subspace in $W^{l+2m}(G)$. We
consider the bounded operator $\mathbf L_B :W_B^{l+2m}(G)\to
W^{l}(G)$ given by
$$
\mathbf L_Bu=\mathbf P u,\quad u\in W_B^{l+2m}(G).
$$

To study problem~(\ref{eqPinG}), (\ref{eqBinG}) with homogeneous
nonlocal conditions, we will need that the following conditions
for the operators $B_{i\mu s}(y, D_y)$ hold (see,
e.g.,~\cite[Ch.~2, \S~1]{LM}).
\begin{condition}\label{condNormBinG}
For all $i=1, \dots, N_0$, the system $\{B_{i\mu0}(y,
D_y)\}_{\mu=1}^m$ is normal on $\bar\Upsilon_i$ and the orders of
the operators $B_{i\mu s}(y, D_y)$ ($s=0,\dots,S_i$) are
$\le$$2m-1$.
\end{condition}

In this subsection, we prove the following result.
\begin{theorem}\label{thL_BFredNoEigen}
Let Condition~\ref{condNoEigenG} hold. Then the operator $\mathbf
L_B$ is Fredholm.

Let the line $\Im\lambda=1-l-2m$ contains an improper eigenvalue
$\lambda_0$ of the operator $\tilde{\mathcal L}_p(\lambda)$ for
some $p$ and Condition~\ref{condNormBinG} hold. Then the image of
the operator $\mathbf L_B$ is not closed (and, therefore, $\mathbf
L_B$ is not Fredholm).
\end{theorem}

Let, to the orbit $\Orb_p$, there correspond model
problem~(\ref{eqPinK0}), (\ref{eqBinK0}) in the angles $K_j=K_j^p$
with the sides $\gamma_{j\sigma}=\gamma_{j\sigma}^p$, $j=1, \dots,
N=N_{1p}$, $\sigma=1, 2$.

The following lemma allows one to reduce nonlocal problems with
nonhomogeneous nonlocal conditions to the corresponding problems
with homogeneous ones.
\begin{lemma}\label{lTraceUp}
Let Condition~\ref{condNormBinG} hold. Then, for any
$f_{j\sigma\mu}\in
H_a^{l+2m-m_{j\sigma\mu}-1/2}(\gamma_{j\sigma})$ with $\supp
f_{j\sigma\mu}\subset\mathcal O_{\varepsilon_1}(0)$
($\varepsilon_1>0$ is fixed), there exists a function $V\in
H_a^{l+2m,N}(K)$ such that $\supp V\subset\mathcal
O_{2\varepsilon_1}(0)$ and
\begin{equation}\label{eqTraceUp'}
 \mathbf B_{j\sigma\mu}(y,D_y)V|_{\gamma_{j\sigma}}=f_{j\sigma\mu},
\end{equation}
\begin{equation}\label{eqTraceUp''}
 \|V\|_{H_a^{l+2m,N}(K)}\le
 c_{\varepsilon_1}\sum\limits_{j,\sigma,\mu}
 \|f_{j\sigma\mu}\|_{H_a^{l+2m-m_{j\sigma\mu}-1/2}(\gamma_{j\sigma})},
\end{equation}
where $c_{\varepsilon_1}>0$ is independent of $f_{j\sigma\mu}$.
\end{lemma}
\begin{proof}
1. Analogously to the proof of Lemma~3.1~\cite{MP} (in which the
authors consider differential operators with constant
coefficients), one can construct functions $V_{j\sigma}\in
H_a^{l+2m}(K_j)$ such that
\begin{equation}\label{eqTraceUp7}
 B_{j\sigma\mu j0}(y, D_y)V_{j\sigma}|_{\gamma_{j\sigma}}=f_{j\sigma\mu},
\end{equation}
\begin{equation}\label{eqTraceUp8}
 \|V_{j\sigma}\|_{H_a^{l+2m}(K_j)}\le
 k_2\sum\limits_{\mu=1}^m
 \|f_{j\sigma\mu}\|_{H_a^{l+2m-m_{j\sigma\mu}-1/2}(\gamma_{j\sigma})}.
\end{equation}
Since $\supp f_{j\sigma\mu}\subset\mathcal O_{\varepsilon_1}(0)$,
one can assume that $\supp V_{j\sigma}\subset\mathcal
O_{2\varepsilon_1}(0)$

\smallskip

2. We denote $\delta=\min|(-1)^\sigma b_j+\omega_{j\sigma ks}\pm
b_k|/2$ ($j, k=1, \dots, N;\ \sigma=1, 2;\ s=1, \dots, S_{j\sigma
k}$) and introduce functions $\zeta_{j\sigma}\in
C_0^\infty(\mathbb R^2)$ such that $\zeta_{j\sigma}(\omega)=1$ for
$|(-1)^\sigma b_j-\omega|<\delta/2$ and
$\zeta_{j\sigma}(\omega)=0$ for $|(-1)^\sigma b_j-\omega|>\delta$.
Since the functions $\zeta_{j\sigma}$ are multipliers in the space
$H_a^{l+2m}(K_j)$, it follows from~(\ref{eqTraceUp7})
and~(\ref{eqTraceUp8}) that the function $V=(\zeta_{11}
V_{11}+\zeta_{12} V_{12}, \dots, \zeta_{N1} V_{N1}+\zeta_{N2}
V_{N2})$ satisfies conditions~(\ref{eqTraceUp'})
and~(\ref{eqTraceUp''}).
\end{proof}

\begin{remark}
One cannot repeat the analogous arguments in Sobolev spaces, since
the functions $\zeta_{j\sigma}$ are not multipliers in the spaces
$W^{l+2m}(K_j)$. Moreover, one can construct functions
$f_{j\sigma\mu}\in W^{l+2m-m_{j\sigma\mu}-1/2}(\gamma_{j\sigma})$
 ($j=1, \dots, N;\ \sigma=1, 2;\ \mu=1, \dots, m$) such that
neither of functions $V\in W^{l+2m,N}(K)$ satisfies
conditions~(\ref{eqTraceUp'}). This explains why the problem with
homogeneous nonlocal conditions is not equivalent to that with
nonhomogeneous conditions (i.e., the former may be Fredholm while
the latter is not Fredholm, see examples in Sec.~\ref{sectEx}).
\end{remark}

For each fixed orbit $\Orb_p$, we denote (as before)
$d_1=\min\{\chi_{j\sigma ks}, 1\}/2$, $d_2=2\max\{\chi_{j\sigma
ks}, 1\}$ and, for any $\varepsilon>0$, put
$d=d(\varepsilon)=2d_2\varepsilon$. The following result will be
used to study the image of the operator ${\mathbf L}_B$ (cf.
Lemma~\ref{lClosedRange}).

\begin{lemma}\label{lL_BClosedRange}
Let Condition~\ref{condNormBinG} hold and the image of ${\mathbf
L}_B$ be closed. Then, for each orbit $\Orb_p$, sufficiently small
$\varepsilon$ and all $U\in W^{l+2m,N}(K^{d})$ satisfying
relations~\eqref{eqConnectUS0lN} and such that
 \begin{equation}\label{eqL_BClosedRange'}
 \mathcal B_{j\sigma\mu}(D_y)U|_{\gamma^{2\varepsilon}_{j\sigma}}=0
 \quad (j=1, \dots, N;\ \sigma=1, 2;\ \mu=1, \dots, m),
\end{equation}
the following estimate holds:\footnote{Under the assumptions of
this lemma, it follows from Lemma~\ref{lWinHa} that $U_j\in
H_a^{l+2m}(K_j^{d})$ for any $a>0$. Therefore, $U_j\in
H_0^{l+2m-1}(K_j^{d})$ and estimate~(\ref{eqL_BClosedRange''}) is
well defined.}
\begin{equation}\label{eqL_BClosedRange''}
  \|U\|_{W^{l+2m,N}(K^\varepsilon)}\le
  c\sum\limits_{j=1}^N\big(\|\mathcal P_j(D_y) U_j\|_{
  W^{l}(K_j^{d})}+
  \|U_j\|_{H_0^{l+2m-1}(K_j^{d})}\big).
\end{equation}
\end{lemma}
\begin{proof}
1. Since the image of ${\mathbf L}_B$ is closed, it follows from
Lemma~\ref{lLFiniteKernel}, compactness of the embedding
$W^{l+2m}(G)\subset W^{l+2m-1}(G)$, and Theorem~7.1~\cite{Kr} that
\begin{equation}\label{eqL_BClosedRange1}
  \|u\|_{W^{l+2m}(G)}\le c(\|\mathbf P(y, D_y) u\|_{W^l(G)}+\|u\|_{W^{l+2m-1}(G)})
\end{equation}
for all $u\in W_B^{l+2m}(G)$. Let us substitute functions $u\in
W_B^{l+2m}(G)$ such that $\supp
u\in\bigcup\limits_{j=1}^{N_{1p}}\mathcal
O_{2\varepsilon_1}(g_j^p)$,
$2\varepsilon_1<\min\{\varepsilon_0,\varkappa_1\}$,
into~(\ref{eqL_BClosedRange1}). By virtue of~\eqref{eqSeparK23'},
for such functions, we have $\mathbf B^2u=0$. Therefore, using
Lemma~3.2~\cite[Ch.~2]{LM}, we see that, provided $\varepsilon_1$
is small enough, the estimate
\begin{equation}\label{eqL_BClosedRange2}
  \|U\|_{W^{l+2m,N}(K)}\le k_1\sum\limits_{j=1}^N\big(\|\mathcal P_j(D_y) U_j\|_{
  W^{l}(K_j)}+\|U_j\|_{W^{l+2m-1}(K_j)}\big),
\end{equation}
holds for all $U\in W^{l+2m,N}(K)$ such that $\supp
U\subset\mathcal O_{2\varepsilon_1}(0)$ and
\begin{equation}\label{eqL_BClosedRange3}
 \mathbf B_{j\sigma\mu}(y, D_y)U|_{\gamma_{j\sigma}}=0\quad
 (j=1, \dots, N;\ \sigma=1, 2;\ \mu=1, \dots, m).
\end{equation}

\smallskip

2. Let us show that, provided $\varepsilon_2<\varepsilon_1 d_1$ is
small enough, estimate~(\ref{eqL_BClosedRange2}) holds for all
$U\in W^{l+2m,N}(K)$ satisfying relations~\eqref{eqConnectUS0lN}
and such that $\supp U\subset\mathcal O_{2\varepsilon_2}(0)$ and
\begin{equation}\label{eqL_BClosedRange4}
 \mathcal B_{j\sigma\mu}(D_y)U|_{\gamma_{j\sigma}}=0\quad
 (j=1,\ \dots,\ N;\ \sigma=1,\ 2;\ \mu=1,\ \dots,\ m).
\end{equation}

We put $\Phi_{j\sigma\mu}=\mathbf B_{j\sigma\mu}(y,\
D_y)U|_{\gamma_{j\sigma}}$. Clearly,
\begin{equation}\label{eqSupportPhi2}
\supp\Phi\subset\mathcal O_{\varepsilon_2/d_1}(0)\subset\mathcal
O_{\varepsilon_1}(0).
\end{equation}

Let us fix some $a$, $0<a<1$, and prove that
\begin{equation}\label{eqL_BClosedRange6}
 \|\Phi_{j\sigma\mu}\|_{H_0^{l+2m-m_{j\sigma\mu}-1/2}(\gamma_{j\sigma})}\le
 k_2 \varepsilon_2^{1-a}\|U\|_{W^{l+2m,N}(K)}.
\end{equation}
By virtue of~(\ref{eqL_BClosedRange4}) and boundedness of the
trace operator in weighted spaces, it suffices to estimate the
terms of the following type:
$$
 \bigl(a_\alpha(y)-a_\alpha(0)\bigr)D^\alpha U_j\quad
 (|\alpha|=m_{j\sigma\mu}),\qquad
 a_\beta(y)D^\beta U_j\quad (|\beta|\le m_{j\sigma\mu}-1),
$$
where $a_\alpha$ and $a_\beta$ are infinitely differentiable
functions. Using the restriction on the support of $U_j$,
Lemma~$3.3'$~\cite{KondrTMMO67}, and Lemma~\ref{lWinHa}, we get
\begin{multline}\notag
\|\bigl(a_\alpha(y)-a_\alpha(0)\bigr)D^\alpha U_j\|_{
H_0^{l+2m-m_{j\sigma\mu}}(K_j)}\\
\le k_3
\varepsilon_2^{1-a}\|\bigl(a_\alpha(y)-a_\alpha(0)\bigr)D^\alpha
U_j\|_{
H_{a-1}^{l+2m-m_{j\sigma\mu}}(K_j)}\\
\le k_4\varepsilon_2^{1-a}\|D^\alpha
U_j\|_{H_{a}^{l+2m-m_{j\sigma\mu}}(K_j)}\le
k_5\varepsilon_2^{1-a}\|U_j\|_{W^{l+2m}(K_j)}.
\end{multline}
Similarly, by using Lemma~\ref{lWinHa}, we obtain
$$
 \|a_\beta(y)D^\beta U_j\|_{H_0^{l+2m-m_{j\sigma\mu}}(K_j)}\le
 k_6\varepsilon_2^{1-a}\|U_j\|_{H_{a-1}^{l+2m-1}(K_j)}\le k_7
 \varepsilon_2^{1-a}\|U_j\|_{W^{l+2m}(K_j)}.
$$
Thus, estimate~(\ref{eqL_BClosedRange6}) is proved.

Further, by virtue of~\eqref{eqSupportPhi2} and
Lemma~\ref{lTraceUp}, there exists a function $V=(V_1, \dots,
V_N)\in H_0^{l+2m,N}(K)$ such that $\supp V\subset\mathcal
O_{2\varepsilon_1}(0)$ and
\begin{gather}
 \mathbf B_{j\sigma\mu}(y,\ D_y)V|_{\gamma_{j\sigma}}=
 \Phi_{j\sigma\mu},\label{eqL_BClosedRange7}\\
 \|V\|_{H_0^{l+2m,N}(K)}\le
 c_{\varepsilon_1}\sum\limits_{j,\sigma,\mu}
 \|\Phi_{j\sigma\mu}\|_{
 H_0^{l+2m-m_{j\sigma\mu}-1/2}(\gamma_{j\sigma})},\label{eqL_BClosedRange8}
\end{gather}
where $c_{\varepsilon_1}$ is independent of $\varepsilon_2$.

Estimating $U-V$ with the help of~(\ref{eqL_BClosedRange2}) and
using inequalities~(\ref{eqL_BClosedRange8})
and~(\ref{eqL_BClosedRange6}), we get
\begin{multline}\notag
 \|U\|_{W^{l+2m,N}(K)}\le
 \|U-V\|_{W^{l+2m,N}(K)}+\|V\|_{W^{l+2m,N}(K)}\le \\
 \le k_8
 \sum\limits_{j=1}^N\big(\|\mathcal P_j(D_y) U_j\|_{
  W^{l}(K_j)}+\|U_j\|_{W^{l+2m-1}(K_j)}+\varepsilon_2^{1-a}\|U_j\|_{W^{l+2m}(K_j)}\big).
\end{multline}
Now, choosing sufficiently small $\varepsilon_2$, we obtain
estimate~(\ref{eqL_BClosedRange2}) valid for all $U\in
W^{l+2m,N}(K)$ with $\supp U\subset\mathcal O_{2\varepsilon_2}(0)$
and satisfying relations~(\ref{eqConnectUS0lN})
and~(\ref{eqL_BClosedRange4}).

\smallskip

3. Let us refuse the assumption $\supp U\subset\mathcal
O_{2\varepsilon_2}(0)$ and prove that, for
$\varepsilon<\varepsilon_2 d_1$ and any $U\in W^{l+2m,N}(K^{d})$
satisfying~\eqref{eqConnectUS0lN} and~\eqref{eqL_BClosedRange'},
estimate~\eqref{eqL_BClosedRange''} holds.

We introduce a function $\psi\in C_0^\infty(\mathbb R^2)$ such
that $\psi(y)=1$ for $|y|\le \varepsilon$,
$\supp\psi\subset\mathcal O_{2\varepsilon}(0)$, and $\psi$ does
not depend on polar angle $\omega$.

Put $\Psi_{j\sigma\mu}=\mathcal B_{j\sigma\mu}(D_y)(\psi
U)|_{\gamma_{j\sigma}}$. Clearly,
\begin{equation}\label{eqSupportPsi}
\supp\Psi_{j\sigma\mu}\subset\mathcal
O_{\varepsilon/d_1}(0)\subset\mathcal O_{\varepsilon_2}(0).
\end{equation}

Let us show that
\begin{equation}\label{eqL_BClosedRange9}
 \|\Psi_{j\sigma\mu}\|_{H_0^{l+2m-m_{j\sigma\mu}-1/2}(\gamma_{j\sigma})}\le
 k_9
\sum\limits_{k=1}^N\big(\|\mathcal P_k(D_y) U_k\|_{
W^{l}(K_k^{d})}+\|U_k\|_{H_0^{l+2m-1}(K_j^{d})}\big).
\end{equation}

Taking into account relations~\eqref{eqL_BClosedRange'}, we can
represent the function $\Psi_{j\sigma\mu}$ as follows:
\begin{equation}\label{eqPhi_jSigmaMu}
\Psi_{j\sigma\mu}=\sum\limits_{k,s}\Psi_{j\sigma\mu
ks}+\sum\limits_{(k,s)\ne(j,0)}J_{j\sigma\mu ks},
\end{equation}
where
\begin{align*}
 \Psi_{j\sigma\mu ks}&=\big([B_{j\sigma\mu ks}(D_y),\psi]U_k\big)\big(\mathcal G_{j\sigma
 ks}y\big)\big|_{\gamma_{j\sigma}},\\
 J_{j\sigma\mu ks}&=\big(\psi(\mathcal G_{j\sigma
ks}y)-\psi(y)\big)\big(B_{j\sigma\mu ks}(D_y)U_k\big)\big(\mathcal
G_{j\sigma ks}y\big)\big|_{\gamma_{j\sigma}}
\end{align*}
with $[\cdot,\cdot]$ denoting the commutator.

Since the expression for $\Psi_{j\sigma\mu ks}$ contains
derivatives of $U_k$ of order $\le$$m_{j\sigma\mu}-1$, it follows
that
\begin{equation}\label{eqL_BClosedRange10}
 \|\Psi_{j\sigma\mu ks}\|_{H_0^{l+2m-m_{j\sigma\mu}-1/2}(\gamma_{j\sigma})}\le
 k_{10}\|U_k\|_{H_0^{l+2m-1}(K_k^{d})}.
\end{equation}

Further, repeating the arguments of item~1 in the proof of
Lemma~\ref{lNon-ClosedRange}, we get
\begin{equation}\label{eqL_BClosedRange11}
\|J_{j\sigma\mu
ks}\|_{H_0^{l+2m-m_{j\sigma\mu}-1/2}(\gamma_{j\sigma})} \le
k_{11}(\|\mathcal
P_k(D_y)U_k\|_{W^l(\{d_1\varepsilon/2<|y|<2d_2\varepsilon\})}+
\|U_k\|_{W^{l+2m-1}(\{d_1\varepsilon/2<|y|<2d_2\varepsilon\})}).
\end{equation}
Now ~\eqref{eqL_BClosedRange9} follows
from~\eqref{eqPhi_jSigmaMu}, \eqref{eqL_BClosedRange10},
and~\eqref{eqL_BClosedRange11}.

\smallskip

4. By virtue of~\eqref{eqSupportPsi} and Lemma~\ref{lTraceUp}
(being applied to the operators $\mathcal B_{j\sigma\mu}(D_y)$),
there exists a function $V=(V_1, \dots, V_N)\in H_0^{l+2m,N}(K)$
such that $\supp V\subset\mathcal O_{2\varepsilon_2}(0)$ and
\begin{equation}\label{eqL_BClosedRange12}
 \mathcal B_{j\sigma\mu}(D_y)V|_{\gamma_{j\sigma}}=
 \Psi_{j\sigma\mu},
\end{equation}
\begin{equation}\label{eqL_BClosedRange13}
 \|V\|_{H_0^{l+2m,N}(K)}\le
 k_{12}\sum\limits_{j,\sigma,\mu}
 \|\Psi_{j\sigma\mu}\|_{
 H_0^{l+2m-m_{j\sigma\mu}-1/2}(\gamma_{j\sigma})}.
\end{equation}
Estimating $\psi U-V$ with the help of~(\ref{eqL_BClosedRange2})
and using Leibniz' formula and
inequalities~(\ref{eqL_BClosedRange13}),
(\ref{eqL_BClosedRange9}), we obtain
\begin{multline}\notag
 \|U\|_{W^{l+2m,N}(K^\varepsilon)}\le \|\psi U\|_{W^{l+2m,N}(K)}\le
\|\psi U-V\|_{W^{l+2m,N}(K)}+\|V\|_{W^{l+2m,N}(K)} \\
 \le k_{11}
 \sum\limits_{j=1}^N\big(\|\mathcal P_j(D_y) U_j\|_{
  W^{l}(K_j^{d})}+
  \|U_j\|_{H_0^{l+2m-1}(K_j^{d})}\big).
\end{multline}
\end{proof}

Lemma~\ref{lL_BClosedRange} allows us to prove that if the line
$\Im\lambda=1-l-2m$ contains an improper eigenvalue, then the
operator ${\mathbf L}_B$, like ${\mathbf L}$, is not Fredholm.

\begin{lemma}\label{lL_BNon-ClosedRange}
Let the line $\Im\lambda=1-l-2m$ contain an improper eigenvalue
$\lambda_0$ of the operator $\tilde{\mathcal L}_p(\lambda)$ for
some $p$ and Condition~\ref{condNormBinG} hold. Then the image of
${\mathbf L}_B$ is not closed.
\end{lemma}
\begin{proof}
1. Assume that the image of ${\mathbf L}_B$ is closed. We denote
by $\varphi^{(0)}(\omega), \dots, \varphi^{(\varkappa-1)}(\omega)$
an eigenvector and associate vectors corresponding to the
eigenvalue $\lambda_0$ (see~\cite{GS}). By virtue of
Remark~2.1~\cite{GurAsympAngle}, the vectors
$\varphi^{(k)}(\omega)$ belong to $W^{l+2m,N}(-b, b)$ and satisfy
\begin{equation}\label{eqL_BNon-ClosedRange1}
 \mathcal P_j(D_y) V_j^k=0,\quad  \mathcal B_{j\sigma\mu}(D_y) V^k=0.
\end{equation}
where $V^k=r^{i\lambda_0}\sum\limits_{s=0}^k\frac{1}{s!}(i\ln
r)^k\varphi^{(k-s)}(\omega)$, $k=0, \dots, \varkappa-1$. Since
$\lambda_0$ is not a proper eigenvalue, it follows that, for some
$k\ge0$, the function $V^k(y)$ is not a vector-polynomial. For
simplicity, we assume that
$V^0=r^{i\lambda_0}\varphi^{(0)}(\omega)$ is not a
vector-polynomial (the case where $k>0$ can be considered
analogously).

Let $\varepsilon$ and $d=d(\varepsilon)$ be the same constants as
in Lemma~\ref{lL_BClosedRange}. We consider the sequence
$U^\delta=r^\delta V^0/\|r^\delta
V^0\|_{W^{l+2m,N}(K^\varepsilon)}$. For any $\delta>0$, the
denominator is finite, but $\|r^\delta
V^0\|_{W^{l+2m,N}(K^\varepsilon)}\to\infty$ as $\delta\to 0$,
since $V^0$ is not a vector-polynomial. However,  $\|r^\delta
V^0\|_{H_0^{l+2m-1,N}(K^{d})}\le c$, where $c>0$ is independent of
$\delta\ge0$; therefore,
\begin{equation}\label{eqL_BNon-ClosedRange2}
 \|U^\delta\|_{H_0^{l+2m-1,N}(K^{d})}\to 0\quad \mbox{as } \delta\to 0.
\end{equation}
By using~(\ref{eqL_BNon-ClosedRange1}), analogously to the proof
of Lemma~\ref{lNon-ClosedRange}, one can check that
\begin{equation}\label{eqL_BNon-ClosedRange3}
 \|\mathcal P_j(D_y)U_j^\delta\|_{W^{l}(K_j^{d})}\to 0\quad \mbox{as } \delta\to 0.
\end{equation}\begin{equation}\label{eqL_BNon-ClosedRange4}
 \|\mathcal B_{j\sigma\mu}(D_y)U^\delta\big|_{\gamma_{j\sigma}^{3\varepsilon}}\|_{
 H_0^{l+2m-m_{j\sigma\mu}-1/2}(\gamma_{j\sigma}^{3\varepsilon})}\to 0\quad
 \mbox{as } \delta\to 0.
\end{equation}

\smallskip

2. We introduce the function $\psi\in C_0^\infty(\mathbb R^2)$
such that $\psi(y)=1$ for $y\in\mathcal O_{2\varepsilon}(0)$ and
$\supp\psi\subset \mathcal O_{3\varepsilon}(0)$.

Applying Lemma~\ref{lTraceUp} to the operators $\mathcal
B_{j\sigma\mu}(D_y)$ and functions $f_{j\sigma\mu}=\psi \mathcal
B_{j\sigma\mu}(D_y)U^\delta|_{\gamma_{j\sigma}}$ (notice that
$\supp f_{j\sigma\mu}\subset\mathcal O_{3\varepsilon}(0)$), we
obtain a function $W^\delta\in H_0^{l+2m,N}(K)$ ($\delta>0$) such
that $\supp W^\delta\subset\mathcal O_{6\varepsilon}(0)$ and
\begin{equation}\label{eqL_BNon-ClosedRange5}
 \mathcal B_{j\sigma\mu}(D_y)W^\delta|_{\gamma_{j\sigma}^{2\varepsilon}}=
 \mathcal B_{j\sigma\mu}(D_y)U^\delta|_{\gamma_{j\sigma}^{2\varepsilon}},
\end{equation}
\begin{equation}\label{eqL_BNon-ClosedRange6}
 \|W^\delta\|_{H_0^{l+2m,N}(K^{6\varepsilon})}\le
 k_1\sum\limits_{j,\sigma,\mu}
 \|\mathcal B_{j\sigma\mu}(D_y)U^\delta|_{\gamma_{j\sigma}^{3\varepsilon}}\|_{
 H_0^{l+2m-m_{j\sigma\mu}-1/2}(\gamma_{j\sigma}^{3\varepsilon})}.
\end{equation}
Moreover, the function $U^\delta-W^\delta$ satisfies
relations~(\ref{eqConnectUS0lN}); therefore we can apply
Lemma~\ref{lL_BClosedRange} to $U^\delta-W^\delta$. Then, from
estimate~(\ref{eqL_BClosedRange''}), using the boundedness of the
embedding $H_0^{l+2m}(K_j^{6\varepsilon})\subset
W^{l+2m}(K_j^{6\varepsilon})$ and
inequality~(\ref{eqL_BNon-ClosedRange6}), we get
\begin{multline}\label{eqL_BNon-ClosedRange7}
  \|U^\delta\|_{W^{l+2m,N}(K^\varepsilon)}\le
 \|U^\delta-W^\delta\|_{W^{l+2m,N}(K^\varepsilon)}+
 \|W^\delta\|_{W^{l+2m,N}(K^\varepsilon)}\\
  \le k_2\sum\limits_{j=1}^N\Big(\|\mathcal P_j(D_y) U^\delta_j\|_{
  W^{l}(K_j^{d})}+\sum\limits_{\sigma,\mu}\|\mathcal B_{j\sigma\mu}(D_y)U^\delta|_{\gamma_{j\sigma}^{3\varepsilon}}\|_{
 H_0^{l+2m-m_{j\sigma\mu}-1/2}(\gamma_{j\sigma}^{3\varepsilon})}\\
  +\|U^\delta_j\|_{H_0^{l+2m-1}(K_j^{d})}\Big).
\end{multline}
However,
assertions~(\ref{eqL_BNon-ClosedRange2})--(\ref{eqL_BNon-ClosedRange4})
contradict estimate~(\ref{eqL_BNon-ClosedRange7}), since
$\|U^\delta\|_{W^{l+2m,N}(K^\varepsilon)}=1$.
\end{proof}

\textit{Proof of Theorem~\ref{thL_BFredNoEigen}.} The first part
of Theorem~\ref{thL_BFredNoEigen} follows from
Theorem~\ref{thLFred}. The second part follows from
Lemma~\ref{lL_BNon-ClosedRange}. \qed

\subsection{The case where the line
$\Im\lambda=1-l-2m$ contains the proper eigenvalue of
$\tilde{\mathcal L}_p(\lambda)$}

It remains to study the case where the line $\Im\lambda=1-l-2m$
contains only the proper eigenvalue. Let
Condition~\ref{condProperEigenG} hold. Then we prove that the
Fredholm property of the operator $\mathbf L_B$, for a fixed $l\ge
1$, is determined by the following condition.
\begin{condition}\label{condAllDxiPj}
For $l\ge1$ and all $p\in\Pi$, system~\eqref{eqSystemBP'}
corresponding the orbit $\Orb_p$ contains all the operators
$D^{\xi}\mathcal P_{j}(D_y)$ ($|\xi|=l-1,$ $j=1, \dots,
N=N_{1p}$).
\end{condition}

\begin{theorem}\label{thL_BFred}
Let Condition~\ref{condProperEigenG} and Consistency
Condition~\ref{condS_1=S} hold. Then
\begin{enumerate}
\item[$1.$] the operator ${\mathbf L}_B: W_B^{2m}(G)\to L_2(G)$ is
Fredholm;
\item[$2.$] if $l\ge1$ and Condition~\ref{condAllDxiPj} holds, then
the operator ${\mathbf L}_B: W_B^{l+2m}(G)\to W^l(G)$ is Fredholm;
\item[$2'.$] if $l\ge1$, Condition~\ref{condAllDxiPj} does not hold, and
Condition~\ref{condNormBinG} holds, then the image of the operator
${\mathbf L}_B: W_B^{l+2m}(G)\to W^l(G)$ is not closed (and,
therefore, ${\mathbf L}_B$ is not Fredholm).
\end{enumerate}
\end{theorem}
\begin{proof}
1. Lemma~\ref{lLFiniteKernel} implies that the kernel of ${\mathbf
L}_B$ is finite-dimensional. Let us study the image $\mathcal
R(\mathbf L_B)$ of the operator $\mathbf L_B$.

\smallskip

2. First, we assume that $l\ge 1$ and Condition~\ref{condAllDxiPj}
holds. We claim that the set
\begin{equation}\label{eqL_BFred1}
 \bigl\{f_0\in W^l(G): \{f_0, 0\}\in\hat{\mathcal S}^l(G,\Upsilon)\bigr\}
\end{equation}
is a closed subset of finite codimension in $W^l(G)$. Indeed, let
$\psi^p$ be the functions appearing in the definition of the space
$\hat{\mathcal S}^l(G,\Upsilon)$ (see
Sec.~\ref{subsectHatLFredB^2=0}). Then, to the vector of
right-hand sides $\{\psi^pf_0,0\}$ in problem~(\ref{eqPinG}),
(\ref{eqBinG}), there corresponds some vector $\{f_j^p,0\}$ of the
right-hand sides in problem~(\ref{eqPinK}), (\ref{eqBinK}). Let
$p\in\Pi$. Clearly, $\mathcal T_{j\sigma\mu}\{f_j^p,0\}=0$.
Moreover, by virtue of Condition~\ref{condAllDxiPj},
relations~\eqref{eqRelPhif} are absent. Thus, due
to~\eqref{eqHatSNorm}, the norm of the function
$\{f_0,0\}\in\hat{\mathcal S}^l(G,\Upsilon)$ in $\hat{\mathcal
S}^l(G,\Upsilon)$ is equivalent to the norm of $f_0$ in $ W^l(G)$,
while the set~\eqref{eqL_BFred1} is the subspace in $W^l(G)$
consisting of functions which satisfy relations~\eqref{eqS1l1}.

Further, since $\hat{\mathcal S}^l(G,\Upsilon)\subset\hat{\mathcal
S}^l(G,\Upsilon)\oplus\mathcal R^l(G,\Upsilon)$, it follows that
the set
\begin{equation}\label{eqL_BFred2}
 \bigl\{f_0\in W^l(G): \{f_0, 0\}\in\hat{\mathcal
S}^l(G,\Upsilon)\oplus\mathcal R^l(G,\Upsilon)\bigr\}
\end{equation}
is also a close subset of finite codimension in $W^l(G)$. On the
other hand, $f_0\in \mathcal R(\mathbf L_B)$ if and only if
$\{f_0, 0\}\in\mathcal R(\hat{\mathbf L})$. Combining this with
the fact that the operator $\hat{\mathbf L}$ is Fredholm implies
that the image of $\mathbf L_B$ is closed and of finite
codimension.

\smallskip

3. Now we assume that $l\ge 1$ but Condition~\ref{condAllDxiPj}
fails. Let us prove that the image of $\mathbf L_B$ is not closed.
To this end, we will use the results of
Sec.~\ref{sectLinKProperEigen}. Since Condition~\ref{condAllDxiPj}
fails, the set of conditions~(\ref{eqRelPhif}) is not empty and,
for some $j, \xi$, the norm~(\ref{eqNormS0lN}) contains the
corresponding term $\|\mathcal T_{j\xi}f\|_{H_0^1(\mathbb R^2)}$.
Therefore, as follows from the proof of Lemma~\ref{lS0NotClosed},
there exists a sequence $f^\delta=\{f_j^\delta, 0\} \in
\hat{\mathcal S}^{l,N}(K,\gamma)$, $\delta>0$, such that $\supp
f^\delta\subset\mathcal O_{\varepsilon}(0)$ and $f^\delta$
converges in $\mathcal W^{l,N}(K,\gamma)$ to
$f^0\notin\hat{\mathcal S}^{l,N}(K,\gamma)$ as $\delta\to0$.

By virtue of Lemma~\ref{lInvL1p0}, for each $f^\delta$, there
exists a function $U^\delta\in W^{l+2m,N}(K^d)$ such that
\begin{equation}\label{eqL_BFred3}
 \mathcal P_j(D_y) U_j^\delta = f_j^\delta,\quad
  \mathcal B_{j\sigma\mu}(D_y) U^\delta = 0,
\end{equation}
\begin{equation}\label{eqL_BFred4}
\|U^\delta\|_{H_0^{l+2m-1,N}(K^{d})}\le c\|f^\delta\|_{\mathcal
W^{l,N}(K,\gamma)}
\end{equation}
($c>0$ is independent of $\delta$) and $U^\delta$ satisfies
relations~(\ref{eqConnectUS0lN}). By virtue of the second relation
in~(\ref{eqL_BFred3}) and relations~(\ref{eqConnectUS0lN}), we can
apply Lemma~\ref{lL_BClosedRange} to the function $U^\delta$. By
using estimate~(\ref{eqL_BClosedRange''}), convergence of
$f^\delta$ to $f^0\notin\hat{\mathcal S}^{l,N}(K,\gamma)$, and
inequality~(\ref{eqL_BFred4}), we arrive at the contradiction (cf.
the proof of Lemma~\ref{lNon-ClosedRange}).

\smallskip

4. In the case where $l=0$, the set of
conditions~(\ref{eqRelPhif}) is empty, since these conditions
appear only for $l\ge 1$. From this, similarly to item~2 of the
proof, we deduce the conclusion of the theorem.
\end{proof}

\section{Examples of Nonlocal Elliptic Problems in Sobolev Spaces}\label{sectEx}

In this section, we consider two examples which illustrate the
results of our work.

\subsection{Example~1}\label{subsectEx1}

\subsubsection{Problem with nonhomogeneous nonlocal conditions}

\begin{figure}[ht]
{ \hfill\epsfbox{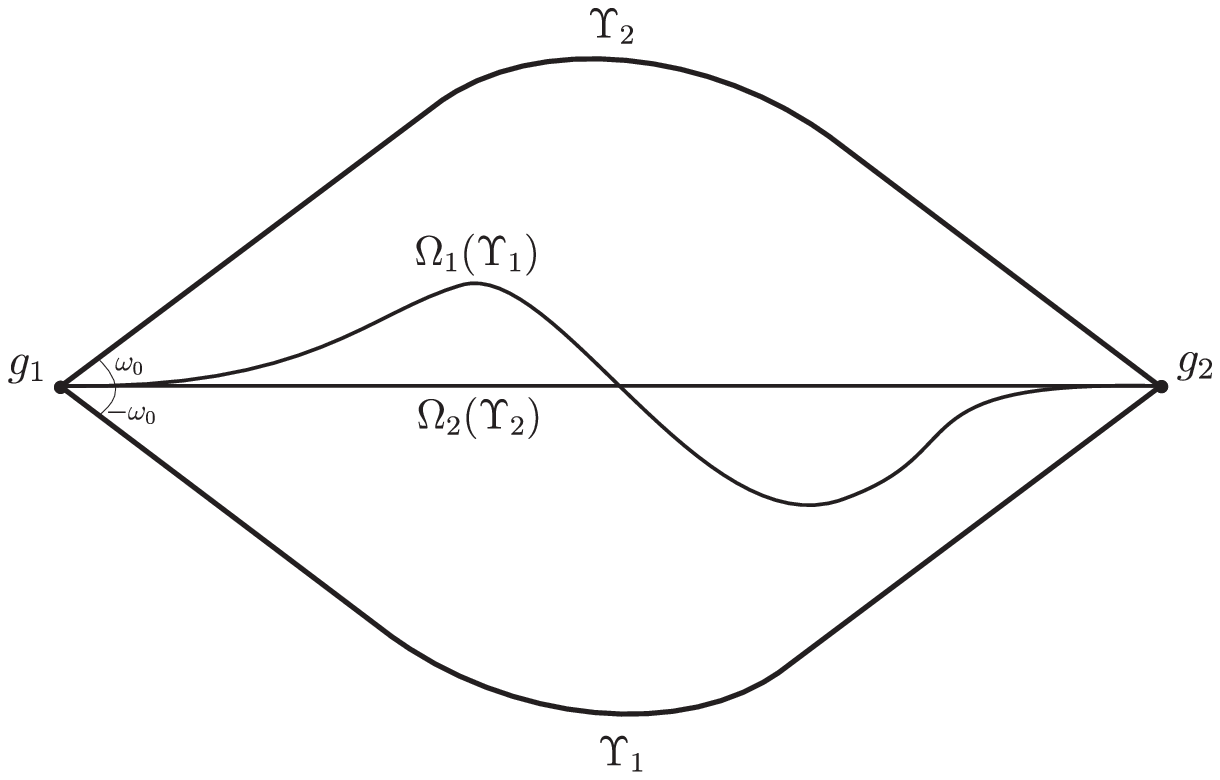}\hfill\ } \caption{Domain $G$ with the
boundary $\partial G=\bar\Upsilon_1\cup\bar\Upsilon_2$.}
   \label{figEx1}
\end{figure}

Let $\partial G\setminus{\mathcal
K}=\bigcup\limits_{i=1}^{2}\Upsilon_i$, where $\Upsilon_i$ are
open (in the topology of $\partial G$)  smooth curves and
$\mathcal K=\bar\Upsilon_1\cap\bar\Upsilon_2=\{g_1, g_2\}$ with
$g_1,g_2$ being the ends of the curves
$\bar\Upsilon_1,\bar\Upsilon_2$. We assume that, in neighborhoods
of the points $g_1, g_2$, the domain $G$ coincides with the plane
angles of the same opening $2\omega_0$, $0<2\omega_0<2\pi$. We
consider the following nonlocal problem in the domain $G$:
\begin{align}
 \Delta u&=f_0(y)\quad (y\in G),\label{eq1Statement1}\\
 u|_{\Upsilon_i}+b_i u\bigl(\Omega_i(y)\bigr)\big|_{\Upsilon_i}&=f_i(y)\quad
 (y\in\Upsilon_i;\ i=1, 2).\label{eq1Statement2}
\end{align}
Here $b_1, b_2\in\mathbb R$ and $\Omega_i$ is an infinitely
differentiable nondegenerate transformation mapping a neighborhood
${\mathcal O}_i$ of the curve $\Upsilon_i$ onto $\Omega({\mathcal
O}_i)$ so that $\Omega(\Upsilon_i)\subset G$, $\Omega_i(g_j)=g_j$,
$j=1, 2$, and, near the points $g_1, g_2$, the transformation
$\Omega_i$ is the rotation of $\Upsilon_i$ by the angle $\omega_0$
inwards $G$ (see Fig.~\ref{figEx1}).

According to Remark~\ref{remS_1=S}, Condition~\ref{condS_1=S}
holds. Clearly, Condition~\ref{condNormBinG} also holds.

\smallskip

To each of the points $g_1, g_2$, there corresponds the same model
problem in the plane angle:
\begin{align}
 \Delta U&=f_0(y)\quad (y\in K),\label{eq1Statement3}\\
 U|_{\gamma_j}+b_j U(\mathcal G_jy)|_{\gamma_j}&=f_j(y)\quad
 (y\in\gamma_j;\ j=1, 2).\label{eq1Statement4}
\end{align}
Here $K=\{y\in{\mathbb R}^2: r>0,\ |\omega|<\omega_0\}$,
$\gamma_j=\{y\in{\mathbb R}^2: r>0,\ \omega=(-1)^j\omega_0\}$, and
$$
 \mathcal G_j=
 \begin{pmatrix}
  \cos\omega_0 & (-1)^j\sin\omega_0\\
  (-1)^{j+1}\sin\omega_0 & \cos\omega_0
 \end{pmatrix}
$$
 is the operator of rotation by the angle $(-1)^{j+1}\omega_0$
around the origin, $j=1, 2$.

The eigenvalues problem corresponding to
problem~(\ref{eq1Statement3}), (\ref{eq1Statement4}) has the
following form:
\begin{gather}
\frac{d^2\varphi(\omega)}{d\omega^2}-\lambda^2\varphi(\omega)=0\quad
(|\omega|<\omega_0),\label{eq1Statement5}\\
\varphi(-\omega_0)+b_1\varphi(0)=0,\quad
\varphi(\omega_0)+b_2\varphi(0)=0.\label{eq1Statement6}
\end{gather}

Let us find the eigenvalues of problem~(\ref{eq1Statement5}),
(\ref{eq1Statement6}).

I. First, we consider the case where $\lambda\ne 0$. Substituting
the general solution
$\varphi(\omega)=c_1e^{\lambda\omega}+c_2e^{-\lambda\omega}$ for
Eq.~(\ref{eq1Statement5}) into nonlocal
condition~(\ref{eq1Statement6}), we get the following system of
equations:
\begin{equation}\label{eq1Statement7}
 \begin{pmatrix}
  e^{-\lambda\omega_0}+b_1 & e^{\lambda\omega_0}+b_1\\
  e^{\lambda\omega_0}+b_2 & e^{-\lambda\omega_0}+b_2
 \end{pmatrix}
 \begin{pmatrix}
  c_1\\
  c_2
 \end{pmatrix}=
 \begin{pmatrix}
  0\\
  0
 \end{pmatrix}.
\end{equation}
Equating the determinant of system~(\ref{eq1Statement7}) with
zero, we get
$$
(e^{-\lambda\omega_0}-e^{\lambda\omega_0})
(e^{\lambda\omega_0}+e^{-\lambda\omega_0}+b_1+b_2)=0.
$$

1. From the equation $e^{-\lambda\omega_0}-e^{\lambda\omega_0}=0$,
we obtain
\begin{equation}\label{eq1Statement8}
 \lambda=\frac{\pi k}{\omega_0}i,\quad k\in{\mathbb Z}\setminus\{0\}.
\end{equation}

2. Let us consider the equation
$e^{\lambda\omega_0}+e^{-\lambda\omega_0}+b_1+b_2=0$. If
$b_1+b_2=0$, then
\begin{equation}\label{eq1Statement9}
 \lambda=\frac{\pi/2+\pi k}{\omega_0}i,\quad k\in{\mathbb Z}.
\end{equation}
If $b_1+b_2\ne0$, then
\begin{equation}\label{eq1Statement10}
\lambda^\pm=\left\{
 \begin{aligned}
&\frac{\ln\left(-\frac{b_1+b_2}{2}\pm\frac{\sqrt{(b_1+b_2)^2-4}
}{2}\right)}{\omega_0} +\frac{2\pi n}{\omega_0}i\
     &\mbox{for }& b_1+b_2<-2,\\
&\frac{\pm\arctg\frac{\sqrt{4-(b_1+b_2)^2} }{b_1+b_2}+2\pi n}{\omega_0}i\
&\mbox{for }& -2<b_1+b_2<0,\\
&\frac{\pm\arctg\frac{\sqrt{4-(b_1+b_2)^2}
}{b_1+b_2}+(2n+1)\pi}{\omega_0}i\ &\mbox{for }& 0<b_1+b_2<2,\\
&\frac{\ln\left(\frac{b_1+b_2}{2}\pm\frac{\sqrt{(b_1+b_2)^2-4}
}{2}\right)}{\omega_0} +\frac{(2n+1)\pi}{\omega_0}i\
       &\mbox{for }& b_1+b_2>2,
  \end{aligned}
  \right.
\end{equation}
$n\in{\mathbb Z}$. For $|b_1+b_2|=2$, we get the eigenvalues from
the series~(\ref{eq1Statement8}).

II. Similarly, one can consider the case where $\lambda=0$ and
verify that $\lambda=0$ is an eigenvalue of
problem~(\ref{eq1Statement5}), (\ref{eq1Statement6}) if and only
if $b_1+b_2=-2$.

\medskip

Let us study the particular case where $\omega_0=\pi/2$, which
implies that $\partial G\in C^\infty$.

I. Let $\lambda\ne 0$.

1. From~(\ref{eq1Statement8}), we get the following pure imaginary
eigenvalues with integer imaginary parts:
\begin{equation}\label{eq1Statement11}
 \lambda_{2k}=2ki,\quad k\in{\mathbb Z}\setminus\{0\}.
\end{equation}

2. If $b_1+b_2=0$, we get from~(\ref{eq1Statement9}) the following
pure imaginary eigenvalues with integer imaginary parts:
\begin{equation}\label{eq1Statement12}
 \lambda_{2k+1}=(2k+1)i,\quad k\in{\mathbb Z}.
\end{equation}
If $b_1+b_2\ne 0$, we get from~(\ref{eq1Statement10}) the
following eigenvalues:
\begin{equation}\label{eq1Statement13}
\lambda_{n}^\pm=\left\{
 \begin{aligned}
&\frac{2\ln\left(-\frac{b_1+b_2}{2}\pm\frac{\sqrt{(b_1+b_2)^2-4}
}{2}\right)}{\pi} +4ni\ &\mbox{for }& b_1+b_2<-2,\\
&\frac{\pm 2\arctg\frac{\sqrt{4-(b_1+b_2)^2} }{b_1+b_2}}{\pi}i+4ni\
&\mbox{for }& -2<b_1+b_2<0,\\
&\frac{\pm 2\arctg\frac{\sqrt{4-(b_1+b_2)^2}
}{b_1+b_2}}{\pi}i+(4n+2)i\ &\mbox{for }& 0<b_1+b_2<2,\\
&\frac{2\ln\left(\frac{b_1+b_2}{2}\pm\frac{\sqrt{(b_1+b_2)^2-4}
}{2}\right)}{\pi}+(4n+2)i\
       &\mbox{for }& b_1+b_2>2,
  \end{aligned}
  \right.
\end{equation}
$n\in{\mathbb Z}$. If $|b_1+b_2|=2$, we get the eigenvalues from
the series~(\ref{eq1Statement11}).

II. The number $\lambda_0=0$ is an eigenvalue of
problem~(\ref{eq1Statement5}), (\ref{eq1Statement6}) if and only
if $b_1+b_2=-2$.

\smallskip

Let us consider the operator $\mathbf L: W^{l+2}(G)\to \mathcal
W^l(G,\Upsilon)$ corresponding to problem~(\ref{eq1Statement1}),
(\ref{eq1Statement2}) with $\omega_0=\pi/2$.
From~(\ref{eq1Statement11})--(\ref{eq1Statement13}) and
Theorem~\ref{thLFred}, we derive the following result.

\begin{theorem}\label{th1StatementLFred}
Suppose that $\omega_0=\pi/2$. Let $l$ be even; then the operator
$\mathbf L: W^{l+2}(G)\to \mathcal W^l(G,\Upsilon)$ is Fredholm if
and only if $b_1+b_2\ne 0$.

Let $l$ be odd; then the operator $\mathbf L: W^{l+2}(G)\to
\mathcal W^l(G,\Upsilon)$ is not Fredholm for any $b_1,
b_2\in\mathbb R$.
\end{theorem}

Notice that if $l$ is even and $b_1=b_2=0$, then the operator
$\mathbf L$ corresponding to the ``local'' boundary-value problem
is not Fredholm (its image is not closed). However, if we add
nonlocal terms with \textit{arbitrary small} coefficient $b_1,
b_2$ (such that $b_1+b_2\ne 0$) in the boundary-value conditions,
the problem becomes Fredholm.

\subsubsection{Problem with homogeneous nonlocal conditions}

Let us study problem~(\ref{eq1Statement1}), (\ref{eq1Statement2})
with homogeneous nonlocal conditions in the case where
$\omega_0=\pi/2$. We denote
$$
W_B^{l+2}(G)=\left\{u\in W^{l+2}(G):
u|_{\Upsilon_i}+b_i u\bigl(\Omega_i(y)\bigr)\big|_{\Upsilon_i}=0,\
i=1, 2\right\}
$$
and introduce the corresponding operator $\mathbf L_B
:W_B^{l+2}(G)\to W^{l}(G)$ given by
$$
 \mathbf L_B u= \Delta u, \quad u\in W_B^{l+2}(G).
$$

The Fredholm solvability of the operator $\mathbf L_B$ is
influenced only by the eigenvalues of
problem~(\ref{eq1Statement5}), (\ref{eq1Statement6}), lying on the
line $\Im\lambda=-(l+1)$, $l\ge0$. Thus, we have to consider only
the eigenvalues~(\ref{eq1Statement11}), (\ref{eq1Statement12}) for
$k\le -1$ and~(\ref{eq1Statement13}) for $|b_1+b_2|>2$, $n\le -1$.
Clearly, the eigenvalues~(\ref{eq1Statement13}) for $|b_1+b_2|>2$
are improper since they are not pure imaginary. Therefore, let us
begin with the question when the
eigenvalues~(\ref{eq1Statement11}) and~(\ref{eq1Statement12}) are
proper.

1. Consider the numbers $\lambda_{2k}=2ki$, $k=-1, -2, \dots$,
which are eigenvalues of problem~(\ref{eq1Statement5}),
(\ref{eq1Statement6}) for any $b_1, b_2$. Let us show that
$\lambda_{2k}$ is a proper eigenvalue if and only if $b_1+b_2\ne
2(-1)^{k+1}$.

To the eigenvalue $\lambda_{2k}$, there corresponds the
eigenvector
$\varphi_{2k}^{(0)}(\omega)=e^{i2k\omega}-e^{-i2k\omega}=2i\sin(2k\omega)$
(and, for $b_1=b_2=(-1)^{k+1}$, there is the second eigenvector
$\psi_{2k}^{(0)}(\omega)=e^{i2k\omega}+e^{-i2k\omega}=2\cos(2k\omega)$).
If an associate vector $\varphi_{2k}^{(1)}$ exists, then it
satisfies the equation
\begin{equation}\label{eq1sectHatLFred1}
 \frac{d^2\varphi_{2k}^{(1)}(\omega)}{d\omega^2}+4k^2\varphi_{2k}^{(1)}(\omega)=
 4ik\varphi_{2k}^{(0)}(\omega)\quad
(|\omega|<\pi/2)
\end{equation}
and nonlocal conditions~(\ref{eq1Statement6}). Substituting the
general solution
$$
 \varphi_{2k}^{(1)}(\omega)=
 c_1e^{i2k\omega}+c_2e^{-i2k\omega}+\omega(e^{i2k\omega}+e^{-i2k\omega})
$$
for Eq.~(\ref{eq1sectHatLFred1}) into nonlocal
conditions~(\ref{eq1Statement6}), we get the following system of
equations for $c_1, c_2$:
$$
 \begin{pmatrix}
  (-1)^k+b_1 & (-1)^k+b_1\\
  (-1)^k+b_2 & (-1)^k+b_2
 \end{pmatrix}
 \begin{pmatrix}
  c_1\\
  c_2
 \end{pmatrix}=
 \begin{pmatrix}
  \pi(-1)^k\\
  -\pi(-1)^k
 \end{pmatrix}.
$$
Clearly, this system is incompatible if and only if $b_1+b_2\ne
2(-1)^{k+1}$. Combining this with the fact that
$r^{-2k}\varphi_{2k}^{(0)}(\omega)$ is a polynomial with respect
to $y_1, y_2$ for $k=-1, -2, \dots$, we see that $\lambda_{2k}$ is
a proper eigenvalue if and only if $b_1+b_2\ne 2(-1)^{k+1}$.

\smallskip

2. Consider the numbers $\lambda_{2k+1}=(2k+1)i,\ k=-1, -2,
\dots$, which are eigenvalues of problem~(\ref{eq1Statement5}),
(\ref{eq1Statement6}) if and only if $b_1+b_2=0$. Let us prove
that, for $b_1+b_2=0$, the eigenvalues $\lambda_{2k+1}$ are proper
ones.

To the eigenvalue $\lambda_{2k+1}$, there corresponds the unique
eigenvector
$\varphi_{2k+1}^{(0)}(\omega)=e^{i(2k+1)\omega}+e^{-i(2k+1)\omega}=2i\sin((2k+1)\omega)$.
If an associate eigenvector $\varphi_{2k+1}^{(1)}$ exists, then it
satisfies the equation
\begin{equation}\label{eq1sectHatLFred2}
 \frac{d^2\varphi_{2k+1}^{(1)}(\omega)}{d\omega^2}+(2k+1)^2\varphi_{2k+1}^{(1)}(\omega)=
 2i(2k+1)\varphi_{2k+1}^{(0)}(\omega)\quad
(|\omega|<\pi/2)
\end{equation}
and nonlocal conditions~(\ref{eq1Statement6}). Substituting the
general solution
$$
\varphi_{2k+1}^{(1)}(\omega)=
c_1e^{i(2k+1)\omega}+c_2e^{-i(2k+1)\omega}+\omega(e^{i(2k+1)\omega}-e^{-i(2k+1)\omega})
$$
for Eq.~(\ref{eq1sectHatLFred2}) into nonlocal
conditions~(\ref{eq1Statement6}), we get the following system of
equations for $c_1, c_2$:
$$
 \begin{pmatrix}
  -i(-1)^k+b_1 & i(-1)^k+b_1\\
  i(-1)^k+b_2 & -i(-1)^k+b_2
 \end{pmatrix}
 \begin{pmatrix}
  c_1\\
  c_2
 \end{pmatrix}=
 \begin{pmatrix}
  -i\pi(-1)^k\\
  -i\pi(-1)^k
 \end{pmatrix}.
$$
Clearly, this system is incompatible for $b_1+b_2=0$. Combining
this with the fact that $r^{-(2k+1)}\varphi_{2k+1}^{(0)}(\omega)$
is a polynomial with respect to $y_1, y_2$ for $k=-1, -2, \dots$,
we see that, for $b_1+b_2=0$, the eigenvalues $\lambda_{2k+1}$ are
proper.

\begin{remark}\label{rAssocVec}
While checking whether an eigenvalue is proper, we sought only for
a first associate vector. Obviously, we can continue this
procedure and find the whole Jordan chain (see, e.g., Example~2.1
in~\cite{GurAsympAngle}); however, we do not do this here, since
the existence of a first associate vector already implies that the
corresponding eigenvalue is improper.
\end{remark}

\medskip

I. Consider the operator $\mathbf L_B :W_B^{2}(G)\to L_2(G)$. The
line $\Im\lambda=-1$ contains either no eigenvalues of
problem~(\ref{eq1Statement5}), (\ref{eq1Statement6}) (if
$b_1+b_2\ne0$) or only the proper eigenvalue $\lambda_{-1}=-i$ (if
$b_1+b_2=0$). Applying either Theorem~\ref{thL_BFredNoEigen} (if
$b_1+b_2\ne0$) or Theorem~\ref{thL_BFred} (if $b_1+b_2=0$), we see
that the \textit{operator $\mathbf L_B :W_B^{2}(G)\to L_2(G)$ is
Fredholm for all $b_1, b_2$.}

\medskip

II. Consider the operator $\mathbf L_B :W_B^{3}(G)\to W^1(G)$.

{\rm(a)} Let $b_1+b_2>2$. Then the line $\Im\lambda=-2$ contains
the proper eigenvalue $\lambda_{-2}=-2i$ and the two improper
eigenvalues
$\lambda_{-2}^\pm=\dfrac{2\ln\left(\frac{b_1+b_2}{2}\pm\frac{\sqrt{(b_1+b_2)^2-4}
}{2}\right)}{\pi}-2i$. Therefore, by
Theorem~\ref{thL_BFredNoEigen}, the operator $\mathbf L_B
:W_B^{3}(G)\to W^1(G)$ is not Fredholm.

{\rm(b)} Let $b_1+b_2=2$. Then the line $\Im\lambda=-2$ contains
only the improper eigenvalue $\lambda_{-2}=-2i$. Therefore, by
Theorem~\ref{thL_BFredNoEigen}, the operator $\mathbf L_B
:W_B^{3}(G)\to W^1(G)$ is not Fredholm.

{\rm(c)} Let $b_1+b_2<2$. Then the line $\Im\lambda=-2$ contains
only the proper eigenvalue $\lambda_{-2}=-2i$. We have to check
Condition~\ref{condAllDxiPj}. Differentiating the expression
$U(y)+b_jU(\mathcal G_jy)$ with respect to $y_2$ twice and
replacing the values of the corresponding function at the point
$\mathcal G_jy$ by the values at $y$, we see that
system~(\ref{eqSystemB}) has the following form:
$$
 \hat{\mathcal B}_1(D_y)U=\frac{\partial^2 U}{\partial y_2^2}+b_1\frac{\partial^2
 U}{\partial y_1^2},\quad
\hat{\mathcal B}_2(D_y)U=\frac{\partial^2 U}{\partial
y_2^2}+b_2\frac{\partial^2
 U}{\partial y_1^2}
$$

\begin{enumerate}
\item[$\rm (c_1)$] Let $b_1\ne b_2$. Then the operators $\hat{\mathcal
B}_1(D_y)U$ and $\hat{\mathcal B}_2(D_y)U$ are linearly
independent and, therefore, both included in
system~(\ref{eqSystemBP'}). Clearly, the operator $\Delta U$ is
not included in this system, since the system
$$
{\mathcal B}_1(D_y)U,\quad \hat{\mathcal B}_2(D_y)U,\quad \Delta
U
$$
is linearly dependent. Hence, Condition~\ref{condAllDxiPj} fails,
and Theorem~\ref{thL_BFred} implies that the operator $\mathbf L_B
:W_B^{3}(G)\to W^1(G)$ is not Fredholm.

\item[$\rm(c_2)$] Let $b_1=b_2$ (and, therefore, $b_1=b_2<1$). Then the operators
$\hat{\mathcal B}_1(D_y)U$ and $\hat{\mathcal B}_2(D_y)U$
coincide. Since $b_1<1$, the system
$$
\hat{\mathcal B}_1(D_y)U,\quad \Delta U
$$
is linearly independent and constitutes
system~(\ref{eqSystemBP'}). Hence, Condition~\ref{condAllDxiPj}
holds, and Theorem~\ref{thL_BFred} implies that the operator
$\mathbf L_B :W_B^{3}(G)\to W^1(G)$ is Fredholm.
\end{enumerate}

Thus, we proved that the \textit{operator $\mathbf L_B
:W_B^{3}(G)\to W^1(G)$ is Fredholm if and only if $b_1=b_2<1$.}

\medskip

III. Consider the operator $\mathbf L_B :W_B^{l+2}(G)\to W^l(G)$
with even $l\ge2$.

{\rm(a)} Let $b_1+b_2\ne 0$. Then the line $\Im\lambda=-(l+1)$
contains no eigenvalues of problem~(\ref{eq1Statement5}),
(\ref{eq1Statement6}). Therefore, by Theorem~\ref{thL_BFred}, the
operator $\mathbf L_B :W_B^{l+2}(G)\to W^l(G)$ is Fredholm.

{\rm(b)} Let $b_1+b_2=0$. Then the line $\Im\lambda=-(l+1)$
contains only the proper eigenvalue $\lambda_{-(l+1)}=-(l+1)i$.
Unlike the case where $l=0$, now we have to check
Condition~\ref{condAllDxiPj}. Differentiating the expression
$U(y)+b_jU(\mathcal G_jy)$  with respect to $y_2$~ $l+1$ times and
replacing the values of the corresponding function at the point
$\mathcal G_jy$ by the values at $y$, we see that
system~(\ref{eqSystemB}) has the form
$$
\hat{\mathcal B}_1(D_y)U=\frac{\partial^{l+1} U}{\partial
y_2^{l+1}}-b_1\frac{\partial^{l+1}
 U}{\partial y_1^{l+1}},\quad
\hat{\mathcal B}_2(D_y)U=\frac{\partial^{l+1} U}{\partial
y_2^{l+1}}+b_2\frac{\partial^{l+1}
 U}{\partial y_1^{l+1}}.
$$
Since $b_2=-b_1$, only the operator $\hat{\mathcal B}_1(D_y)U$ is
included in system~(\ref{eqSystemBP'}).

Let us show that the system consisting of the operator
$\hat{\mathcal B}_1(D_y)U$ and
$$
\frac{\partial^{l-1}}{\partial y_1^{\xi_1}\partial
y_2^{\xi_2}}\Delta U \equiv \frac{\partial^{l+1} U}{\partial
y_1^{\xi_1+2}\partial y_2^{\xi_2}}+\frac{\partial^{l+1}
U}{\partial y_1^{\xi_1}\partial y_2^{\xi_2+2}}\quad
(\xi_1+\xi_2=l-1)
$$
is linearly independent. To this end, we associate with each
derivative $\dfrac{\partial^{l+1} U}{\partial y_1^{s}\partial
y_2^{l+1-s}}$, $s=0, \dots, l+1$, the vector
$$
 (0, \dots, 0, 1, 0, \dots, 0)
$$
of length $l+2$ such that its $(s+1)$st component is equal to one
while all the remaining components are equal to zero. Then the
operator $\hat{\mathcal B}_1(D_y)U$ is associated with the vector
\begin{equation}\label{eq1sectHatLFred3}
  (1, 0, \dots, 0, -b_1)
\end{equation}
and the operators $\dfrac{\partial^{l-1}}{\partial
y_1^{\xi_1}\partial y_2^{\xi_2}}\Delta$, $\xi_1=0, \dots, l-1$,
are associated with the vectors
\begin{equation}\label{eq1sectHatLFred4}
  (0, \dots, 1, 0, 1, \dots, 0)
\end{equation}
such that their $(\xi_1+1)$st and $(\xi_1+3)$rd components are
equal to one while all the remaining components are equal to zero.
Thus, we have to show that the rank of the
$\big((l+1)\times(l+2)\big)$ order matrix
$$
 A=\begin{pmatrix}
  1 & 0 & 0 & 0 & \dots & 0 & 0 & -b_1\\
  1 & 0 & 1 & 0 & \dots & 0 & 0 & 0\\
  0 & 1 & 0 & 1 & \dots & 0 & 0 & 0\\
  \vdots & \vdots & \vdots & \vdots & \ddots & \vdots & \vdots & \vdots\\
  0 & 0 & 0 & 0 & \dots & 1 & 0 & 0\\
  0 & 0 & 0 & 0 & \dots & 0 & 1 & 0\\
  0 & 0 & 0 & 0 & \dots & 1 & 0 & 1\\
 \end{pmatrix}
$$
consisting of the rows~(\ref{eq1sectHatLFred3}),
(\ref{eq1sectHatLFred4}) is equal to $l+1$. We denote by $A'$ the
matrix obtained from the matrix $A$ by deleting the last column.
Decomposing the determinant of $A'$ by the first row, we see that
$\det A'=\det A_l$, where
$$
 A_l=\begin{pmatrix}
  0 & 1 & 0 & \dots & 0 & 0 & 0\\
  1 & 0 & 1 & \dots & 0 & 0 & 0\\
  0 & 1 & 0 & \dots & 0 & 0 & 0\\
  \vdots & \vdots & \vdots & \ddots & \vdots & \vdots & \vdots\\
  0 & 0 & 0 & \dots & 0 & 1 & 0\\
  0 & 0 & 0 & \dots & 1 & 0 & 1\\
  0 & 0 & 0 & \dots & 0 & 1 & 0
 \end{pmatrix}
$$
is a tridiagonal matrix of order $l\times l$. By induction, one
can easily check that
\begin{equation}\label{eq1sectHatLFred5}
  \det A_l=\left\{
   \begin{aligned}
    0\quad &\mbox{for}&l&=2n-1,\\
    1\quad &\mbox{for}&l&=4n,\\
    -1\quad &\mbox{for}&l&=4n-2,\\
   \end{aligned}
  \right.
\end{equation}
$n\ge 1$. From~(\ref{eq1sectHatLFred5}), it follows that $|\det
A'|=|\det A_l|=1$. Therefore, the system
$$
\hat{\mathcal B}_1(D_y)U,\quad \frac{\partial^{l-1}}{\partial
y_1^{\xi_1}\partial y_2^{\xi_2}}\Delta U\ (\xi_1+\xi_2=l-1)
$$
is linearly independent, and Theorem~\ref{thL_BFred} implies that
the operator $\mathbf L_B :W_B^{l+2}(G)\to W^l(G)$ is Fredholm.

Thus we proved that the \textit{operator $\mathbf L_B
:W_B^{l+2}(G)\to W^l(G)$ with even $l\ge2$ is Fredholm for any
$b_1, b_2$.}

\medskip

IV. Finally, we consider the operator $\mathbf L_B
:W_B^{l+2}(G)\to W^l(G)$ with odd $l\ge3$. First, we assume that
$l+1=4n$ for some $n\ge1$.

{\rm(a)} Let $b_1+b_2<-2$. Then the line $\Im\lambda=-(l+1)=-4n$
contains the proper eigenvalue $\lambda_{-4n}=-4ni$ and the two
improper eigenvalues
$\lambda_{-4n}^\pm=\dfrac{2\ln\left(\frac{b_1+b_2}{2}\pm\frac{\sqrt{(b_1+b_2)^2-4}
}{2}\right)}{\pi}-4ni$. Therefore, by
Theorem~\ref{thL_BFredNoEigen}, the operator $\mathbf L_B
:W_B^{l+2}(G)\to W^l(G)$ is not Fredholm.

{\rm(b)} Let $b_1+b_2=-2$. Then the line $\Im\lambda=-(l+1)=-4n$
contains only the improper eigenvalue $\lambda_{-4n}=-4ni$.
Therefore, by Theorem~\ref{thL_BFredNoEigen}, the operator
$\mathbf L_B :W_B^{l+2}(G)\to W^l(G)$ is not Fredholm.

{\rm(c)} Let $b_1+b_2>-2$. Then the line $\Im\lambda=-(l+1)=-4n$
contains only the proper eigenvalue $\lambda_{-2}=-4ni$. We have
to check Condition~\ref{condAllDxiPj}. Differentiating the
expression $U(y)+b_jU(\mathcal G_jy)$  with respect to $y_2$~
$l+1$ times and replacing the values of the corresponding function
at the point $\mathcal G_jy$ by the values at $y$, we see that
system~(\ref{eqSystemB}) has the form
$$
\hat{\mathcal B}_1(D_y)U=\frac{\partial^{l+1} U}{\partial
y_2^{l+1}}+b_1\frac{\partial^{l+1}
 U}{\partial y_1^{l+1}},\quad
\hat{\mathcal B}_2(D_y)U=\frac{\partial^{l+1} U}{\partial
y_2^{l+1}}+b_2\frac{\partial^{l+1}
 U}{\partial y_1^{l+1}}.
$$

\begin{enumerate}
\item[$\rm(c_1)$] Let $b_1\ne b_2$. Then the operators $\hat{\mathcal
B}_1(D_y)U$ and $\hat{\mathcal B}_2(D_y)U$ are linearly
independent and, therefore, both included in
system~(\ref{eqSystemBP'}). Let us show that the system
$$
\hat{\mathcal B}_1(D_y)U,\quad \hat{\mathcal B}_2(D_y)U,\quad
\frac{\partial^{l-1}}{\partial y_1^{\xi_1}\partial
y_2^{\xi_2}}\Delta U\ (\xi_1+\xi_2=l-1)
$$
is linearly dependent. (Notice that, unlike the case where $l=1$,
this system now contains all the $l+1$ order derivatives of $U$.)
Since $\hat{\mathcal B}_1(D_y)U$ and $\hat{\mathcal B}_2(D_y)U$
are linearly independent, it suffices to show that the system
$$
\frac{\partial^{l+1}U}{\partial y_2^{l+1}},\quad
\frac{\partial^{l+1}U}{\partial y_1^{l+1}},\quad
\frac{\partial^{l-1}}{\partial y_1^{\xi_1}\partial
y_2^{\xi_2}}\Delta U\ (\xi_1+\xi_2=l-1)
$$
is linearly dependent. We consider the corresponding matrix
$$
 A=\begin{pmatrix}
  1 & 0 & 0 & 0 & \dots & 0 & 0 & 0\\
  0 & 0 & 0 & 0 & \dots & 0 & 0 & 1\\
  1 & 0 & 1 & 0 & \dots & 0 & 0 & 0\\
  0 & 1 & 0 & 1 & \dots & 0 & 0 & 0\\
  \vdots & \vdots & \vdots & \vdots & \ddots & \vdots & \vdots & \vdots\\
  0 & 0 & 0 & 0 & \dots & 1 & 0 & 0\\
  0 & 0 & 0 & 0 & \dots & 0 & 1 & 0\\
  0 & 0 & 0 & 0 & \dots & 1 & 0 & 1\\
 \end{pmatrix}
$$
of order $(l+2)\times(l+2)$. Decomposing the determinant of $A$ by
the first row and then decomposing the determinant of the matrix
which we obtained by the first row again, we see that $\det A=\det
A_l$. Since $l$ is odd, it follows from~\eqref{eq1sectHatLFred5}
that $\det A=0$. Therefore, Condition~\ref{condAllDxiPj} fails,
and Theorem~\ref{thL_BFred}, implies that the operator $\mathbf
L_B :W_B^{l+2}(G)\to W^l(G)$ is not Fredholm.
\item[$\rm(c_2)$] Let $b_1=b_2$ (and, therefore, $b_1=b_2>-1$).
Then system~(\ref{eqSystemBP'}) contains only the operator
$\hat{\mathcal B}_1(D_y)U$. Let us show that the system
$$
\hat{\mathcal B}_1(D_y)U,\quad \frac{\partial^{l-1}}{\partial
y_1^{\xi_1}\partial y_2^{\xi_2}}\Delta U\ (\xi_1+\xi_2=l-1)
$$
is linearly independent. We consider the corresponding matrix
$$
 A=\begin{pmatrix}
  1 & 0 & 0 & 0 & \dots & 0 & 0 & b_1\\
  1 & 0 & 1 & 0 & \dots & 0 & 0 & 0\\
  0 & 1 & 0 & 1 & \dots & 0 & 0 & 0\\
  \vdots & \vdots & \vdots & \vdots & \ddots & \vdots & \vdots & \vdots\\
  0 & 0 & 0 & 0 & \dots & 1 & 0 & 0\\
  0 & 0 & 0 & 0 & \dots & 0 & 1 & 0\\
  0 & 0 & 0 & 0 & \dots & 1 & 0 & 1\\
 \end{pmatrix}
$$
of order $(l+1)\times(l+2)$. Deleting the second column from $A$,
decomposing the determinant of the matrix which we obtained by the
first row, and using relation~(\ref{eq1sectHatLFred5}), we get
\begin{equation}\notag
 \begin{vmatrix}
  1 &  0 & 0 & \dots & 0 & 0 & b_1\\
  1 &  1 & 0 & \dots & 0 & 0 & 0\\
  0 &  0 & 1 & \dots & 0 & 0 & 0\\
  \vdots & \vdots & \vdots & \ddots & \vdots & \vdots & \vdots\\
  0 &  0 & 0 & \dots & 1 & 0 & 0\\
  0 &  0 & 0 & \dots & 0 & 1 & 0\\
  0 &  0 & 0 & \dots & 1 & 0 & 1
 \end{vmatrix}
 =1-b_1\begin{vmatrix}
  1 &  1 & 0 & \dots & 0 & 0 & 0\\
  0 &  0 & 1 & \dots & 0 & 0 & 0\\
  0 &  1 & 0 & \dots & 0 & 0 & 0\\
  \vdots & \vdots & \vdots & \ddots & \vdots & \vdots & \vdots\\
  0 &  0 & 0 & \dots & 0 & 1 & 0\\
  0 &  0 & 0 & \dots & 1 & 0 & 1\\
  0 &  0 & 0 & \dots & 0 & 1 & 0
 \end{vmatrix}=1-b_1\det A_{l-1}=1+b_1\ne0
\end{equation}
since $b_1>-1$. Therefore, Condition~\ref{condAllDxiPj} holds, and
Theorem~\ref{thL_BFred} implies that the operator $\mathbf L_B
:W_B^{l+2}(G)\to W^l(G)$ is Fredholm.
\end{enumerate}

Thus, we proved that the \textit{operator $\mathbf L_B
:W_B^{l+2}(G)\to W^l(G)$ with $l+1=4n$, $n\ge1$, is Fredholm if
and only if $b_1=b_2>-1$.}

Analogously, by using~(\ref{eq1sectHatLFred5}) and
Theorem~\ref{thL_BFred}, one can show that the \textit{operator
$\mathbf L_B :W_B^{l+2}(G)\to W^l(G)$ with $l+1=4n+2$, $n\ge1$, is
Fredholm if and only if $b_1=b_2<1$.}

The following theorem summarize the results obtained.
\begin{theorem}\label{th1sectHatLFred}
Let $l$ be even; then the operator $\mathbf L_B: W_B^{l+2}(G)\to
W^l(G)$ is Fredholm for any $b_1, b_2\in\mathbb R$.

Let $l$ be odd and $l=4n+1,\ n=0, 1, 2, \dots$; then the operator
$\mathbf L_B: W_B^{l+2}(G)\to W^l(G)$ is Fredholm if and only if
$b_1=b_2<1$.

Let $l$ be odd and $l=4n+3,\ n=0, 1, 2, \dots$; then the operator
$\mathbf L_B: W_B^{l+2}(G)\to W^l(G)$ is Fredholm if and only if
$b_1=b_2>-1$.
\end{theorem}

Notice that, for $\omega_0=\pi/2$ and $b_1=b_2=0$, we have the
``local'' Dirichlet problem in a smooth domain with homogeneous
boundary-value conditions. In this case, as is well known, the
operator $\mathbf L_B: W_B^{l+2}(G)\to W^l(G)$ corresponding to
problem~(\ref{eq1Statement1}), (\ref{eq1Statement2}) with
homogeneous boundary-value conditions is not only Fredholm for any
$l\ge0$ but also invertible.

\subsection{Example~2}

\subsubsection{Problem with nonhomogeneous nonlocal conditions}

\begin{figure}[ht]
{ \hfill\epsfbox{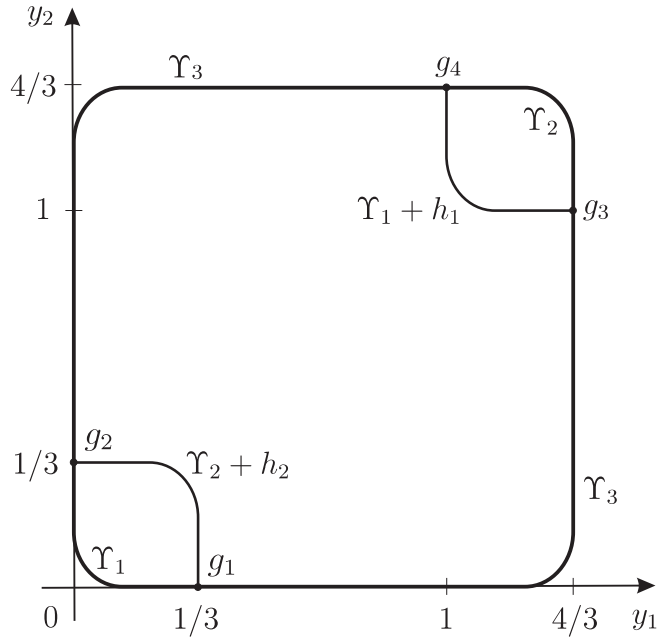}\hfill\ } \caption{Domain $G$ with smooth
boundary $\partial
G=\bar\Upsilon_1\cup\bar\Upsilon_2\cup\bar\Upsilon_3$.}
   \label{figEx2}
\end{figure}

Let $G\subset{\mathbb R}^2$ be a domain such that its boundary
$\partial G\in C^\infty$ coincides, outside the disks
$B_{1/8}((i4/3,j4/3))$ $(i,j=0,1)$, with the boundary of the
square $(0,4/3)\times(0,4/3)$. We denote
$\Upsilon_1=\{y\in\partial G: y_1<1/3,\ y_2<1/3\}$,
$\Upsilon_2=\{y\in\partial G: y_1>1,\ y_2>1\}$,
$\Upsilon_3=\partial G
\setminus(\bar\Upsilon_1\cup\bar\Upsilon_2)$. Thus, we have
$\mathcal K=\{g_1, \dots, g_4\}$, where $g_1=(1/3, 0),\ g_2=(0,
1/3),\ g_3=(4/3, 1),\ g_4=(1, 4/3)$ (see Fig.~\ref{figEx2}).

We consider the following nonlocal elliptic problem in the domain
$G$:
\begin{equation}\label{eq2Statement1}
\Delta u=f_0(y)\quad(y\in G),
\end{equation}
\begin{equation}\label{eq2Statement2}
\begin{aligned}
u(y)|_{\Upsilon_i}+b_i u(y+h_i)|_{\Upsilon_i}&=f_i(y)\quad(y\in\Upsilon_i;\ i=1, 2),\\
 u(y)|_{\Upsilon_3}&=f_3(y)\quad(y\in\Upsilon_3),
\end{aligned}
\end{equation}
where $h_1=(1, 1)$, $h_2=(-1, -1)$, and $b_1, b_2 \in{\mathbb R}$.
Clearly, $\mathcal K=\Orb_1\cup\Orb_2$, where the orbit $\Orb_1$
consists of the points $g_1$ and $g_3=g_1+h_1$ and the orbit
$\Orb_2$ consists of the points $g_2$ and $g_4=g_2+h_2$.

According to Remark~\ref{remS_1=S}, Condition~\ref{condS_1=S}
holds. Clearly, Condition~\ref{condNormBinG} also holds.

\smallskip

First, we assume that $b_1^2+b_2^2\ne0$ (for definiteness, we
suppose that $b_1\ne0$).

To each of the orbits $\Orb_1,\ \Orb_2$, there corresponds the
same model problem in the plane angles:
\begin{equation}\label{eq2Statement3}
 \Delta U_j=f_j(y)\quad (y\in K),
\end{equation}
\begin{equation}\label{eq2Statement4}
 \begin{matrix}
 U_1|_{\gamma_1}=f_{11}(y)\ (y\in\gamma_1),\quad
 U_1|_{\gamma_2}+b_1 U_2(\mathcal G y)|_{\gamma_2}=f_{12}(y)\
 (y\in\gamma_2),\\
 U_2|_{\gamma_1}=f_{21}(y)\ (y\in\gamma_1),\quad
 U_2|_{\gamma_2}+b_2 U_1(\mathcal G y)|_{\gamma_2}=f_{22}(y)\
 (y\in\gamma_2).
 \end{matrix}
\end{equation}
Here $K=\{y\in{\mathbb R}^2: r>0,\ |\omega|<\pi/2\}$,
$\gamma_j=\{y\in{\mathbb R}^2: r>0,\ \omega=(-1)^j\pi/2\}$, and $
 \mathcal G=
 \begin{pmatrix}
  0 & 1\\
  -1 & 0
 \end{pmatrix}
$ is the operator of rotation by the angle $-\pi/2$.

The eigenvalues problem corresponding to
problem~(\ref{eq2Statement3}), (\ref{eq2Statement4}) has the
following form:
\begin{equation}\label{eq2Statement5}
\frac{d^2\varphi_j(\omega)}{d\omega^2}-\lambda^2\varphi_j(\omega)=0\quad
(|\omega|<\pi/2;\ j=1, 2),
\end{equation}
\begin{equation}\label{eq2Statement6}
\begin{matrix}
\varphi_1(-\pi/2)=0,\quad \varphi_1(\pi/2)+b_1\varphi_2(0)=0,\\
\varphi_2(-\pi/2)=0,\quad \varphi_2(\pi/2)+b_2\varphi_1(0)=0.
\end{matrix}
\end{equation}

One can find the eigenvalues of problem~(\ref{eq2Statement5}),
(\ref{eq2Statement6}) by straightforward computations
(see~\cite{SkRJMP}). They are
\begin{align}
\lambda_{2k}&=2ki,\quad k\in{\mathbb Z}\setminus\{0\}\quad
(\text{for all }b_1,b_2,\
b_1^2+b_2^2\ne0),\label{eq2Statement7}\\
 \lambda_{2k+1}&=(2k+1)i,\quad k\in{\mathbb Z}\quad (\mbox{for }
 b_2=0,\ b_1\ne0),\label{eq2Statement8}
\end{align}
and
\begin{equation}\label{eq2Statement9}
 \lambda_{n}^\pm=\left\{
 \begin{aligned}
&\frac{2}{\pi}\ln\left|\frac{\sqrt{-b_1b_2}}{2}\pm\frac{\sqrt{4-b_1b_2}
}{2}\right| +(2n+1)i\ &\mbox{for }& b_1b_2<0,\\
&\left(\pm \frac{2}{\pi}\arctg\sqrt{4(b_1b_2)^{-1}-1}+2n\right)i\
&\mbox{for }& 0<b_1b_2<4,\\
&\frac{2}{\pi}\ln\left(\frac{\sqrt{b_1b_2}}{2}\pm\frac{\sqrt{b_1b_2-4}}{2}\right)+2ni
&\mbox{for }& b_1b_2\ge4,
  \end{aligned}
  \right.
\end{equation}
$n\in{\mathbb Z}$. If $b_1b_2=4$, then there is one more
eigenvalue $\lambda_0=0$.

\begin{remark}\label{rb2=0}
If $b_2=0$, we can consider the other statement of nonlocal
problem different from~(\ref{eq2Statement1}),
(\ref{eq2Statement2}), namely:
\begin{equation}\label{eq2Statement1'}
\Delta u=f(y)\quad(y\in G),
\end{equation}
\begin{equation}\label{eq2Statement2'}
\begin{aligned}
u(y)|_{\Upsilon_1}+b_1 u(y+h_1)|_{\Upsilon_1}&=f_1(y)\quad(y\in\Upsilon_1),\\
u(y)|_{\bar\Upsilon_2\cup\Upsilon_3}&=f_2(y)\quad(y\in\bar\Upsilon_2\cup\Upsilon_3).
\end{aligned}
\end{equation}
In this case, we have $\mathcal K=\{g_1, g_2\}$ (notice that
Condition~\ref{condS_1=S} now fails). Solutions to
problem~(\ref{eq2Statement1'}), (\ref{eq2Statement2'}) may have
singularities only near the points $g_1, g_2$, while solutions to
problem~(\ref{eq2Statement1}), (\ref{eq2Statement2}) may have
those near $g_1, \dots, g_4$.

To each of the points $g_1, g_2$, there corresponds the same model
``\textit{local}'' problem:
\begin{gather}
 \Delta U_1=f_1(y)\quad (y\in K),\label{eq2Statement3'}\\
 U_1|_{\gamma_1}=f_{1}(y)\ (y\in\gamma_1),\quad
 U_1|_{\gamma_2}=f_{2}(y)\ (y\in\gamma_2).\label{eq2Statement4'}
\end{gather}

The eigenvalues problem for problem~(\ref{eq2Statement3'}),
(\ref{eq2Statement4'}) has the following form:
\begin{gather}
\frac{d^2\varphi_1(\omega)}{d\omega^2}-\lambda^2\varphi_1(\omega)=0\quad
(|\omega|<\pi/2),\label{eq2Statement5'}\\
\varphi_1(-\pi/2)=\varphi_1(\pi/2)=0.\label{eq2Statement6'}
\end{gather}

The eigenvalues of problem~(\ref{eq2Statement5'}),
(\ref{eq2Statement6'}) are as follows:
\begin{equation}\label{eq2Statement7'}
 \lambda_{k}=ki,\quad k\in{\mathbb Z}\setminus\{0\}.
\end{equation}
They coincide with the eigenvalues of
problem~(\ref{eq2Statement5}), (\ref{eq2Statement6}) for $b_2=0$.
Therefore, according to Theorem~\ref{thLFred},
problem~(\ref{eq2Statement1'}), (\ref{eq2Statement2'}) is Fredholm
if and only if problem~(\ref{eq2Statement1}),
(\ref{eq2Statement2}) is Fredholm.
\end{remark}

Let us consider the operator $\mathbf L: W^{l+2}(G)\to \mathcal
W^l(G,\Upsilon)$ corresponding to problem~(\ref{eq2Statement1}),
(\ref{eq2Statement2}).
From~(\ref{eq2Statement7})--(\ref{eq2Statement9}) and
Theorem~\ref{thLFred}, we derive the following result.

\begin{theorem}\label{th2StatementLFred}
Let $l$ be even; then the operator $\mathbf L: W^{l+2}(G)\to
\mathcal W^l(G,\Upsilon)$ is Fredholm if and only if $b_1b_2>0$.

Let $l$ be odd; then the operator $\mathbf L: W^{l+2}(G)\to
\mathcal W^l(G,\Upsilon)$ is not Fredholm for any $b_1,
b_2\in\mathbb R$.
\end{theorem}

Notice that Theorem~\ref{th2StatementLFred} is proved under the
assumption that $b_1^2+b_2^2\ne0$; but the operator $\mathbf L:
W^{l+2}(G)\to \mathcal W^l(G,\Upsilon)$ with $b_1=b_2=0$
corresponding to problem~(\ref{eq2Statement1}),
(\ref{eq2Statement2}) is not Fredholm either. This follows from
the fact that, to each of the points $g_1, \dots, g_4\in\mathcal
K$, there corresponds model problem~(\ref{eq2Statement5'}),
(\ref{eq2Statement6'}) with the eigenvalues~\eqref{eq2Statement7'}
lying on the lines $-(l+1)$, $l\ge0$.

\subsubsection{Problem with homogeneous nonlocal conditions}

Let us study problem~(\ref{eq2Statement1}), (\ref{eq2Statement2})
with homogeneous nonlocal conditions. We denote
$$
W_B^{l+2}(G)=\left\{u\in W^{l+2}(G): u|_{\Upsilon_i}+b_i
u(y+h_i)|_{\Upsilon_i}=0,\ i=1, 2;\ u|_{\Upsilon_3}=0\right\}
$$
and introduce the corresponding operator $\mathbf L_B
:W_B^{l+2}(G)\to W^{l}(G)$ given by
$$
 \mathbf L_B u= \Delta u, \quad u\in W_B^{l+2}(G).
$$
First, we assume that $b_1^2+b_2^2\ne0$ (for definiteness, we
again suppose that $b_1\ne0$).
\begin{remark}
Problem~(\ref{eq2Statement1'}), (\ref{eq2Statement2'}) with
homogeneous nonlocal conditions is equivalent to
problem~(\ref{eq2Statement1}), (\ref{eq2Statement2}) with $b_2=0$.
Therefore, there is no need to study
problem~(\ref{eq2Statement1'}), (\ref{eq2Statement2'})
independently.
\end{remark}

The Fredholm property of the operator $\mathbf L_B$ is influenced
only by the eigenvalues of problem~(\ref{eq2Statement5}),
(\ref{eq2Statement6}), lying on the line $\Im\lambda=-(l+1)$,
$l\ge0$. Thus, we have to consider only the eigenvalues
$\lambda_{2k}$, $\lambda_{2k+1}$ (if $b_2=0$), and $\lambda_n^\pm$
(if $b_1b_2<0$ or $b_1b_2\ge4$) for $k, n\le -1$. Clearly, the
eigenvalues $\lambda_n^\pm$ (if $b_1b_2<0$ of $b_1b_2\ge4$) are
improper, since they are not pure imaginary. Therefore, let us
begin  with the question when the eigenvalues $\lambda_{2k}$ and
$\lambda_{2k+1}$ (if $b_2=0$) are proper.

1.  Consider the numbers $\lambda_{2k}=2ki$, $k=-1, -2, \dots$,
which are eigenvalues of problem~(\ref{eq2Statement5}),
(\ref{eq2Statement6}) for any $b_1, b_2$. Let us show that
$\lambda_{2k}$ is a proper eigenvalue.

To the eigenvalue $\lambda_{2k}$, there correspond the two
linearly independent eigenvectors
\begin{align*}
\bigl(\varphi_{1,2k}^{(0)}(\omega),\
\varphi_{2,2k}^{(0)}(\omega)\bigr)
&=\bigl(e^{i2k\omega}-e^{-i2k\omega},\ 0\bigr)=\bigl(2i\sin(2k\omega),\ 0\bigr),\\
\bigl(\psi_{1,2k}^{(0)}(\omega),\
\psi_{2,2k}^{(0)}(\omega)\bigr)&=\bigl(0,\
e^{i2k\omega}-e^{-i2k\omega}\bigr)=\bigl(0,\
2i\sin(2k\omega)\bigr).
\end{align*}

If an associate vector $(\varphi_{1,2k}^{(1)},\
\varphi_{2,2k}^{(1)})$ corresponding to the first of the
eigenvectors exists, then it satisfies the equations
\begin{equation}\label{eq2sectHatLFred1}
 \begin{aligned}
 \frac{d^2\varphi_{1,2k}^{(1)}(\omega)}{d\omega^2}+4k^2\varphi_{1,2k}^{(1)}(\omega)&=
 4ik(e^{i2k\omega}-e^{-i2k\omega})\quad
(|\omega|<\pi/2),\\
 \frac{d^2\varphi_{2,2k}^{(1)}(\omega)}{d\omega^2}+4k^2\varphi_{2,2k}^{(1)}(\omega)&=
 0\quad
(|\omega|<\pi/2)
\end{aligned}
\end{equation}
and nonlocal conditions~(\ref{eq2Statement6}). Substituting the
general solution
\begin{align*}
 \varphi_{1,2k}^{(1)}(\omega)&=
 c_1e^{i2k\omega}+c_2e^{-i2k\omega}+\omega(e^{i2k\omega}+e^{-i2k\omega}),\\
 \varphi_{2,2k}^{(1)}(\omega)&=
 c_3e^{i2k\omega}+c_4e^{-i2k\omega}
\end{align*}
for Eqs.~(\ref{eq2sectHatLFred1}) into nonlocal
conditions~(\ref{eq2Statement6}), we get the following system of
equations for $c_1, \dots, c_4$:
$$
 \begin{pmatrix}
  (-1)^k & (-1)^k & 0      & 0\\
  (-1)^k & (-1)^k & b_1    & b_1\\
  0      & 0      & (-1)^k & (-1)^k\\
  b_2    & b_2    & (-1)^k & (-1)^k
 \end{pmatrix}
 \begin{pmatrix}
  c_1\\
  c_2\\
  c_3\\
  c_4
 \end{pmatrix}=
 \begin{pmatrix}
  \pi(-1)^k\\
  -\pi(-1)^k\\
  0\\
  0
 \end{pmatrix}.
$$
It is easy to see that this system is incompatible; therefore, the
first eigenvector has no associate ones. Similarly, one can check
that the second eigenvector also has no associate ones. Combining
this with the fact that $r^{-2k}\varphi_{j,2k}^{(0)}(\omega)$ and
$r^{-2k}\psi_{j,2k}^{(0)}(\omega)$ ($j=1, 2$) are polynomials with
respect to $y_1, y_2$ for $k=-1, -2, \dots$, we see that
$\lambda_{2k}$ is a proper eigenvalue.

\smallskip

2. Consider the numbers $\lambda_{2k+1}=(2k+1)i$, $k=-1, -2,
\dots$, which are eigenvalues of problem~(\ref{eq2Statement5}),
(\ref{eq2Statement6}) with $b_2=0$ (we remind that $b_1\ne0$). Let
us show that $\lambda_{2k+1}$ is an improper eigenvalue.

To the eigenvalue $\lambda_{2k+1}$, there corresponds the only
eigenvector
$$
\bigl(\varphi_{1,2k+1}^{(0)}(\omega),\
\varphi_{2,2k+1}^{(0)}(\omega)\bigr)
=(e^{i(2k+1)\omega}+e^{-i(2k+1)\omega},\ 0)=(2\cos((2k+1)\omega),\
0).
$$

If an associate eigenvector $(\varphi_{1,2k+1}^{(1)},\
\varphi_{2,2k+1}^{(1)})$ exists, then it satisfies the equations
\begin{equation}\label{eq2sectHatLFred2}
 \begin{aligned}
\frac{d^2\varphi_{1,2k+1}^{(1)}(\omega)}{d\omega^2}+(2k+1)^2\varphi_{1,2k+1}^{(1)}(\omega)&=
 2(2k+1)i(e^{i(2k+1)\omega}+e^{-i(2k+1)\omega})\quad
(|\omega|<\pi/2),\\
 \frac{d^2\varphi_{2,2k+1}^{(1)}(\omega)}{d\omega^2}+(2k+1)^2\varphi_{2,2k+1}^{(1)}(\omega)&=
 0\quad
(|\omega|<\pi/2)
\end{aligned}
\end{equation}
and nonlocal conditions~(\ref{eq2Statement6}). Substituting the
general solution
\begin{align*}
 \varphi_{1,2k}^{(1)}(\omega)&=
 c_1e^{i(2k+1)\omega}+c_2e^{-i(2k+1)\omega}+\omega(e^{i(2k+1)\omega}-e^{-i(2k+1)\omega}),\\
 \varphi_{2,2k}^{(1)}(\omega)&=
 c_3e^{i(2k+1)\omega}+c_4e^{-i(2k+1)\omega}
\end{align*}
for Eqs.~(\ref{eq2sectHatLFred2}) into nonlocal
conditions~(\ref{eq2Statement6}), we get the following system of
equations for $c_1, \dots, c_4$:
$$
 \begin{pmatrix}
  i(-1)^{k+1} & i(-1)^k & 0           & 0\\
  i(-1)^k & i(-1)^{k+1} & b_1         & b_1\\
  0       & 0           & i(-1)^{k+1} & i(-1)^k\\
  0       & 0           & i(-1)^k     & i(-1)^{k+1}
 \end{pmatrix}
 \begin{pmatrix}
  c_1\\
  c_2\\
  c_3\\
  c_4
 \end{pmatrix}=
 \begin{pmatrix}
  \pi i(-1)^{k+1}\\
  \pi i(-1)^{k+1}\\
  0\\
  0
 \end{pmatrix}.
$$
It is easy to see that this system is compatible. Therefore,
$\lambda_{2k+1}$ is an improper eigenvalue.

\medskip

I. Consider the operator $\mathbf L_B :W_B^{2}(G)\to L_2(G)$. The
line $\Im\lambda=-1$ contains either no eigenvalues of
problem~(\ref{eq2Statement5}), (\ref{eq2Statement6}) (if
$b_1b_2>0$) or the improper eigenvalue $\lambda_{-1}$ (if $b_2=0$)
or $\lambda_{-1}^\pm$ (if $b_1b_2<0$). Therefore, by
Theorem~\ref{thL_BFredNoEigen}, the \textit{operator $\mathbf L_B
:W_B^{2}(G)\to L_2(G)$ is Fredholm if and only if $b_1b_2>0$.}

\medskip

II. Consider the operator $\mathbf L_B :W_B^{3}(G)\to W^1(G)$.

{\rm(a)}. Let $b_1b_2\ge4$. Then the line $\Im\lambda=-2$ contains
the proper eigenvalue $\lambda_{-2}$ and the two improper
eigenvalues $\lambda_{-1}^\pm$. Therefore, by
Theorem~\ref{thL_BFredNoEigen}, the operator $\mathbf L_B
:W_B^{3}(G)\to W^1(G)$ is not Fredholm.

{\rm(b)} Let $b_1b_2<4$. Then the line $\Im\lambda=-2$ contains
only the proper eigenvalue $\lambda_{-2}=-2i$. Let us show that
Condition~\ref{condAllDxiPj} fails.

Differentiating the expressions $U_1(y)+b_1U_2(\mathcal Gy)$ and
$U_2(y)+b_2U_1(\mathcal Gy)$ with respect to $y_2$ twice and
replacing the values of the corresponding functions at the point
$\mathcal Gy$ by the values at $y$, we see that
system~(\ref{eqSystemB}) has the following form:
 \begin{align*}
 \hat{\mathcal B}_{11}(D_y)U&=\frac{\partial^2 U_1}{\partial y_2^2},\quad
 \hat{\mathcal B}_{12}(D_y)U=\frac{\partial^2 U_1}{\partial y_2^2}+b_1\frac{\partial^2
 U_2}{\partial y_1^2},\\
 \hat{\mathcal B}_{21}(D_y)U&=\frac{\partial^2 U_2}{\partial y_2^2},\quad
 \hat{\mathcal B}_{22}(D_y)U=\frac{\partial^2 U_2}{\partial
y_2^2}+b_2\frac{\partial^2 U_1}{\partial y_1^2}.
 \end{align*}
Since $b_1\ne0$, the operators $\hat{\mathcal B}_{11}(D_y)U$,
$\hat{\mathcal B}_{12}(D_y)U$, and $\hat{\mathcal B}_{21}(D_y)U$
are linearly independent and, therefore, included in
system~(\ref{eqSystemBP'}). But the system consisting of these
three operators and the operators $\Delta U_1$ and $\Delta U_2$ is
linearly dependent. Therefore, Condition~\ref{condAllDxiPj} fails,
and Theorem~\ref{thL_BFred} implies that the operator $\mathbf L_B
:W_B^{3}(G)\to W^1(G)$ is not Fredholm.

Thus, we proved that the \textit{operator $\mathbf L_B
:W_B^{3}(G)\to W^1(G)$ is not Fredholm for any $b_1,b_2$
($b_1^2+b_2^2\ne0$).}

\medskip

III. Consider the operator $\mathbf L_B :W_B^{l+2}(G)\to W^l(G)$
with even $l\ge2$. The line $\Im\lambda=-(l+1)$ contains either no
eigenvalues of problem~(\ref{eq2Statement5}),
(\ref{eq2Statement6}) (if $b_1b_2>0$) or the improper eigenvalue
$\lambda_{-(l+1)}$ (if $b_2=0$) of $\lambda_{-1-l/2}^\pm$ (if
$b_1b_2<0$). Therefore, by Theorem~\ref{thL_BFredNoEigen}, the
\textit{operator $\mathbf L_B :W_B^{l+2}(G)\to W^l(G)$ with even
$l\ge2$ is Fredholm if and only if $b_1b_2>0$.}

\medskip

IV. Consider the operator $\mathbf L_B :W_B^{l+2}(G)\to W^l(G)$
with odd $l\ge3$.

{\rm(a)} Let $b_1b_2\ge4$. Then the line $\Im\lambda=-(l+1)$
contains the proper eigenvalue $\lambda_{-(l+1)}$ and the two
improper eigenvalues $\lambda_{-1/2-l/2}^\pm$. Therefore, by
Theorem~\ref{thL_BFredNoEigen}, the operator $\mathbf L_B
:W_B^{l+2}(G)\to W^l(G)$ is not Fredholm.

{\rm(b)} Let $b_1b_2<4$. Then the line $\Im\lambda=-(l+1)$
contains only the proper eigenvalue $\lambda_{-(l+1)}=-(l+1)i$.
Let us show that Condition~\ref{condAllDxiPj} fails.
Differentiating the expressions $U_1(y)+b_1U_2(\mathcal Gy)$ and
$U_2(y)+b_2U_1(\mathcal Gy)$ with respect to $y_2$~  $l+1$ times
and replacing the values of the corresponding functions at the
point $\mathcal Gy$ by the values at $y$, we see that
system~(\ref{eqSystemB}) has the following form:
 \begin{align*}
 \hat{\mathcal B}_{11}(D_y)U&=\frac{\partial^{l+1} U_1}{\partial
 y_2^{l+1}},\quad
\hat{\mathcal B}_{12}(D_y)U=\frac{\partial^{l+1} U_1}{\partial
y_2^{l+1}}+b_1\frac{\partial^{l+1} U_2}{\partial y_1^{l+1}},\\
 \hat{\mathcal B}_{21}(D_y)U&=\frac{\partial^{l+1} U_2}{\partial
 y_2^{l+1}},\quad
 \hat{\mathcal B}_{22}(D_y)U=\frac{\partial^{l+1} U_2}{\partial
y_2^{l+1}}+b_2\frac{\partial^{l+1} U_1}{\partial y_1^{l+1}}.
 \end{align*}
Since $b_1\ne0$, the operators $\hat{\mathcal B}_{11}(D_y)U$,
$\hat{\mathcal B}_{12}(D_y)U$, and $\hat{\mathcal B}_{21}(D_y)U$
are linearly independent and, therefore, included in
system~(\ref{eqSystemBP'}). Let us show that the system of these
three operators and
\begin{align*}
\frac{\partial^{l-1}}{\partial y_1^{\xi_1}\partial
y_2^{\xi_2}}\Delta U_1 &\equiv \frac{\partial^{l+1} U_1}{\partial
y_1^{\xi_1+2}\partial y_2^{\xi_2}}+\frac{\partial^{l+1}
U_1}{\partial y_1^{\xi_1}\partial
y_2^{\xi_2+2}}\quad  (\xi_1+\xi_2=l-1),\\
 \frac{\partial^{l-1}}{\partial y_1^{\xi_1}\partial
y_2^{\xi_2}}\Delta U_2 &\equiv \frac{\partial^{l+1} U_2}{\partial
y_1^{\xi_1+2}\partial y_2^{\xi_2}}+\frac{\partial^{l+1}
U_2}{\partial y_1^{\xi_1}\partial y_2^{\xi_2+2}}\quad
(\xi_1+\xi_2=l-1)
\end{align*}
is linearly dependent. To this end, we associate with each
derivative $\dfrac{\partial^{l+1} U_1}{\partial y_1^{s}\partial
y_2^{l+1-s}}$, $s=0, \dots, l+1$, the vector
$$
 (0, \dots, 0, 1, 0, \dots, 0)
$$
of length $2l+4$ such that its $(s+1)$st component is equal to one
while all the remaining components are equal to zero. Further, we
associate with each derivative $\dfrac{\partial^{l+1}
U_2}{\partial y_1^{s}\partial y_2^{l+1-s}}$, $s=0, \dots, l+1$,
the vector
$$
 (0, \dots, 0, 1, 0, \dots, 0)
$$
of length $2l+4$ such that its $(l+2+s+1)$st component is equal to
one while all the remaining components are equal to zero. Thus, it
suffices to show that the rank of the $(2l+3)\times(2l+4)$ order
matrix
$$
 A=\begin{pmatrix}
  1 & 0 & 0 & \dots & 0\rule{4mm}{0mm}\vrule  & 0 & 0 & 0 & \dots & 0\\
  1 & 0 & 0 & \dots & 0\rule{4mm}{0mm}\vrule  & 0 & 0 & 0 & \dots & b_1\\
  0 & 0 & 0 & \dots & 0\rule{4mm}{0mm}\vrule  & 1 & 0 & 0 & \dots & 0\\
  \hline
  1 & 0 & 1 & \dots & 0\rule{4mm}{0mm}\vrule  & 0 & 0 & 0 & \dots & 0\\
  0 & 1 & 0 & \dots & 0\rule{4mm}{0mm}\vrule  & 0 & 0 & 0 & \dots & 0\\
  \vdots & \vdots & \vdots & \ddots & \vdots\rule{5mm}{0mm}\vrule & \vdots & \vdots & \vdots & \ddots &
  \vdots\\
  0 & 0 & 0 & \dots & 1\rule{4mm}{0mm}\vrule & 0 & 0 & 0 & \dots & 0\\
  \hline
  0 & 0 & 0 & \dots & 0\rule{4mm}{0mm}\vrule & 1 & 0 & 1 & \dots & 0\\
  0 & 0 & 0 & \dots & 0\rule{4mm}{0mm}\vrule & 0 & 1 & 0 & \dots & 0\\
  \vdots & \vdots & \vdots & \ddots & \vdots\rule{5mm}{0mm}\vrule & \vdots & \vdots & \vdots & \ddots &
  \vdots\\
  0 & 0 & 0 & \dots & 0\rule{4mm}{0mm}\vrule & 0 & 0 & 0 & \dots & 1
 \end{pmatrix}
$$
is less than $2l+3$. (In the matrix $A$, the first three rows
correspond to the operators $\hat{\mathcal B}_{11}(D_y)U$,
$\hat{\mathcal B}_{12}(D_y)U$, and $\hat{\mathcal B}_{21}(D_y)U$
respectively, the next $l+2$ rows correspond to the operators
$\dfrac{\partial^{l-1}}{\partial y_1^{\xi_1}\partial
y_2^{\xi_2}}\Delta U_1$, and the last $l+2$ rows correspond to the
operators $\dfrac{\partial^{l-1}}{\partial y_1^{\xi_1}\partial
y_2^{\xi_2}}\Delta U_2$.)

If we delete the $1$st, $(l+3)$rd, or $(2l+4)$th columns from $A$,
then the $1$st row, the $3$rd row, or the difference between the
$1$st and $2$nd rows in the square matrix which we obtained is
equal to zero. Let us denote by $\hat A$ the matrix which is
obtained from $A$ by deleting any other column. Then, decomposing
the determinant of $\hat A$ consecutively by the first three rows,
we see that $|\det \hat A|=|b_1\det A'|$, where $A'$ is the
$2l\times 2l$ order matrix that is obtained from the
$2l\times(2l+1)$ order matrix
$$
 A''=\begin{pmatrix}
  0 & 1 & \dots & 0\rule{4mm}{0mm}\vrule  & 0 & 0 & \dots & 0 & 0 \\
  1 & 0 & \dots & 0\rule{4mm}{0mm}\vrule  & 0 & 0 & \dots & 0 & 0\\
  \vdots & \vdots & \ddots & \vdots\rule{5mm}{0mm}\vrule & \vdots & \vdots & \ddots
& \vdots & \vdots\\
  0 & 0 & \dots & 1\rule{4mm}{0mm}\vrule  & 0 & 0 & \dots & 0 & 0\\
  \hline
  0 & 0 & \dots & 0\rule{4mm}{0mm}\vrule  & 0 & 1 & \dots & 0 & 0\\
  0 & 0 & \dots & 0\rule{4mm}{0mm}\vrule  & 1 & 0 & \dots & 0 & 0\\
  \vdots & \vdots & \ddots & \vdots\rule{5mm}{0mm}\vrule  & \vdots & \vdots & \ddots
& \vdots & \vdots\\
  0 & 0 & \dots & 0\rule{4mm}{0mm}\vrule  & 0 & 0 & \dots & 1 & 0
 \end{pmatrix}
$$
by deleting the corresponding column. Notice that the last $l$
rows of $A''$ constitute the matrix $(0\ A_l)$ and, by virtue
of~(\ref{eq1sectHatLFred5}), are linearly dependent. Therefore,
after deleting any column from $A''$, we again obtain a degenerate
matrix. Hence, $\det \hat A=0$ and the rank of the matrix $A$ is
less than $2l+3$. Thus, Condition~\ref{condAllDxiPj} does fail,
and Theorem~\ref{thL_BFred} implies that the operator $\mathbf L_B
:W_B^{l+2}(G)\to W^l(G)$ is not Fredholm.

So, we proved that \textit{the operator $\mathbf L_B
:W_B^{l+2}(G)\to W^l(G)$ with odd $l\ge3$ is not Fredholm for any
$b_1,b_2$.}

\medskip

We considered the case where $b_1^2+b_2^2\ne0$. If $b_1=b_2=0$,
one can similarly show that the corresponding operator $\mathbf
L_B: W_B^{l+2}(G)\to W^l(G)$ is Fredholm for any $l\ge0$. However,
we omit the proof of this fact since, for $b_1=b_2=0$, we have the
``local'' Dirichlet problem in a smooth domain. As is well known,
such a problem is not only Fredholm but uniquely solvable for any
$l\ge0$.

\smallskip

The following theorem summarize the results obtained.

\begin{theorem}\label{th2sectHatLFred}
Let $l$ be even; then the operator $\mathbf L_B: W_B^{l+2}(G)\to
W^l(G)$ is Fredholm if and only if $b_1b_2>0$ or $b_1=b_2=0$.

Let $l$ be odd; then the operator $\mathbf L_B: W_B^{l+2}(G)\to
W^l(G)$ is Fredholm if and only if $b_1=b_2=0$.
\end{theorem}

\medskip

The author is grateful to A.~L.~Skubachevskii for attention to
this work.


\begin{thebibliography}{99}
\bibitem{Feller}
W.~Feller, ``Diffusion processes in one dimension,'' {\it Trans.
Amer. Math. Soc., \bf 77}, 1--30 (1954).
\bibitem{Sam}
A.~A.~Samarskii, ``On some problems of theory of differential
equations,'' {\it Differentsial'nye Uravneniya, \bf 16}, No.~11,
1925--1935 (1980); English transl.: {\it Diff. Equ., \bf 16}
(1980).
\bibitem{OnSk}
G.~G.~Onanov, A.~L.~Skubachevskii, ``Differential equations with
displaced arguments in stationary problems in the mechanics of a
deformed body,'' {\it Prikladnaya Mekhanika, \bf 15}, 39--47
(1979); English transl.: {\it Soviet Applied Mech., \bf 15}
(1979).
\bibitem{Sommerfeld}
A.~Sommerfeld, ``Ein Beitrag zur hydrodinamischen Erkl\"arung der
turbulenten Flussigkeitsbewegungen,'' {\it Proc. Intern. Congr.
Math. (Rome, 1908), Reale Accad. Lincei. Roma. \bf 3}, 116--124
(1909).
\bibitem{Tamarkin}
J.~D.~Tamarkin,  {\it Some General Problems of the Theory of
Ordinary Linear Differential Equations and Expansion of an
Arbitrary Function in Series of Fundamental Functions}, Petrograd,
1917. Abridged English transl.: {\it Math. Z., \bf 27}, 1--54
(1928).
\bibitem{Picone}
M.~Picone, ``Equazione integrale traducente il pi\`u generale
problema lineare per le equazioni differenziali lineari ordinarie
di qualsivoglia ordine,'' {\it Academia nazionale dei Lincei. Atti
dei convegni., \bf 15}, 942--948 (1932).
\bibitem{Carleman}
T.~Carleman, ``Sur la th\'eorie des equations integrales et ses
applications,'' {\it Verhandlungen des Internat. Math. Kongr.
Z\"urich., \bf 1}, 132--151 (1932).
\bibitem{BitzSam}
A.~V.~Bitsadze, A.~A.~Samarskii, ``On some simple generalizations
of linear elliptic boundary value problems,'' {\it Dokl. Akad.
Nauk SSSR., \bf 185}, 739--740 (1969); English transl.: {\it
Soviet Math. Dokl., \bf 10} (1969).
\bibitem{ZhEid}
S.~D.~Eidelman, N.~V.~Zhitarashu, ``Nonlocal boundary value
problems for elliptic equations,'' {\it Mat. Issled., \bf 6},
No.~2~(20), 63--73 (1971) [in Russian].
\bibitem{RSh}
Ya.~A.~Roitberg, Z.~G.~Sheftel', ``Nonlocal problems for elliptic
equations and systems,'' {\it Sib. Mat. Zh., \bf 13} (1972),
165--181; English transl.: {\it Siberian Math. J., \bf 13} (1972).
\bibitem{Kishk}
K.~Yu.~Kishkis, ``The index of a Bitsadze--Samarskii Problem for
harmonic functions,'' {\it Differentsial'nye Uravneniya., \bf 24},
No.~1, 105--110 (1988); English transl.: {\it Diff. Equ., \bf 24}
(1988).
\bibitem{GM}
A.~K.~Gushchin, V.~P.~Mikhailov, ``On solvability of nonlocal
problems for elliptic equations of second order,'' {\it Mat. sb.,
\bf 185}, 121--160 (1994); English transl.: {\it Math. Sb., 185}
(1994).
\bibitem{SkMs83}
A.~L.~Skubachevskii, ``Nonlocal elliptic problems with a
parameter,'' {\it Mat. Sb., \bf 121~\rm(\bf 163\rm)}, 201--210
(1983); English transl.: {\it Math. USSR Sb., \bf 49} (1984).
\bibitem{SkMs86}
A.~L.~Skubachevskii, ``Elliptic problems with nonlocal conditions
near the boundary,'' {\it Mat. Sb., \bf 129~\rm(\bf 171\rm)},
279--302 (1986); English transl.: {\it Math. USSR Sb., \bf 57}
(1987).
\bibitem{SkDu90}
A.~L.~Skubachevskii, ``Model nonlocal problems for elliptic
equations in dihedral angles,'' {\it Differentsial'nye Uravneniya,
\bf 26}, 119--131 (1990); English transl.: {\it Differ. Equ., \bf
26} (1990).

\bibitem{SkDu91}
A.~L.~Skubachevskii, ``Truncation-function method in the theory of
nonlocal problems,'' {\it Differentsial'nye Uravneniya, \bf 27},
128--139 (1991); English transl.: {\it Diff. Equ., \bf 27} (1991).
\bibitem{SkBook}
A.~L.~Skubachevskii, {\it Elliptic Functional Differential
Equations and Applications}, Basel--Boston--Berlin, Birkh\"auser,
1997.
\bibitem{KovSk}
O.~A.~Kovaleva, A.~L.~Skubachevskii, ``Solvability of nonlocal
elliptic problems in weighted spaces,'' {\it Mat. Zametki., \bf
67}, 882--898 (2000); English transl.: {\it Math. Notes., \bf 67}
(2000).
\bibitem{SkRJMP}
A.~L.~Skubachevskii, ``Regularity of solutions for some nonlocal
elliptic problem,'' {\it Russ. J. Math. Phys., \bf 8}, 365--374
(2001).
\bibitem{GurGiessen}
P.~L.~Gurevich, ``Nonlocal problems for elliptic equations in
dihedral angles and the Green formula,'' {\it Mitteilungen aus dem
Mathem. Seminar Giessen, Math. Inst. Univ. Giessen, Germany, \bf
247}, 1--74 (2001).
\bibitem{KondrTMMO67}
V.~A.~Kondrat'ev, ``Boundary value problems for elliptic equations
in domains with conical or angular points,'' {\it Trudy Moskov.
Mat. Obshch., \bf 16}, 209--292 (1967); English transl.: {\it
Trans. Moscow Math. Soc., \bf 16} (1967).
\bibitem{LM}
J.~L.~Lions, E~Magenes, {\it Non-Homogeneous Boundary Value
Problems and Applications, Vol.~$1$}, Springer--Verlag, New
York--Heidelberg--Berlin, 1972.
\bibitem{GS}
I.~C.~Gohberg, E.~I.~Sigal, ``An operator generalization of the
logarithmic residue theorem and the theorem of Rouch\'e,'' {\it
Mat. Sb., \bf 84~\rm(\bf 126\rm)}, 607--629 (1971); English
transl.: {\it Math. USSR Sb., \bf 13} (1971).

\bibitem{RS}
F.~Riesz, Sz.-Nagy, {\it Le\c{c}ons d'Analyse Fonctionnelle},
Deuxi\`eme  \'edition, Budapest, 1953.
\bibitem{Stein}
E.~M.~Stein {\it Singular Integrals and Differentiability
Properties of Functions}, Princeton Univ. Press, Princeton, 1970.

\bibitem{AV}
M.~S.~Agranovich, M.~I.~Vishik, ``Elliptic problems with a
parameter and parabolic problems of general type,'' {\it Uspekhi
Mat. Nauk, \bf 19}, 53--161 (1964); English transl.: {\it Russian
Math. Surveys, \bf 19} (1964).

\bibitem{Volevich}
L.~R.~Volevich ``Solvavility of boundary-value problems for
general elliptic systems,'' {\it Mat. sb., \bf 68}, No.~3,
373--416 (1965).
\bibitem{Kr}
S.~G.~Krein {\it Linear Equations in Banach Spaces}, Nauka,
Moscow, 1971 [in Russian]; English transl.: Birkh\"auser, Boston,
1982.

\bibitem{GurAsympAngle}
P.~L.~Gurevich, ``Asymptotics of solutions for nonlocal elliptic
problems in plane angles,'' {\it Trudy seminara imeni
I.~G.~Petrovskogo, \bf 23}, 2003; English transl. {\it J. Math.
Sci., New York.}

\bibitem{MP}
V.~G.~Maz'ya, B.~A.~Plamenevskii, ``$L_p$-estimates of solutions
of elliptic boundary value problems in domains with edges,'' {\it
Trudy Moskov. Mat. Obshch., \bf 37}, 49--93 (1978). English
transl.: {\it Trans. Moscow Math. Soc., \bf 37} (1980).

\end{thebibliography}
\end{document}